# An Extreme-Value Approach for Testing the Equality of Large U-Statistic Based Correlation Matrices

CHENG ZHOU[1,*] FANG HAN[2,**] XIN-SHENG ZHANG[1,†] HAN LIU[3,‡]

[1] *Department of Statistics, Management School, Fudan University, Shanghai, China. E-mail:* [*] chengzhmike@gmail.com; [†] xszhang@fudan.edu.cn

[2] *Department of Statistics, University of Washington, Seattle, WA 98195, USA. E-mail:* [**] fanghan@uw.edu

[3] *Department of Operations Research and Financial Engineering, Princeton University, Princeton, NJ 08544, USA. E-mail:* [‡] hanliu@princeton.edu

There has been an increasing interest in testing the equality of large Pearson's correlation matrices. However, in many applications it is more important to test the equality of large rank-based correlation matrices since they are more robust to outliers and nonlinearity. Unlike the Pearson's case, testing the equality of large rank-based statistics has not been well explored and requires us to develop new methods and theory. In this paper, we provide a framework for testing the equality of two large U-statistic based correlation matrices, which include the rank-based correlation matrices as special cases. Our approach exploits extreme value statistics and the Jackknife estimator for uncertainty assessment and is valid under a fully nonparametric model. Theoretically, we develop a theory for testing the equality of U-statistic based correlation matrices. We then apply this theory to study the problem of testing large Kendall's tau correlation matrices and demonstrate its optimality. For proving this optimality, a novel construction of least favorable distributions is developed for the correlation matrix comparison.

## 1. Introduction

Let $\boldsymbol{X} = (X_1, \ldots, X_d)^T$ and $\boldsymbol{Y} = (Y_1, \ldots, Y_d)^T$ be two $d$-dimensional random vectors. We denote $\boldsymbol{X}_1, \ldots, \boldsymbol{X}_{n_1}$ with $\boldsymbol{X}_k = (X_{k1}, \ldots, X_{kd})^T$ to be $n_1$ independent samples of $\boldsymbol{X}$ and $\boldsymbol{Y}_1, \ldots, \boldsymbol{Y}_{n_2}$ with $\boldsymbol{Y}_k = (Y_{k1}, \ldots, Y_{kd})^T$ to be $n_2$ independent samples of $\boldsymbol{Y}$. Letting $n := \max\{n_1, n_2\}$, we aim to test the equality of U-statistic based correlation matrices (e.g, Kendall's tau or Spearman's rho) of $\boldsymbol{X}$ and $\boldsymbol{Y}$. We consider the high dimensional regime that $d, n \to \infty$ and $d/n$ does not necessarily go to zero as $n \to \infty$. This problem has important applications, including portfolio selection (Markowitz, 1991), high dimensional discriminant analysis (Han et al., 2013; Mai and Zou, 2013) and gene selection (Ho et al., 2008; Hu et al., 2009, 2010).

When $d/n \to 0$, Anderson (2003) and Muirhead (2009) study the problem of testing the equality of two Pearson's correlation matrices. Major test criteria include the likelihood ratio (Anderson, 2003), spectral norm of difference (Roy, 1957) and Frobenius norm of difference (Nagao, 1973). When $d/n \not\to 0$, the likelihood ratio test and the tests in Roy (1957) and Nagao (1973) perform poorly, as Pearson's sample correlation matrices no longer converge to their population counterparts under the spectral norm (Bai and Yin, 1993). A line of research aims to correct the aforementioned tests or proposing new methods. For the likelihood ratio test, Bai et al. (2009) introduce a corrected LRT test which works when $d/n \to c \in (0, 1)$, and Jiang et al. (2012) generalize it to the case when $d < n$ and $c = 1$. Based on the spectral norm of difference, Han et al. (2018) use the bootstrap method to generalize Roy's test in high dimension. As a generalization of Nagao's proposal, Schott (2007) and Li and Chen (2012) propose new test statistics based on an unbiased estimator of the Frobenius norm of the matrix difference, and Srivastava and Yanagihara (2010) propose another test statistic based on the difference of two Frobenius norms. Recently, Cai et al. (2013) propose a method based on the sup-norm of the matrix difference and prove its rate optimality under a sparse alternative.

In many applications, it is more meaningful to test the equality of two rank-based correlation matrices but instead of the Pearson's correlation matrices. In particular, Embrechts et al. (2003) point out that the Pearson's correlation coefficient "might prove very misleading" in measuring the dependence and advocate the usage of rank correlation coefficients, such as Kendall's tau (Kendall, 1938) or Spearman's rho (Spearman, 1904). Though testing the equality of high dimensional rank-based correlation matrices is of fundamental importance, there has been very little work in this area. To bridge this gap, this paper proposes a unified framework for testing the equality of two large U-statistic based correlation matrices $\mathbf{U}_1$ and $\mathbf{U}_2$, which include rank-based correlation matrices as special examples. More specifically, let $\mathbf{U}_1 = (u_{1,ij})$ be a type of correlation matrix of $\boldsymbol{X}$ and all the elements of $\mathbf{U}_1$ can be estimated by U-statistics[1]. Similarly to $\mathbf{U}_1$, we define $\mathbf{U}_2 = (u_{2,ij})$ to be the same kind of U-statistic based correlation matrix of $\boldsymbol{Y}$. In this paper, we aim to test the hypothesis

$$\mathbf{H_0} : \mathbf{U}_1 = \mathbf{U}_2 \quad \text{v.s.} \quad \mathbf{H_1} : \mathbf{U}_1 \neq \mathbf{U}_2. \tag{1.1}$$

Testing (1.1) plays an important role in many fields. For example, testing the equality of two Kendall's tau correlation matrices $\mathbf{U}_1^\tau$ and $\mathbf{U}_2^\tau$,

$$\mathbf{H_0^\tau} : \mathbf{U}_1^\tau = \mathbf{U}_2^\tau \quad \text{v.s.} \quad \mathbf{H_1^\tau} : \mathbf{U}_1^\tau \neq \mathbf{U}_2^\tau, \tag{1.2}$$

can be used to test the model of copula discriminant analysis (Han et al., 2013; Mai and Zou, 2013).

There are 4 major contributions of this paper. First, for the first time in the literature, we develop a unified framework for testing the equality of two large U-statistic based correlation matrices. This framework builds upon a fully nonparametric model and enables us to conduct homogeneity tests using a wide range of correlation measures. Secondly, as a special example, we examine the problem of testing the equality of two large Kendall's tau matrices and prove the minimax optimality of the proposed method. Thirdly, we further propose alternative approaches for testing $\mathbf{U}_1^\tau = \mathbf{U}_2^\tau$, which attain better empirical performance than the Jackknife based one. Finally, to develop a theory of testing the equality of general U-statistic based correlation matrices, we develop an upper bound of Jackknife variance estimation error, which enables us to obtain the explicit rate of convergence. For Kendall's tau matrices, we prove an upper bound of the traditional plug-in variance estimation error and an upper bound of the variance difference between two Kendall's tau correlation coefficients. These upper bounds allow us to exploit the extreme value theory under the dependent setting to prove theorems in this paper. Their constructions are nontrivial and are of independent technical interest. To prove the optimality of the proposed testing methods for Kendall's tau matrices, we construct a collection of least favorable distributions with regard to the test hypothesis. This construction technique is novel and tailored for testing the equality of correlation matrices. In contrast, the construction in Cai et al. (2013) only perturbs the diagonal elements of covariance matrices, which does not affect the resulting correlation matrices.

[1]Such U-statistic based correlation measures are quite general. For example, $u_{1,ij}$ can represent the Kendall's tau correlation coefficient between $X_i$ and $X_j$.

## 1.1. More Related Works

Apart from the Pearson's correlation coefficient and general U-statistic based correlation measurements studied in this paper, existing literature also considers other measures of dependence. These include the distance correlation (Székely et al., 2007) and randomized dependence coefficient (Lopez-Paz et al., 2013). To the best of our knowledge, there is no work discussing testing the equality of dependence structure with regard to these dependent measures.

Our work is closely related to the random matrix theory on rank correlation matrices. Bai and Zhou (2008), Zhou (2007), Bao et al. (2013), and Han et al. (2017) study the theoretical properties of large rank-based correlation matrices. Specifically, for these random matrices, Bai and Zhou (2008) prove the Marchenko-Pastur law for the limiting spectral distribution, Zhou (2007) and Han et al. (2017) prove the extreme value type I distribution for the entry-wise maximum, and Bao et al. (2013) derive the limiting distributions of traces of all higher moments. Most of these results hold only under the independence setting, i.e., the entries of $\boldsymbol{X}$ are independent of each other. In contrast, our work focuses on the dependent setting.

Our work is also related to the robust testing, where the test statistics are robust estimators of the Pearson's covariance/correlation coefficients. These include S-estimators and some robust dispersion estimators. We refer to O'Brien (1992), Aslam and Rocke (2005) and the references therein for details. Our work is also related to the adaptive estimation of a large correlation/covariance matrix (Cai and Liu, 2011) or a large Gaussian (copula) graphical model. See, for example, Bickel and Levina (2008), Zhao et al. (2014), Ravikumar et al. (2011) and Liu et al. (2012).

### 1.2. Notation

We denote $\|\mathbf{v}\|_2 = \big(\sum_{j=1}^d v_j^2\big)^{1/2}$ as the Euclidean norm of a vector $\mathbf{v} = (v_1, \ldots, v_d)^T \in \mathbb{R}^d$. For a matrix $\mathbf{A} = (a_{ij}) \in \mathbb{R}^{d\times q}$, we define its spectral norm $\|\mathbf{A}\|_2 := \sup_{\|\mathbf{x}\|_2 \le 1} \|\mathbf{A}\mathbf{x}\|_2$ and Frobenius norm $\|\mathbf{A}\|_F := \sqrt{\sum_{i,j} a_{ij}^2}$. We define the matrix entrywise sup-norm as $\|\mathbf{A}\|_{\max} := \max\big\{|a_{ij}|\big\}$. We use $\mathrm{Rank}(\mathbf{A})$ to denote the rank of $\mathbf{A}$. If $\mathbf{A}$ is a square matrix, we define $\mathrm{Diag}(\mathbf{A})$ to be a diagonal matrix with the same main diagonal as $\mathbf{A}$. We use $\mathbf{I}_d$ to denote an identity matrix of size $d$. For two sequences of real numbers $\{a_n\}$ and $\{b_n\}$, we write $a_n = O(b_n)$ if there exists a constant $C$ such that $|a_n| \le C|b_n|$ holds for all sufficiently large $n$, write $a_n = o(b_n)$ if $a_n/b_n \to 0$, and write $a_n \asymp b_n$ if there exist constants $C \ge c > 0$ such that $c|b_n| \le |a_n| \le C|b_n|$ for all sufficiently large $n$. For a square matrix $\mathbf{\Sigma} \in \mathbb{R}^{d\times d}$, we use $\lambda_{\min}(\mathbf{\Sigma})$ and $\lambda_{\max}(\mathbf{\Sigma})$ to denote the minimal and maximal eigenvalues of $\mathbf{\Sigma}$. For a set $B$, we use $|B|$ to denote its cardinality.

### 1.3. Paper Organization

The rest of this paper is organized as follows. Section 2 formalizes the problem, describes a general testing procedure and analyzes the theoretical properties (e.g., size and power) of the proposed test. In Section 3 we focus on testing large Kendall's tau matrices, for which we consider two models: a fully nonparametric model and a semiparametric Gaussian copula model. Under certain modelling assumptions, for Kendall's tau matrices we propose additional tests which have better empirical performance compared to the general testing procedure. Section 4 provides thorough numerical results on both simulated and real data. In Section 5, we discuss potential future work. Appendix A contains the proof of the main theorem. We put the proofs of all other results in Supplementary Material of this paper.

## 2. A General Procedure for Testing U-Statistic Based Matrices

This section presents a generic testing method for U-statistic based matrix comparison. In Section 2.1 we describe the proposed testing procedure. In Section 2.2 we analyze its asymptotic size and power. In Section 2.3 we consider comparing a row or column of U-statistic based matrices.

Before presenting the testing procedure, we introduce some notations for U-statistics. For $i, j = 1, \ldots, q$, let $\Phi_{ij}$ be a U-statistic's kernel function defined as

$$\Phi_{ij} : \underbrace{\mathbb{R}^d \times \cdots \times \mathbb{R}^d}_{m} \to \mathbb{R} \ \text{ with the symmetric property} : \Phi_{ij} = \Phi_{ji}, \tag{2.1}$$

where $m$ is the kernel order. Thus, we have a family of functions $\big\{\Phi_{ij}, 1 \le i, j \le q\big\}$. Furthermore, each $\Phi_{ij}$ is a symmetric Borel measurable function with the kernel order $m$ fixed[2]. We assume that $\Phi_{ij}$ is uniformly bounded. Many useful U-statistics satisfy these conditions. We set

$$\widehat{u}_{1,ij} := \binom{n_1}{m}^{-1} \sum_{1\le \ell_1 < \cdots < \ell_m \le n_1} \Phi_{ij}(\boldsymbol{X}_{\ell_1}, \cdots, \boldsymbol{X}_{\ell_m}),$$
$$\widehat{u}_{2,ij} := \binom{n_2}{m}^{-1} \sum_{1\le \ell_1 < \cdots < \ell_m \le n_2} \Phi_{ij}(\boldsymbol{Y}_{\ell_1}, \cdots, \boldsymbol{Y}_{\ell_m}).$$

We then define the following U-statistic based matrices $\widehat{\mathbf{U}}_a \in \mathbb{R}^{q\times q}$ for $a = 1, 2$:

$$\widehat{\mathbf{U}}_1 := \big(\widehat{u}_{1,ij}\big)_{1\le i,j\le q} \quad \text{and} \quad \widehat{\mathbf{U}}_2 := \big(\widehat{u}_{2,ij}\big)_{1\le i,j\le q}. \tag{2.2}$$

Correspondingly, we use $\mathbf{U}_a := (u_{a,ij})_{1\le i,j\le q}$ to denote the expectation of $\widehat{\mathbf{U}}_a$, i.e., $u_{a,ij} = \mathbb{E}[\widehat{u}_{a,ij}]$. We can view $\mathbf{U}_1$ and $\mathbf{U}_2$ as a type of correlation matrices of $\boldsymbol{X}$ and $\boldsymbol{Y}$. We are interested in testing the equality of $\mathbf{U}_1$ and $\mathbf{U}_2$, which includes testing the equality of two large Kendall's tau or Spearman's rho correlation matrices.

We note that $q$ is the row and column number of $\mathbf{U}_a$ and $\widehat{\mathbf{U}}_a$, while $d$ is the dimension of $\boldsymbol{X}$ and $\boldsymbol{Y}$. $q$ and $d$ can be different. Therefore, the framework considered in this paper is quite general. For example, it allows $\mathbf{U}_a$ to represent the dependence structure on a dimension reduced data, where the dimension reduction step is incorporated in the kernel function $\{\Phi_{ij}, 1 \le i, j \le q\}$.

**Remark 2.1.** We can relax $\Phi_{ij}$ to be an asymmetric kernel function without loss of generality. Specifically, an asymmetric kernel $\Phi(\cdot)$ gives a U-statistic

$$\widehat{u} = \frac{1}{m!}\binom{n_1}{m}^{-1} \sum \Phi(\boldsymbol{X}_{\ell_1}, \cdots, \boldsymbol{X}_{\ell_m}),$$

where the summation is taken over all combinations of distinct elements $\{\ell_1, \ldots, \ell_m\}$ from $\{1, \ldots, n_1\}$. Using the Hoeffding's method (Hoeffding, 1948), $\widehat{u}$ is also a U-statistic of the symmetric kernel $\Phi^0(\cdot)$:

$$\Phi^0(\mathbf{x}_1, \mathbf{x}_2, \ldots, \mathbf{x}_m) = \frac{1}{m!}\sum \Phi(\mathbf{x}_{\alpha_1}, \cdots, \mathbf{x}_{\alpha_m}),$$

where the summation is taken over all permutations of $\{1, \ldots, m\}$. For example, to construct an unbiased[3] estimator for Spearman's rho, El Maache and Lepage (2003) recommends to use the U-statistic with the kernel

$$\Phi_{ij}(\boldsymbol{X}_1, \boldsymbol{X}_2, \boldsymbol{X}_3) = 2^{-1}\sum\sum\sum_{1\le\alpha\ne\beta\ne\gamma\le 3} \operatorname{sign}(X_{\alpha i} - X_{\beta i}) \operatorname{sign}(X_{\alpha j} - X_{\gamma j}).$$

[2] We assume each $\Phi_{ij}$ has the same fixed kernel order $m$ for presentation clearness. It is straightforward to extend to the setting that $m$'s are uniformly bounded.

[3] The Spearman rank-order correlation, i.e., the sample correlation between the rank values of two variables, is a biased estimator of the population Spearman's rho.

The Kendall's tau matrix is an example of the U-statistic based matrix defined in (2.2). More specifically, we set

$$\Phi_{ij}(\boldsymbol{X}_k, \boldsymbol{X}_\ell) = \mathrm{sign}(X_{ki} - X_{\ell i})\,\mathrm{sign}(X_{kj} - X_{\ell j}),$$
$$\Phi_{ij}(\boldsymbol{Y}_k, \boldsymbol{Y}_\ell) = \mathrm{sign}(Y_{ki} - Y_{\ell i})\,\mathrm{sign}(Y_{kj} - Y_{\ell j}),$$

and $q = d$. The Kendall's tau sample correlation coefficients $\widehat{\tau}_{1,ij}$ and $\widehat{\tau}_{2,ij}$ are then defined as

$$\widehat{\tau}_{1,ij} := \frac{2}{n_1(n_1-1)} \sum_{1\le k<\ell\le n_1} \mathrm{sign}(X_{ki} - X_{\ell i})\,\mathrm{sign}(X_{kj} - X_{\ell j}),$$
$$\widehat{\tau}_{2,ij} := \frac{2}{n_2(n_2-1)} \sum_{1\le k<\ell\le n_2} \mathrm{sign}(Y_{ki} - Y_{\ell i})\,\mathrm{sign}(Y_{kj} - Y_{\ell j}).$$

Their population counterparts are $\tau_{a,ij} := \mathbb{E}[\widehat{\tau}_{a,ij}]$ for $a = 1, 2$. We then write sample and population Kendall's tau matrices as

$$\widehat{\mathbf{U}}_a^\tau = (\widehat{\tau}_{a,ij}) \qquad \text{and} \qquad \mathbf{U}_a^\tau = (\tau_{a,ij}), \tag{2.3}$$

where $a = 1, 2$. In Section 3, we consider testing the large Kendall's tau matrices.

## 2.1. A General Testing Procedure

For testing (1.1) in high dimensions, we use the sup-norm criterion. Such a choice is motivated by the fact that the sup-norm is very sensitive to perturbations on a small number of entries compared to the null hypothesis. We then propose the test statistic:

$$M_n := \max_{1\le i,j\le q} M_{ij} \quad \text{with} \quad M_{ij} := \frac{(\widehat{u}_{1,ij} - \widehat{u}_{2,ij})^2}{\widehat{\sigma}^2(\widehat{u}_{1,ij}) + \widehat{\sigma}^2(\widehat{u}_{2,ij})} \quad \text{for} \quad 1 \le i,j \le q. \tag{2.4}$$

In (2.4), $\widehat{\sigma}^2(\widehat{u}_{1,ij})$ is a Jackknife estimator of $\widehat{u}_{1,ij}$'s variance and is defined as

$$\widehat{\sigma}^2(\widehat{u}_{1,ij}) := \frac{m^2(n_1-1)}{n_1(n_1-m)^2} \sum_{\alpha=1}^{n_1} (q_{1\alpha,ij} - \widehat{u}_{1,ij})^2, \tag{2.5}$$

with

$$q_{1\alpha,ij} := \binom{n_1-1}{m-1}^{-1} \sum_{\substack{1\le \ell_1<\cdots<\ell_{m-1}\le n_1 \\ \ell_j\ne\alpha, j=1,\cdots,m-1}} \Phi_{ij}(\boldsymbol{X}_\alpha, \boldsymbol{X}_{\ell_1}, \ldots, \boldsymbol{X}_{\ell_{m-1}}).$$

The definition of $\widehat{\sigma}^2(\widehat{u}_{2,ij})$ is similar for $\boldsymbol{Y}$.

For a given significance level $0 < \alpha < 1$, we construct the test to be

$$\mathrm{T}_\alpha := \mathbb{1}\left\{M_n \ge G^-(\alpha) + 4\log q - \log(\log q)\right\}, \tag{2.6}$$

where $G^{-}(\alpha) := -\log(8\pi) - 2\log\big(-\log(1-\alpha)\big)$. We reject $\mathbf{H_0}$ in (1.1) if and only if $\mathrm{T}_\alpha = 1$.

In some applications, our interest is to compare a particular row or column of matrices, i.e., we aim at testing the hypothesis:

$$\mathbf{H_{0,i}} : \mathbf{u}_{1,i\star} = \mathbf{u}_{2,i\star} \quad \text{v.s.} \quad \mathbf{H_{1,i}} : \mathbf{u}_{1,i\star} \neq \mathbf{u}_{2,i\star}, \tag{2.7}$$

where $\mathbf{u}_{1,i\star}$ and $\mathbf{u}_{2,i\star}$ are the $i$-th rows of $\mathbf{U}_1$ and $\mathbf{U}_2$. To test this hypothesis, we construct a similar test statistic $M_{n,i} = \max_{1\le j\le q} M_{i,j}$, and the according test is

$$\mathrm{T}_{\alpha,i} = \mathbb{1}\big\{M_{n,i} > G'^{-}(\alpha) + 2\log q - \log\log q\big\}, \tag{2.8}$$

where $G'^{-}(\alpha) := -\log(\pi) - 2\log(-\log(1-\alpha))$. We reject $\mathbf{H_{0,i}}$ if and only if $\mathrm{T}_{\alpha,i} = 1$.

## 2.2. Theoretical Properties

Our main theoretical result is to characterize the limiting null distribution of $M_n$. We further analyze the power of the proposed test under a sparse alternative.

We introduce three assumptions that will be used later. Assumption **(A1)** specifies the sparsity of $\mathbf{U} = \mathbf{U}_1 = \mathbf{U}_2$. Assumption **(A2)** specifies the scaling of $q, n$. Assumption **(A3)** is a technical condition that we impose for obtaining the limiting distribution of $M_n$. In Section 3.2, we will show that Assumption **(A3)** can be further relaxed under a semiparametric Gaussian copula model.

In detail, for a fixed constant $\alpha_0 > 0$, we define

$$\mathrm{supp}_j(\alpha_0) := \big\{1 \le i \le q : |u_{1,ij}| \ge (\log q)^{-1-\alpha_0} \;\text{ or }\; |u_{2,ij}| \ge (\log q)^{-1-\alpha_0}\big\}.$$

$\mathrm{supp}_j(\alpha_0)$ is the set of indices $i$ such that either the $i$-th variable of $\boldsymbol{X}$ is highly correlated ($|u_{a,ij}| > (\log q)^{-1-\alpha_0}$) with the $j$-th variable of $\boldsymbol{X}$, or the $i$-th variable of $\boldsymbol{Y}$ is highly correlated with the $j$-th variable of $\boldsymbol{Y}$. We then introduce Assumption **(A1)** as follows.

> **(A1)**. We assume that there exits a subset $\Gamma \subset \{1, 2, \ldots, q\}$ with $|\Gamma| = o(q)$ and a constant $\alpha_0 > 0$ such that for all $\gamma > 0$, we have
>
> $$\max_{1\le j\le q, j\notin\Gamma} \big|\mathrm{supp}_j(\alpha_0)\big| = o(q^\gamma).$$

Before stating Assumption **(A2)**, we need some additional notations. Set

$$\Psi_{ij}(\boldsymbol{X}_{\ell_1}, \ldots, \boldsymbol{X}_{\ell_m}) := \Phi_{ij}(\boldsymbol{X}_{\ell_1}, \ldots, \boldsymbol{X}_{\ell_m}) - u_{1,ij}.$$

For $\ell = 1, \ldots, n_1$, we also denote $g_{ij}(\boldsymbol{X}_\ell)$ and $h_{ij}(\boldsymbol{X}_\ell)$ as

$$g_{ij}(\boldsymbol{X}_\ell) := \mathbb{E}[\Phi_{ij}(\boldsymbol{X}_{\ell_1}, \ldots, \boldsymbol{X}_{\ell_m})|\boldsymbol{X}_\ell], \qquad h_{ij}(\boldsymbol{X}_\ell) := \mathbb{E}[\Psi_{ij}(\boldsymbol{X}_{\ell_1}, \ldots, \boldsymbol{X}_{\ell_m})|\boldsymbol{X}_\ell], \tag{2.9}$$

where $\{\ell_1,\dots,\ell_m\}$ is an arbitrary subset of $\{1,\dots,n_1\}$ with distinct elements and contains $\ell$. $g_{ij}(\boldsymbol{Y}_\ell)$ and $h_{ij}(\boldsymbol{Y}_\ell)$ are similarly defined for $\ell=1,\cdots,n_2$. We then denote $\zeta_{1,ij}$ to be the variance of $g_{ij}(\boldsymbol{X}_\ell)$, i.e.,

$$\zeta_{1,ij} := \mathbb{E}\Big[\mathbb{E}\big[\Psi_{ij}(\boldsymbol{X}_{\ell_1},\dots,\boldsymbol{X}_{\ell_m})|\boldsymbol{X}_\ell\big]^2\Big] = \mathrm{Var}\big(g_{ij}(\boldsymbol{X}_\ell)\big). \tag{2.10}$$

Similarly, we define $\zeta_{2,ij} := \mathrm{Var}(g_{ij}(\boldsymbol{Y}_\ell))$.

With these introduced notations, we are now ready to state Assumption **(A2)**.

> **(A2)**. We assume $n_1 \asymp n_2 \asymp n$ and $\log q = O(n^{1/3-\epsilon})$ for an arbitrary $0<\epsilon<1/3$. We also assume $\zeta_{a,ij} > r_a > 0$ for $a=1,2$, where $r_1$ and $r_2$ are constants which are irrelevant to $i$ and $j$.

The condition that $\zeta_{a,ij} > r_a > 0$ is mild. It is used to exclude the degenerate cases of U-statistics and has been widely used for analyzing U-statistics.

To describe Assumption **(A3)**, we write

$$S = \big\{(i,j):\ 1\le i,\ j\le q\big\} \quad \text{and} \quad S_0 = \big\{(i,j):\ 1\le i\le q,\ i\in \mathrm{supp}_j(\alpha_0)\big\}. \tag{2.11}$$

By the definition of $S_0$, for any $(i,j)\in S\setminus S_0$, we have $|u_{a,ij}| \le (\log q)^{-1-\alpha_0}$. Moreover, we use $u_{1,ijk\ell}$ and $u_{2,ijk\ell}$ to denote $\mathbb{E}\big[g_{ij}(\boldsymbol{X}_\ell)g_{k\ell}(\boldsymbol{X}_\ell)\big]$ and $\mathbb{E}\big[g_{ij}(\boldsymbol{Y}_\ell)g_{k\ell}(\boldsymbol{Y}_\ell)\big]$.

> **(A3)** Assume $u_{a,ijk\ell} = O((\log q)^{-1-\alpha_0})$ for any $(i,j)\neq(k,\ell)\in S\setminus S_0$ and $a=1,2$.

Under fully nonparametric models, we note that $u_{a,ijk\ell}$ is estimable (Klüppelberg and Kuhn, 2009). Thus it is possible to verify Assumption **(A3)** in applications. When we test the equality of two Kendall's tau correlation matrices $\mathbf{U}_1^\tau$ and $\mathbf{U}_2^\tau$, under a semiparametric Gaussian copula model Assumption **(A3)** can be replaced by a simplified condition which is easier to be verified. More details are provided in Section 3.2.

Under the above assumptions, our main theoretical result quantifies the limiting distribution of the extreme value statistic $M_n$.

**Theorem 2.2.** Assuming **(A1)**, **(A2)**, **(A3)** hold, under $\mathbf{H_0}$ of (1.1), we have

$$\mathbb{P}\big(M_n - 4\log q + \log(\log q) \le x\big) \to \exp\Big(-\frac{1}{\sqrt{8\pi}}\exp\big(-\frac{x}{2}\big)\Big), \tag{2.12}$$

for any $x\in\mathbb{R}$, as $n,q\to\infty$. Furthermore, (2.12) holds uniformly for all random vectors $\boldsymbol{X}$ and $\boldsymbol{Y}$ satisfying Assumptions **(A1)**, **(A2)** and **(A3)**.

***Proof.*** We list a sketch of this proof. The detailed proof is in Appendix A. The proof proceeds in three steps.

**Step (i) (Sketch).** We set $\widehat{\sigma}^2(\widehat{u}_{a,ij})$ as the Jackknife variance estimator of $\widehat{u}_{a,ij}$ and $\sigma^2(\widehat{u}_{a,ij})$ as the true variance of $\widehat{u}_{a,ij}$. We then analyze the estimation error of Jackknife variance estimator by providing an upper bound of $|n_a\widehat{\sigma}^2(\widehat{u}_{a,ij}) - m^2\zeta_{a,ij}|$, where $\zeta_{a,ij}$ is defined in (2.10). The central limit theorem for U-statistics (Lemma D.3 in Supplement

D of Supplementary Material) implies that $m^2\zeta_{a,ij}$ is the limit of $n_a\sigma^2(\widehat{u}_{a,ij})$ as $n_a$ goes to infinity. This motivates us to define

$$M_{ij} := \frac{(\widehat{u}_{1,ij} - \widehat{u}_{2,ij})^2}{\widehat{\sigma}^2(\widehat{u}_{1,ij}) + \widehat{\sigma}^2(\widehat{u}_{2,ij})} \quad \text{and} \quad \widetilde{M}_{ij} := \frac{(\widehat{u}_{1,ij} - \widehat{u}_{2,ij})^2}{m^2\zeta_{1,ij}/n_1 + m^2\zeta_{2,ij}/n_2}. \tag{2.13}$$

In $\widetilde{M}_{ij}$, we use $m^2\zeta_{a,ij}/n_1$ to replace $\widehat{\sigma}^2(\widehat{u}_{a,ij})$ of $M_{ij}$.

By using the obtained upper bound of $|n_a\widehat{\sigma}^2(\widehat{u}_{a,ij}) - m^2\zeta_{a,ij}|$, we prove $\max_{1\le i,j\le q} M_{ij}$ and $\max_{1\le i,j\le q} \widetilde{M}_{ij}$ have the same limiting distribution, i.e., it suffices to prove that

$$\lim_{n,q\to\infty} \mathbb{P}\big(\widetilde{M}_n - 4\log q + \log(\log q) \le x\big) = \exp\big(-\exp(-x/2)/\sqrt{8\pi}\big), \tag{2.14}$$

where $\widetilde{M}_n := \max_{1\le i,j\le q} \widetilde{M}_{ij}$.

**Step (ii) (Sketch).** We use the Hoeffding decomposition (Lemma D.4 in Supplement D of Supplementary Material) to decompose the U-statistic $\widetilde{u}_{a,ij} := \widehat{u}_{a,ij} - u_{a,ij}$. By the definition of $\widetilde{u}_{a,ij}$, we have $\mathbb{E}[\widetilde{u}_{a,ij}] = 0$. By the Hoeffding decomposition, we decompose $\widetilde{u}_{a,ij}$ into two pieces. One is the sum of independent and identically distributed (i.i.d.) random variables and the other is the residual term. In detail, decompose $\widetilde{u}_{a,ij}$ as

$$\widetilde{u}_{1,ij} = \frac{m}{n_1}\sum_{\alpha=1}^{n_1} h_{ij}(\boldsymbol{X}_\alpha) + \binom{n_1}{m}^{-1}\Delta_{n_1,ij},\ \widetilde{u}_{2,ij} = \frac{m}{n_2}\sum_{\alpha=1}^{n_2} h_{ij}(\boldsymbol{Y}_\alpha) + \binom{n_2}{m}^{-1}\Delta_{n_2,ij}, \tag{2.15}$$

where we set

$$\Delta_{n_1,ij} = \sum_{1\le \ell_1<\ell_2<\ldots<\ell_m\le n_1} \Big(\Phi_{ij}(\boldsymbol{X}_{\ell_1},\ldots,\boldsymbol{X}_{\ell_m}) - u_{1,ij} - \sum_{k=1}^m h_{ij}(\boldsymbol{X}_{\ell_k})\Big),$$
$$\Delta_{n_2,ij} = \sum_{1\le \ell_1<\ell_2<\ldots<\ell_m\le n_2} \Big(\Phi_{ij}(\boldsymbol{Y}_{\ell_1},\ldots,\boldsymbol{Y}_{\ell_m}) - u_{2,ij} - \sum_{k=1}^m h_{ij}(\boldsymbol{Y}_{\ell_k})\Big).$$

Apparently, $m\sum_{\alpha=1}^{n_1} h_{ij}(\boldsymbol{X}_\alpha)/n_1$ and $m\sum_{\alpha=1}^{n_2} h_{ij}(\boldsymbol{Y}_\alpha)/n_2$ are terms for the sum of i.i.d. random variables and $\binom{n_a}{m}^{-1}\Delta_{n_a,ij}$ is the residual term. We then use $m\sum_{\alpha=1}^{n_1} h_{ij}(\boldsymbol{X}_\alpha)/n_1$ and $m\sum_{\alpha=1}^{n_2} h_{ij}(\boldsymbol{Y}_\alpha)/n_2$ as the approximations of $\widetilde{u}_{1,ij}$ and $\widetilde{u}_{2,ij}$ and define

$$T_{ij} := \frac{\sum_{\alpha=1}^{n_1} h_{ij}(\boldsymbol{X}_\alpha)/n_1 - \sum_{\alpha=1}^{n_2} h_{ij}(\boldsymbol{Y}_\alpha)/n_2}{\sqrt{\zeta_{1,ij}/n_1 + \zeta_{2,ij}/n_2}} \quad \text{and} \quad T_n := \max_{1\le i,j\le q} (T_{ij})^2. \tag{2.16}$$

We then prove that the small residual term $\binom{n_a}{m}^{-1}\Delta_{n_a,ij}$ is negligible for our theorem, i.e., to obtain Theorem 2.2, it suffices to prove that as $n$, $q\to\infty$, we have

$$\mathbb{P}\big(T_n - 4\log q + \log(\log q) \le x\big) \to \exp\big(-\exp(-x/2)/\sqrt{8\pi}\big). \tag{2.17}$$

**Step (iii) (Sketch).** In the last step, we derive the limiting distribution of $T_n$ to prove (2.17). $T_n$ is the maximum of $(T_{ij})^2$ over $\{1 \le i, j \le q\}$ and $T_{ij}$ is not independent of each other. Therefore, we cannot straightforwardly exploit the extreme value theorem under the independent setting to obtain the limiting distribution of $T_n$. To solve this problem, we exploit the normal approximation to get the extreme value distribution of $\{(T_{ij})^2\}_{1\leq i,j\leq q}$ under the setting that $T_{ij}$ can be dependent of each other. The detailed proof of this theorem is in Appendix A. □

Theorem 2.2 justifies the size of the proposed test $\mathrm{T}_\alpha$ in (2.6). It shows that under $\mathbf{H_0}$ of (1.1), $M_n - 4\log q + \log(\log q)$ converges weakly to an extreme value Type I distribution with the distribution function $F(t) = \exp\big(-\exp(t/2)/\sqrt{8\pi}\big)$.

**Remark 2.3.** Theorem 2.2 provides a unified framework for testing the equality of two large U-statistic based matrices, which include ranked-based correlation matrices as special examples. Our test method exploits the Jackknife strategy and extreme value statistics, and it works under a fully nonparametric model. Technically, for proving Theorem 2.2, we develop a set of tools for analyzing the Jackknife variance estimator defined in (2.5), which is technically nontrivial and is of independent interest for analyzing U-statistics in more general settings.

Next, we analyze the power of $\mathrm{T}_\alpha$. To this end, we first introduce an alternative hypothesis characterized by the following set of matrix pairs

$$\mathbb{A}(C) = \left\{(\mathbf{U}_1, \mathbf{U}_2) : \max_{1\leq i,j\leq q} \frac{|u_{1,ij} - u_{2,ij}|}{\sqrt{m^2\zeta_{1,ij}/n_1 + m^2\zeta_{2,ij}/n_2}} \geq C\sqrt{\log q}\right\},$$

where $C > 0$ is a constant. The setting that only one entry of $\mathbf{U}_1$ and $\mathbf{U}_2$ differentiates large enough will make $(\mathbf{U}_1, \mathbf{U}_2) \in \mathbb{A}(C)$ for some constant $C$. The next theorem shows that the null hypothesis is asymptotically distinguishable from $\mathbb{A}(4)$ by $\mathrm{T}_\alpha$, i.e., we can use $\mathrm{T}_\alpha$ to reject $\mathbf{H_0}$ in (1.1) with an overwhelming probability if $(\mathbf{U}_1, \mathbf{U}_2) \in \mathbb{A}(4)$.

**Theorem 2.4.** (Power of the Test $\mathrm{T}_\alpha$) If **(A2)** is satisfied, as $n, q \to \infty$ we have

$$\inf_{(\mathbf{U}_1,\mathbf{U}_2)\in\mathbb{A}(4)} \mathbb{P}(\mathrm{T}_\alpha = 1) \to 1. \tag{2.18}$$

**Remark 2.5.** From the above theorem, for big enough $C > 0$, only one entry of $\mathbf{U}_1 - \mathbf{U}_2$ has a magnitude more than $C\sqrt{\log q/n}$ is enough for the test $\mathrm{T}_\alpha$ to correctly reject $\mathbf{H_0}$ of (1.1). We don't impose Assumptions (**A1**) and (**A3**) to obtain such results.

## 2.3. Testing Rows or Columns of Two U-statistic Based Matrices

In some applications, instead of testing the equality of two full matrices, we are interested in testing the equality of a particular row or column of the given matrix pair. This requires us to test the hypothesis in (2.7). For simplicity, we only present the result for row comparison here. The application to the column comparison is straightforward.

To test the hypothesis in (2.7), we define the test statistic as

$$M_{n,i} = \max_{1\leq j\leq q} M_{ij}.$$

The following theorem derives the limiting distribution of $M_{n,i}$ under the null hypothesis.

**Theorem 2.6.** If the null hypothesis $\mathbf{H_{0,i}}$ in (2.7) and conditions in Theorem 2.2 hold, we have

$$\mathbb{P}(M_{n,i} - 2\log q + \log\log q \le x) \to \exp\Big(-\frac{1}{\sqrt{\pi}}\exp\big(-\frac{x}{2}\big)\Big), \tag{2.19}$$

for any given $x \in \mathbb{R}$, as $n, q \to \infty$.

The above theorem can be proved in a similar way to Theorem 2.2.

**Remark 2.7.** For analyzing the power of $T_{\alpha,i}$, we define the following set of vector pairs,

$$\mathbb{A}_{i\star}(C) = \Big\{(\mathbf{u}_{1,i\star}, \mathbf{u}_{2,i\star}) : \max_{1\le j\le q} \frac{|u_{1,ij} - u_{2,ij}|}{\sqrt{m^2\zeta_{1,ij}/n_1 + m^2\zeta_{2,ij}/n_2}} \ge C\sqrt{\log q}\Big\}.$$

This allows us to yield a similar result to Theorem 2.4.

# 3. Applications to Testing Large Kendall's tau Correlation Matrix

In this section, we focus on testing the equality of two Kendall's tau matrices $\mathbf{U}_1^\tau$ and $\mathbf{U}_2^\tau$. This section contains two parts. In the first part, we assume the samples are from a fully nonparametric model. Under this model, in addition to the general Jackknife-based approach outlined in the previous section, we introduce two additional methods for testing (1.2) and analyze their theoretical properties (e.g., size and power). In the second part, we assume the samples are generated from a Gaussian copula model, under which we can relax Assumption **(A3)** to a much simplified form.

Kendall's tau provides a way to describe the nonlinear relationship between two random variables. As it is rank-based, it is especially suitable to analyze data from heavy-tailed or corrupted distributions. In this section, we aim to test the equality of two Kendall's tau matrices. More specifically, we set

$$\Phi_{ij}(\boldsymbol{X}_k, \boldsymbol{X}_\ell) := \operatorname{sign}(X_{ki} - X_{\ell i})\operatorname{sign}(X_{kj} - X_{\ell j}),$$
$$\Phi_{ij}(\boldsymbol{Y}_k, \boldsymbol{Y}_\ell) := \operatorname{sign}(Y_{ki} - Y_{\ell i})\operatorname{sign}(Y_{kj} - Y_{\ell j}),$$

and $q = d$. We aim to test whether $\mathbf{U}_1^\tau = \mathbf{U}_2^\tau$.

## 3.1. Methods and Theory under Fully Nonparametric Models

Section 3.1 contains two parts. The first part introduces two additional test procedures tailored for testing the equality of Kendall's tau matrices $\mathbf{U}_1^\tau$ and $\mathbf{U}_2^\tau$. The second part presents the theoretical properties of all the three tests. In addition, we further prove the rate-optimality of the proposed tests. Our technical contributions include providing an upper bound of the traditional plug-in variance estimation error, which enables us to establish the explicit rate of convergence of the plug-in variance estimator. We also prove an upper bound of the variance difference between two Kendall's tau correlation coefficients. These bounds allow us to derive the limiting distribution of the additional test statistic. The construction of these bounds requires the nontrivial usage of special structures of variance estimators and is of independent interest themselves. Moreover, for proving our test methods' optimality for the Kendall's tau matrix comparison, we construct a collection of least favourable multivariate normal distributions with regard to the test hypothesis. This novel construction technique is developed for correlation matrix comparison and is one of our technical contributions.

Recall that $\mathbf{U}_a^\tau$ and $\widehat{\mathbf{U}}_a^\tau$, defined in (2.3), are symmetric and we have $\mathrm{Diag}(\mathbf{U}_a^\tau) = \mathrm{Diag}(\widehat{\mathbf{U}}_a^\tau) = \mathbf{I}_d$. Therefore, we don't need to compare the main diagonals of $\mathbf{U}_a^\tau$. Hence, we reset $S = \{(i,j) : 1 \leq i < j \leq d\}$ for testing the equality of large Kendall's tau correlation matrices. The Jackknife-based statistic $M_n$ in Section 2.1 then becomes

$$M_n^{\tau,\mathrm{jack}} := \max_{(i,j)\in S} \frac{(\widehat{\tau}_{1,ij} - \widehat{\tau}_{2,ij})^2}{\widehat{\sigma}^2(\widehat{\tau}_{1,ij}) + \widehat{\sigma}^2(\widehat{\tau}_{2,ij})}. \tag{3.1}$$

Here, we still use $\widehat{\sigma}^2(\cdot)$ to denote the Jackknife variance estimator. Accordingly, we obtain $\mathrm{T}_\alpha^{\tau,\mathrm{jack}}$:

$$\mathrm{T}_\alpha^{\tau,\mathrm{jack}} := \mathbb{1}\left\{M_n^{\tau,\mathrm{jack}} \geq G^-(\alpha) + 4\log d - \log(\log d)\right\}.$$

#### *3.1.1. Three Procedures to Compare Kendall's tau Matrices*

In this section, we present two additional methods for comparing two Kendall's tau matrices. We start with the introduction of a plug-in method, which directly estimates the variances of $\{\widehat{\tau}_{a,ij}\}_{a=1,2}$ and plugs them into the test statistic. For this, recall that the Kendall's tau sample correlation between two random variable $U$ and $V$ is set as

$$\widehat{\tau} = \frac{2}{n(n-1)} \sum_{1\leq i<j\leq n} \mathrm{sign}(U_i - U_j)\,\mathrm{sign}(V_i - V_j),$$

where $U_1, \ldots, U_n$ and $V_1, \ldots, V_n$ are $n$ random samples from $U$ and $V$. Let $\Pi_c$ be the probability of the event that among two members drawn from the sample without replacement, they are concordant with each other. In other words, we have

$$\Pi_c = \mathbb{P}\big((U_2 - U_1)(V_2 - V_1) > 0\big). \tag{3.2}$$

Kruskal (1958) prove that the variance of $\widehat{\tau}$ can be written as

$$\frac{8}{n(n-1)}\Pi_c(1-\Pi_c) + 16\frac{1}{n}\frac{n-2}{n-1}(\Pi_{cc} - \Pi_c^2), \tag{3.3}$$

where $\Pi_{cc}$ is the probability of the event that among three members drawn from the sample without replacement, the second and third are concordant with the first. In other

words, we have

$$\Pi_{cc} = \mathbb{P}\Big(\big[(U_2 - U_1)(V_2 - V_1) > 0\big] \cap \big[U_3 - U_1)(V_3 - V_1) > 0\big]\Big). \tag{3.4}$$

As $n \to \infty$, the quantity in (3.3) multiplied by $n$ has the limit $16(\Pi_{cc} - \Pi_c^2)$. Motivated by this result, we propose the following plug-in variance estimator

$$\widehat{\sigma}^2_{\text{plug}}(\widehat{\tau}) = \frac{16}{n}(\widehat{\Pi}_{cc} - \widehat{\Pi}_c^2), \tag{3.5}$$

as an alternative to the Jackknife based one for estimating the variance of $\widehat{\tau}$. Here $\widehat{\Pi}_{cc}$ and $\widehat{\Pi}_c$ are the corresponding U-statistics to estimate $\Pi_{cc}$ and $\Pi_c$.[4] We replace $\widehat{\sigma}^2(\cdot)$ in $M_n^{\tau,\text{jack}}$ with $\widehat{\sigma}^2_{\text{plug}}(\cdot)$ to construct $M_n^{\tau,\text{plug}}$:

$$M_n^{\tau,\text{plug}} := \max_{1\le i<j\le d} \frac{(\widehat{\tau}_{1,ij} - \widehat{\tau}_{2,ij})^2}{\widehat{\sigma}^2_{\text{plug}}(\widehat{\tau}_{1,ij}) + \widehat{\sigma}^2_{\text{plug}}(\widehat{\tau}_{2,ij})}. \tag{3.6}$$

Accordingly, we construct the plug-in type test $\mathrm{T}_\alpha^{\tau,\text{plug}}$ as follows:

$$\mathrm{T}_\alpha^{\tau,\text{plug}} := \mathbb{1}\left\{M_n^{\tau,\text{plug}} \ge G^{-}(\alpha) + 4\log d - \log(\log d)\right\}.$$

In Section 3.1.2 we will provide the theoretical justification for this plug-in procedure.

Both the theoretical and numerical results indicate that the variance estimation error is also a key factor influencing the test statistics' powers. Up to now, we consider two kinds of variance estimation procedures (Jackknife based and plug-in based) for testing the equality of two Kendall's tau matrices. To exploit the sparsity of $\mathbf{U}^\tau$, we next propose to use the exact variance under the uncorrelated condition ($\tau = 0$). We name this procedure as "pseudo method". It calculates the variance of $\widehat{\tau}_{a,ij}$ by assuming $\tau_{a,ij} = 0$. We set $\widetilde{\sigma}^2_{1,\text{ps}}$ and $\widetilde{\sigma}^2_{2,\text{ps}}$ as the variances of $\sqrt{n_1}\widehat{\tau}_1$ and $\sqrt{n_2}\widehat{\tau}_2$, under $\tau_1 = 0$ and $\tau_2 = 0$. We also set

$$\sigma_{a,\text{ps}} := \lim_{n_a\to\infty} \widetilde{\sigma}_{a,\text{ps}} \qquad \text{for} \qquad a = 1, 2. \tag{3.7}$$

The test statistic becomes

$$M_n^{\tau,\text{ps}} := \max_{1\le i<j\le d} \frac{(\widehat{\tau}_{1,ij} - \widehat{\tau}_{2,ij})^2}{\sigma^2_{1,\text{ps}}/n_1 + \sigma^2_{2,\text{ps}}/n_2}. \tag{3.8}$$

Similarly, we construct the test $\mathrm{T}_\alpha^{\tau,\text{ps}}$:

$$\mathrm{T}_\alpha^{\tau,\text{ps}} := \mathbb{1}\left\{M_n^{\tau,\text{ps}} \ge G^{-}(\alpha) + 4\log d - \log(\log d)\right\}. \tag{3.9}$$

For example, if $\boldsymbol{X}$ and $\boldsymbol{Y}$ are generated from continuous Gaussian copula model, we have

$$\widetilde{\sigma}^2_{1,\text{ps}} = \frac{2(2n_1+5)}{9(n_1-1)},\ \widetilde{\sigma}^2_{2,\text{ps}} = \frac{2(2n_2+5)}{9(n_2-1)} \quad \text{and} \quad \sigma^2_{1,\text{ps}} = \sigma^2_{2,\text{ps}} = \frac{4}{9}.$$

**Remark 3.1.** As long as $|\widetilde{\sigma}^2_{a,\text{ps}} - \sigma^2_{a,\text{ps}}| = o((\log d)^{-1-\epsilon})$ with an arbitrary $\epsilon > 0$, we can show that replacing $\sigma^2_{a,\text{ps}}$ with $\widetilde{\sigma}^2_{a,\text{ps}}$ still gives a valid test. Details are provided in the proof of Theorem 3.3.

[4] By the definition of $\Pi_{cc}$ in (3.4), we should build a U-statistic with an asymmetric kernel to estimate it.

*3.1.2. Theoretical Properties of Three Testing Procedures*

We now present the theoretical properties (size, power, and optimality) of the three tests introduced in the former sections. More specifically, we prove their validity under the null hypothesis and conduct power analysis similarly to Theorems 2.2 and 2.4. Furthermore, we show that these tests are rate optimal against the sparse alternative.

In the beginning, the following theorem gives the limiting distribution for plug-in and Jackknife based test statistics.

**Theorem 3.2.** Assuming **(A1)**, **(A2)** and **(A3)** hold, under $\mathbf{H}_0^\tau$ of (1.2), we have

$$\mathbb{P}\big(M_n^{\tau,\text{jack}} - 4\log d + \log(\log d) \le x\big) \to \exp\Big(-\frac{1}{\sqrt{8\pi}}\exp(-\frac{x}{2})\Big), \tag{3.10}$$

$$\mathbb{P}\big(M_n^{\tau,\text{plug}} - 4\log d + \log(\log d) \le x\big) \to \exp\Big(-\frac{1}{\sqrt{8\pi}}\exp(-\frac{x}{2})\Big), \tag{3.11}$$

for any $x \in \mathbb{R}$, as $n, d \to \infty$. Furthermore, the results hold uniformly for all $\boldsymbol{X}$ and $\boldsymbol{Y}$ satisfying **(A1)**, **(A2)** and **(A3)**.

The following theorem gives the limiting distribution of the pseudo method. It holds under an additional meta-elliptical (defined in Supplement E of Supplementary Material) distributional assumption on the data.

**Theorem 3.3.** We assume that $\boldsymbol{X}$ and $\boldsymbol{Y}$ belong to the meta-elliptical distribution (Fang et al., 2002)[5]. If Assumptions **(A1)**, **(A2)** and **(A3)** hold, under $\mathbf{H}_0^\tau$ in (1.2), we have

$$\mathbb{P}\big(M_n^{\tau,\text{ps}} - 4\log d + \log(\log d) \le x\big) \to \exp\Big(-\frac{1}{\sqrt{8\pi}}\exp(-\frac{x}{2})\Big), \tag{3.12}$$

for any $x \in \mathbb{R}$, as $n, d \to \infty$. Furthermore, the result holds uniformly for all $\boldsymbol{X}$ and $\boldsymbol{Y}$ satisfying **(A1)**, **(A2)**, **(A3)**.

We now analyze the powers of $\mathrm{T}_\alpha^{\tau,\text{jack}}$, $\mathrm{T}_\alpha^{\tau,\text{plug}}$ and $\mathrm{T}_\alpha^{\tau,\text{ps}}$. Similarly to Theorem 2.4, we define

$$\mathbb{U}(C) = \Big\{(\mathbf{U}_1^\tau, \mathbf{U}_2^\tau) : \max_{1\le i<j\le d} \frac{|\tau_{1,ij} - \tau_{2,ij}|}{\sqrt{4\zeta_{1,ij}/n_1 + 4\zeta_{2,ij}/n_2}} \ge C\sqrt{\log d}\Big\},$$

$$\mathbb{V}(C) = \Big\{(\mathbf{U}_1^\tau, \mathbf{U}_2^\tau) : \max_{1\le i<j\le d} \frac{|\tau_{1,ij} - \tau_{2,ij}|}{\sqrt{\sigma_{1,\text{ps}}^2/n_1 + \sigma_{2,\text{ps}}^2/n_2}} \ge C\sqrt{\log d}\Big\}.$$

They are Kendall's tau versions of $\mathbb{A}(C)$ in Theorem 2.4.

---

[5] Detailed introduction of the meta-elliptical distribution family is provided in Supplement E of Supplementary Material.

**Theorem 3.4.** (Power Analysis) Assuming **(A2)** holds, we have

$$\inf_{(\mathbf{U}_1^\tau,\mathbf{U}_2^\tau)\in\mathbb{U}(4)} \mathbb{P}(\mathrm{T}_\alpha^{\tau,\mathrm{jack}}=1)\to 1, \tag{3.13}$$

$$\inf_{(\mathbf{U}_1^\tau,\mathbf{U}_2^\tau)\in\mathbb{U}(4)} \mathbb{P}(\mathrm{T}_\alpha^{\tau,\mathrm{plug}}=1)\to 1, \tag{3.14}$$

as $n,d\to\infty$. If $\boldsymbol{X}$ and $\boldsymbol{Y}$ belong to the meta-elliptical family and **(A1)**, **(A2)** are satisfied, as $n,d\to\infty$, we have

$$\inf_{(\mathbf{U}_1^\tau,\mathbf{U}_2^\tau)\in\mathbb{V}(4)} \mathbb{P}(\mathrm{T}_\alpha^{\tau,\mathrm{ps}}=1)\to 1. \tag{3.15}$$

Theorem 3.4 implies that just one entry of $\mathbf{U}_1^\tau-\mathbf{U}_2^\tau$ has a magnitude no smaller than $C\sqrt{\log d/n}$ is enough for the introduced tests to correctly reject $\mathbf{H}_\mathbf{0}^\tau$.

Next, we show that all the three proposed methods are rate optimal by matching the obtained rates of convergence to a lower bound for correlation matrix comparison. We adopt the general framework used in Baraud (2002) to obtain the lower bound for testing the equality of correlation matrices. The core of the proof is the construction of collections of least favourable multivariate normal distributions with regard to the test hypothesis. Our work is related to Cai et al. (2013) which prove the lower bound for testing the equality of covariance matrices. However, their construction technique is developed for covariance matrices but not the correlation matrices. Specifically, they only perturb the diagonal elements of the covariance, which does not affect the resulting correlation matrices. To test correlation matrices, we need to develop a novel construction by perturbing the off-diagonal elements of the correlation matrices. Details are provided in the proof of Theorem 3.5.

**Theorem 3.5.** Let $\alpha,\beta>0$ and $\alpha+\beta<1$. Assuming that $\log d/n=o(1)$, there exits a sufficiently small positive number $c_0$, such that for any distribution family that contains Gaussian as a subfamily, and all large enough $n$ and $d$, we have

$$\inf_{(\mathbf{U}_1^\tau,\mathbf{U}_2^\tau)\in\mathbb{U}(c_0)} \sup_{T_\alpha\in\mathcal{T}_\alpha} \mathbb{P}(T_\alpha=1)\le 1-\beta, \tag{3.16}$$

where $\mathcal{T}_\alpha$ represents all level $\alpha$ tests for testing the equality of two correlation matrices.

Cai et al. (2013) give a similar result for testing the equality of two covariance matrices. They show that the rate $C\sqrt{\log d/n}$ is optimal for comparing covariance matrices under conditions that $\boldsymbol{X}$ and $\boldsymbol{Y}$ have sub-Gaussian-type or polynomial-type tails. In comparison, the lower bound result in Theorem 3.5 illustrates that our proposed methods are rate optimal under the fully nonparametric model. In particular, we don't impose assumptions on the marginal distributions.

## 3.2. Methods and Theory under Semiparametric Gaussian Copula Models

In this section, we assume that $\boldsymbol{X}$ and $\boldsymbol{Y}$ are $d$-dimensional random vectors from the Gaussian copula with latent correlation matrices $\mathbf{\Sigma}_a=(\sigma_{a,ij})$, $a=1,2$ and $\mathrm{Diag}(\mathbf{\Sigma}_a)=$ $\mathbf{I}_d$[6]. Under the Gaussian copula model, the technical assumption **(A3)** in Section 2.2 can be replaced by a much simplified condition. Specifically, for $r \in (0,1)$, we define

$$\Omega(r) := \{1 \leq i \leq d : |\tau_{1,ij}| > r \ \text{ or } \ |\tau_{2,ij}| > r \ \text{ for some } \ j \neq i\}. \tag{3.17}$$

We describe the technical assumption **(A4)** as follows:

> **(A4)**. For some $r < 1$ and a sequence of numbers $\Omega_{d,r} = o(d)$, we have $|\Omega(r)| \leq \Omega_{d,r}$.

After introducing Assumption **(A4)**, we then discuss its relationship with Assumptions **(A1)** and **(A3)**. For **(A1)**, although it has similar form to **(A4)**, they are essentially different. Assumption **(A1)** is related to the largest eigenvalues of $\mathbf{U}^\tau_\gamma$. In fact, bounded $\lambda_{\max}(\mathbf{U}^\tau_\gamma)$ implies $\max_{1\leq j\leq d} \mathrm{supp}_j(\alpha_0) \leq C(\log d)^{2+2\alpha_0}$. On the contrary, Assumption **(A4)** is related to $\lambda_{\min}(\mathbf{U}^\tau_\gamma)$. For example, if the correlation between two Gaussian random variables goes to 1, the corresponding correlation matrix will be asymptotically degenerated with the least eigenvalue infinite small. In proof, we first use Assumption **(A4)** to select largest sub-matrix of $\mathbf{U}^\tau_\gamma$ so that all its entries' absolute values are less than $r$. We then use Assumption **(A1)** to exclude the influence of entires with $|\tau_{\gamma,ij}| \geq (\log d)^{-1-\alpha_0}$ on the asymptotic results.

Assumptions **(A4)** and **(A3)** are highly related. However, Assumption **(A3)** cannot be straightforwardly implied by Assumptions **(A1)**, **(A2)** and **(A4)**. In fact, the relationship between **(A3)** and **(A4)** is complicated. To see the exact relationship, we need some additional definitions.

First, we have $S = \{(i,j) : 1 \leq i < j \leq d\}$. We then define

$$C_0 = \{(i,j) : i \in \Omega(r) \bigcup \Gamma\} \bigcup \{(i,j) : j \in \Omega(r) \bigcup \Gamma\} \qquad \text{and} \qquad B_0 = S_0 \bigcup C_0,$$

where $\Gamma$ is defined in in Assumption **(A1)** and $S_0$ is defined in (2.11). Furthermore, we denote $A$ to be the biggest subset of $S \setminus B_0$, such that any two pairs $(i,j) \neq (k,\ell) \in A$ must satisfy a condition $(\star)$. More detailed description of condition $(\star)$ will be provided in the proof of Theorem 3.6. Essentially, it specifies that, for any $(i,j) \neq (k,\ell) \in S \setminus B_0$, there exits an $i_1 \in \{i,j,k,\ell\}$ such that for any $j_1 \in \{i,j,k,\ell\} \setminus i_1$, we have $|\tau_{a,i_1j_1}| = O((\log d)^{-1-\alpha_0})$. We also define $\tau_{a,ijk\ell}$ as the Kendall's tau version of $u_{a,ijk\ell}$ in Assumption **(A3)**.

Under Assumptions **(A1)**, **(A2)** and **(A4)**, we can prove that for any $(i,j) \neq (k,\ell) \in A$, we have $|\tau_{a,ijk\ell}| = O((\log d)^{-1-\alpha_0})$, which is essentially Assumption **(A3)** with $u_{a,ijk\ell}$ replaced by $\tau_{a,ijk\ell}$. The only difference is that these conditions hold on $A$ but instead of $S \setminus S_0$ as in Assumption **(A3)**. Theorem 3.6 below specifics that Assumptions **(A1)**, **(A2)** and **(A4)** can be used to replace Assumptions **(A1)**, **(A2)** and **(A3)** when we test the equality of Kendall's tau correlation matrices under the Gaussian copula model.

**Theorem 3.6.** *Let $\boldsymbol{X}$ and $\boldsymbol{Y}$ be Gaussian copula random vectors with latent correlation matrices $\boldsymbol{\Sigma}_a$, $a = 1$ or $2$ and $\mathrm{Diag}(\boldsymbol{\Sigma}_a) = \mathbf{I}_d$. We assume that the smallest eigenvalue* of any 4 by 4 principal sub-matrix of $\mathbf{\Sigma}_a$ is uniformly bounded away from 0. Assuming **(A1)**, **(A2)** and **(A4)** hold, under $\mathbf{H_0^\tau}$ of (1.2), we have

$$\mathbb{P}\Big(M_n^{\tau,\text{jack}} - 4\log d + \log(\log d) \leq x\Big) \to \exp\Big(-\frac{1}{\sqrt{8\pi}}\exp(-\frac{x}{2})\Big),$$
$$\mathbb{P}\Big(M_n^{\tau,\text{plug}} - 4\log d + \log(\log d) \leq x\Big) \to \exp\Big(-\frac{1}{\sqrt{8\pi}}\exp(-\frac{x}{2})\Big),$$
$$\mathbb{P}\Big(M_n^{\tau,\text{ps}} - 4\log d + \log(\log d) \leq x\Big) \to \exp\Big(-\frac{1}{\sqrt{8\pi}}\exp(-\frac{x}{2})\Big),$$

for any $x \in \mathbb{R}$, as $n, d \to \infty$. Furthermore, these limiting results hold uniformly for all $\boldsymbol{X}$ and $\boldsymbol{Y}$ satisfying **(A1)**, **(A2)** and **(A4)**.

[6] Detailed definition of the Gaussian copula is put in Supplement E of Supplement Material.

***Proof.*** Recall that $M_n^{\tau,\text{jack}}$, $M_n^{\tau,\text{plug}}$ and $M_n^{\tau,\text{ps}}$ in (3.1), (3.6) and (3.8) are defined by taking maximum over $S$. The main idea is to show that it is sufficient to use a version of these quantities taking the maximum over the smaller set $A$ as defined before. The proof is technical and left to Supplement A.6 of Supplementary Material. □

**Remark 3.7.** In Supplement E of Supplementary Material, we show that $\mathbf{U}_a^\tau$ and $\mathbf{\Sigma}_a$ are related in terms of $\sigma_{a,ij} = \sin(\tau_{a,ij}\pi/2)$. Hence, testing (1.2) is equivalent to testing

$$\mathbf{H_0} : \mathbf{\Sigma}_1 = \mathbf{\Sigma}_2 \;\; \text{v.s.} \;\; \mathbf{H_1} : \mathbf{\Sigma}_1 \neq \mathbf{\Sigma}_2,$$

under the Gaussian copula model.

**Remark 3.8.** To test the row or column of Kendall's tau matrices, if any of the conditions of Theorems 3.2, 3.3 and 3.6 hold, we get the same limiting result as in (2.19).

# 4. Experiments

In this section, we demonstrate numerical performances of proposed methods on simulated and real data sets. In particular, we compare proposed methods with the state-of-the-art method in the literature.

## 4.1. Numerical Simulations

We compare proposed methods with the sample covariance based method (denoted by $\mathrm{T}_\alpha^{\mathrm{CLX}}$) in Cai et al. (2013). To test our methods under various covariance structures, we introduce the following matrices.

- (Block matrix $\mathbf{\Sigma}^*$) Let $\mathbf{R}^* = (r_{ij}^*) \in \mathbb{R}^{d\times d}$ with $r_{ij}^* = 0.6$ for $5(k-1)+1 \leq i \neq j \leq 5k$ and $k = 1, \ldots, \lfloor d/5 \rfloor$. For other entries in $\mathbf{R}^*$, we set $r_{ii}^* = 1$ and $r_{ij}^* = 0$ when $i \neq j$. Let $\mathbf{D}$ as a diagonal matrix with each nonzero entry following independent uniform distribution on the interval $(0.5, 1.5)$. We then set $\mathbf{\Sigma}^* = \mathbf{D}\mathbf{R}^*\mathbf{D}$.

**Table 1.**. Empirical sizes of **Model 1**, **2** and **3** under $\alpha = 0.05$ based on 2000 repetitions.

| | | $\mathbf{\Sigma}^*$ | | | | | | | $\mathbf{\Sigma}'$ | | | | | | | $\mathbf{\Sigma}^\star$ | | | | | | |
|---|---|---|---|---|---|---|---|---|---|---|---|---|---|---|---|---|---|---|---|---|---|---|
| n | d | 50 | 100 | 200 | 300 | 500 | 700 | 1000 | 50 | 100 | 200 | 300 | 500 | 700 | 1000 | 50 | 100 | 200 | 300 | 500 | 700 | 1000 |
| | | | | | | | | | | | | **Model 1** | | | | | | | | | | |
| 500 | $T_\alpha^{\tau,\text{plug}}$ | 0.05 | 0.05 | 0.07 | 0.07 | 0.07 | 0.08 | 0.09 | 0.05 | 0.06 | 0.06 | 0.07 | 0.07 | 0.08 | 0.09 | 0.05 | 0.05 | 0.05 | 0.06 | 0.06 | 0.08 | 0.07 |
| | $T_\alpha^{\tau,\text{jack}}$ | 0.05 | 0.05 | 0.06 | 0.06 | 0.06 | 0.06 | 0.07 | 0.05 | 0.05 | 0.05 | 0.05 | 0.06 | 0.06 | 0.08 | 0.04 | 0.05 | 0.05 | 0.05 | 0.05 | 0.06 | 0.06 |
| | $T_\alpha^{\tau,\text{ps}}$ | 0.05 | 0.04 | 0.05 | 0.05 | 0.04 | 0.05 | 0.05 | 0.05 | 0.05 | 0.05 | 0.05 | 0.05 | 0.04 | 0.05 | 0.04 | 0.04 | 0.04 | 0.05 | 0.04 | 0.05 | 0.05 |
| | $T_\alpha^{\text{CLX}}$ | 0.04 | 0.04 | 0.04 | 0.05 | 0.04 | 0.05 | 0.05 | 0.05 | 0.05 | 0.04 | 0.05 | 0.04 | 0.04 | 0.05 | 0.03 | 0.03 | 0.05 | 0.05 | 0.05 | 0.06 | 0.06 |
| | | | | | | | | | | | | **Model 2** | | | | | | | | | | |
| | $T_\alpha^{\tau,\text{plug}}$ | 0.05 | 0.06 | 0.06 | 0.07 | 0.07 | 0.08 | 0.09 | 0.05 | 0.06 | 0.07 | 0.07 | 0.08 | 0.08 | 0.09 | 0.05 | 0.05 | 0.05 | 0.06 | 0.06 | 0.08 | 0.08 |
| | $T_\alpha^{\tau,\text{jack}}$ | 0.05 | 0.05 | 0.05 | 0.06 | 0.06 | 0.06 | 0.07 | 0.05 | 0.05 | 0.06 | 0.06 | 0.07 | 0.07 | 0.08 | 0.05 | 0.05 | 0.05 | 0.05 | 0.06 | 0.06 | 0.06 |
| | $T_\alpha^{\tau,\text{ps}}$ | 0.04 | 0.05 | 0.05 | 0.05 | 0.05 | 0.05 | 0.04 | 0.04 | 0.04 | 0.04 | 0.04 | 0.05 | 0.05 | 0.05 | 0.03 | 0.04 | 0.04 | 0.04 | 0.05 | 0.05 | 0.05 |
| | $T_\alpha^{\text{CLX}}$ | 0.00 | 0.00 | 0.00 | 0.00 | 0.00 | 0.00 | 0.00 | 0.00 | 0.00 | 0.00 | 0.00 | 0.00 | 0.00 | 0.00 | 0.00 | 0.00 | 0.00 | 0.00 | 0.00 | 0.00 | 0.00 |
| | | | | | | | | | | | | **Model 3** | | | | | | | | | | |
| | $T_\alpha^{\tau,\text{plug}}$ | 0.05 | 0.06 | 0.06 | 0.07 | 0.07 | 0.08 | 0.09 | 0.05 | 0.05 | 0.06 | 0.07 | 0.08 | 0.08 | 0.09 | 0.05 | 0.05 | 0.05 | 0.06 | 0.06 | 0.08 | 0.08 |
| | $T_\alpha^{\tau,\text{jack}}$ | 0.05 | 0.05 | 0.05 | 0.06 | 0.06 | 0.06 | 0.07 | 0.05 | 0.05 | 0.05 | 0.06 | 0.07 | 0.07 | 0.08 | 0.04 | 0.04 | 0.04 | 0.05 | 0.05 | 0.06 | 0.06 |
| | $T_\alpha^{\tau,\text{ps}}$ | 0.04 | 0.05 | 0.05 | 0.05 | 0.05 | 0.05 | 0.05 | 0.04 | 0.04 | 0.04 | 0.05 | 0.05 | 0.05 | 0.05 | 0.03 | 0.03 | 0.04 | 0.04 | 0.04 | 0.05 | 0.05 |
| | $T_\alpha^{\text{CLX}}$ | 0.00 | 0.01 | 0.02 | 0.05 | 0.08 | 0.16 | 0.22 | 0.00 | 0.01 | 0.02 | 0.05 | 0.10 | 0.16 | 0.25 | 0.00 | 0.01 | 0.02 | 0.04 | 0.08 | 0.15 | 0.23 |
| | | | | | | | | | | | | **Model 1** | | | | | | | | | | |
| 200 | $T_\alpha^{\tau,\text{plug}}$ | 0.06 | 0.08 | 0.09 | 0.12 | 0.14 | 0.16 | 0.19 | 0.06 | 0.08 | 0.11 | 0.13 | 0.15 | 0.16 | 0.19 | 0.05 | 0.06 | 0.08 | 0.09 | 0.12 | 0.15 | 0.14 |
| | $T_\alpha^{\tau,\text{jack}}$ | 0.05 | 0.05 | 0.06 | 0.08 | 0.10 | 0.10 | 0.12 | 0.05 | 0.06 | 0.07 | 0.09 | 0.10 | 0.11 | 0.12 | 0.03 | 0.05 | 0.05 | 0.07 | 0.08 | 0.10 | 0.10 |
| | $T_\alpha^{\tau,\text{ps}}$ | 0.04 | 0.04 | 0.04 | 0.04 | 0.05 | 0.05 | 0.05 | 0.05 | 0.04 | 0.05 | 0.06 | 0.05 | 0.05 | 0.05 | 0.04 | 0.04 | 0.05 | 0.05 | 0.05 | 0.05 | 0.04 |
| | $T_\alpha^{\text{CLX}}$ | 0.04 | 0.04 | 0.04 | 0.04 | 0.05 | 0.05 | 0.06 | 0.04 | 0.04 | 0.04 | 0.05 | 0.05 | 0.05 | 0.06 | 0.04 | 0.04 | 0.05 | 0.05 | 0.06 | 0.05 | 0.06 |
| | | | | | | | | | | | | **Model 2** | | | | | | | | | | |
| | $T_\alpha^{\tau,\text{plug}}$ | 0.06 | 0.08 | 0.09 | 0.12 | 0.14 | 0.16 | 0.19 | 0.06 | 0.08 | 0.11 | 0.13 | 0.15 | 0.16 | 0.20 | 0.05 | 0.06 | 0.08 | 0.09 | 0.12 | 0.15 | 0.15 |
| | $T_\alpha^{\tau,\text{jack}}$ | 0.05 | 0.05 | 0.06 | 0.08 | 0.10 | 0.10 | 0.12 | 0.05 | 0.06 | 0.07 | 0.09 | 0.10 | 0.11 | 0.12 | 0.03 | 0.05 | 0.05 | 0.07 | 0.08 | 0.10 | 0.10 |
| | $T_\alpha^{\tau,\text{ps}}$ | 0.04 | 0.04 | 0.04 | 0.04 | 0.05 | 0.05 | 0.05 | 0.04 | 0.04 | 0.05 | 0.06 | 0.05 | 0.05 | 0.05 | 0.04 | 0.04 | 0.05 | 0.05 | 0.05 | 0.05 | 0.05 |
| | $T_\alpha^{\text{CLX}}$ | 0.00 | 0.00 | 0.00 | 0.00 | 0.00 | 0.00 | 0.00 | 0.00 | 0.00 | 0.00 | 0.00 | 0.00 | 0.00 | 0.00 | 0.00 | 0.00 | 0.00 | 0.00 | 0.00 | 0.00 | 0.00 |
| | | | | | | | | | | | | **Model 3** | | | | | | | | | | |
| | $T_\alpha^{\tau,\text{plug}}$ | 0.06 | 0.08 | 0.09 | 0.12 | 0.14 | 0.16 | 0.19 | 0.06 | 0.08 | 0.11 | 0.13 | 0.15 | 0.16 | 0.20 | 0.05 | 0.06 | 0.08 | 0.09 | 0.12 | 0.15 | 0.15 |
| | $T_\alpha^{\tau,\text{jack}}$ | 0.05 | 0.05 | 0.06 | 0.08 | 0.10 | 0.10 | 0.12 | 0.05 | 0.06 | 0.07 | 0.09 | 0.10 | 0.11 | 0.12 | 0.03 | 0.05 | 0.05 | 0.07 | 0.08 | 0.10 | 0.10 |
| | $T_\alpha^{\tau,\text{ps}}$ | 0.04 | 0.04 | 0.04 | 0.04 | 0.05 | 0.05 | 0.05 | 0.04 | 0.04 | 0.05 | 0.06 | 0.05 | 0.05 | 0.05 | 0.04 | 0.04 | 0.05 | 0.05 | 0.05 | 0.05 | 0.05 |
| | $T_\alpha^{\text{CLX}}$ | 0.00 | 0.01 | 0.02 | 0.06 | 0.12 | 0.23 | 0.37 | 0.00 | 0.00 | 0.02 | 0.07 | 0.14 | 0.26 | 0.37 | 0.00 | 0.01 | 0.01 | 0.08 | 0.15 | 0.28 | 0.39 |

**Table 2.**. Empirical powers of **Model 1**, **2** and **3** under $\alpha = 0.05$ based on 2000 repetitions. We set $\zeta = 0.2$ for $n = 500$ and $\zeta = 0.3$ for $n = 200$.

| | | $\boldsymbol{\Sigma}^*$ | | | | | | | $\boldsymbol{\Sigma}'$ | | | | | | | $\boldsymbol{\Sigma}^\star$ | | | | | | |
|---|---|---|---|---|---|---|---|---|---|---|---|---|---|---|---|---|---|---|---|---|---|---|
| n | d | 50 | 100 | 200 | 300 | 500 | 700 | 1000 | 50 | 100 | 200 | 300 | 500 | 700 | 1000 | 50 | 100 | 200 | 300 | 500 | 700 | 1000 |
| | | | | | | | | | | | | **Model 1** | | | | | | | | | | |
| 500 | $T_\alpha^{\tau,plug}$ | 0.85 | 0.81 | 0.76 | 0.76 | 0.76 | 0.74 | 0.71 | 0.86 | 0.82 | 0.78 | 0.76 | 0.72 | 0.72 | 0.71 | 0.86 | 0.81 | 0.75 | 0.74 | 0.71 | 0.70 | 0.67 |
| | $T_\alpha^{\tau,jack}$ | 0.84 | 0.80 | 0.75 | 0.75 | 0.75 | 0.73 | 0.70 | 0.86 | 0.81 | 0.77 | 0.76 | 0.72 | 0.71 | 0.70 | 0.87 | 0.81 | 0.75 | 0.73 | 0.70 | 0.69 | 0.64 |
| | $T_\alpha^{\tau,ps}$ | 0.83 | 0.79 | 0.74 | 0.73 | 0.73 | 0.71 | 0.68 | 0.85 | 0.81 | 0.75 | 0.74 | 0.70 | 0.68 | 0.67 | 0.85 | 0.79 | 0.72 | 0.70 | 0.67 | 067 | 0.61 |
| | $T_\alpha^{CLX}$ | 0.83 | 0.78 | 0.74 | 0.74 | 0.72 | 0.72 | 0.67 | 0.86 | 0.80 | 0.77 | 0.75 | 0.71 | 0.68 | 0.66 | 0.80 | 0.77 | 0.71 | 0.70 | 0.66 | 0.66 | 0.60 |
| | | | | | | | | | | | | **Model 2** | | | | | | | | | | |
| | $T_\alpha^{\tau,plug}$ | 0.78 | 0.72 | 0.70 | 0.68 | 0.66 | 0.65 | 0.62 | 0.77 | 0.72 | 0.68 | 0.66 | 0.63 | 0.63 | 0.60 | 0.79 | 0.71 | 0.63 | 0.62 | 0.58 | 0.58 | 0.55 |
| | $T_\alpha^{\tau,jack}$ | 0.77 | 0.72 | 0.69 | 0.68 | 0.67 | 0.66 | 0.61 | 0.76 | 0.72 | 0.67 | 0.66 | 0.62 | 0.61 | 0.59 | 0.79 | 0.70 | 0.64 | 0.62 | 0.57 | 0.57 | 0.54 |
| | $T_\alpha^{\tau,ps}$ | 0.76 | 0.70 | 0.67 | 0.66 | 0.65 | 0.64 | 0.61 | 0.75 | 0.69 | 0.66 | 0.63 | 0.60 | 0.60 | 0.58 | 0.75 | 0.67 | 0.58 | 0.57 | 0.53 | 0.52 | 0.52 |
| | $T_\alpha^{CLX}$ | 0.12 | 0.09 | 0.07 | 0.06 | 0.03 | 0.02 | 0.02 | 0.09 | 0.06 | 0.05 | 0.04 | 0.02 | 0.02 | 0.02 | 0.07 | 0.04 | 0.03 | 0.02 | 0.01 | 0.01 | 0.01 |
| | | | | | | | | | | | | **Model 3** | | | | | | | | | | |
| | $T_\alpha^{\tau,plug}$ | 0.84 | 0.81 | 0.78 | 0.77 | 0.74 | 0.72 | 0.71 | 0.85 | 0.81 | 0.80 | 0.78 | 0.72 | 0.72 | 0.71 | 0.87 | 0.82 | 0.74 | 0.75 | 0.70 | 0.70 | 0.68 |
| | $T_\alpha^{\tau,jack}$ | 0.84 | 0.80 | 0.78 | 0.76 | 0.73 | 0.71 | 0.70 | 0.85 | 0.81 | 0.78 | 0.77 | 0.71 | 0.70 | 0.68 | 0.87 | 0.82 | 0.77 | 0.74 | 0.69 | 0.69 | 0.65 |
| | $T_\alpha^{\tau,ps}$ | 0.82 | 0.80 | 0.75 | 0.75 | 0.72 | 0.69 | 0.68 | 0.84 | 0.81 | 0.77 | 0.75 | 0.70 | 0.69 | 0.67 | 0.85 | 0.79 | 0.73 | 0.72 | 0.66 | 0.66 | 0.61 |
| | $T_\alpha^{CLX}$ | 0.00 | 0.01 | 0.02 | 0.04 | 0.08 | 0.16 | 0.23 | 0.00 | 0.00 | 0.04 | 0.05 | 0.10 | 0.15 | 0.24 | 0.01 | 0.01 | 0.03 | 0.04 | 0.09 | 0.16 | 0.25 |
| | | | | | | | | | | | | **Model 1** | | | | | | | | | | |
| 200 | $T_\alpha^{\tau,plug}$ | 0.86 | 0.81 | 0.79 | 0.78 | 0.76 | 0.74 | 0.73 | 0.80 | 0.76 | 0.74 | 0.71 | 0.70 | 0.68 | 0.68 | 0.81 | 0.72 | 0.66 | 0.64 | 0.61 | 0.58 | 0.56 |
| | $T_\alpha^{\tau,jack}$ | 0.85 | 0.80 | 0.79 | 0.73 | 0.73 | 0.71 | 0.69 | 0.79 | 0.74 | 0.71 | 0.67 | 0.66 | 0.64 | 0.63 | 0.79 | 0.70 | 0.62 | 0.60 | 0.56 | 0.50 | 0.52 |
| | $T_\alpha^{\tau,ps}$ | 0.82 | 0.77 | 0.75 | 0.68 | 0.68 | 0.64 | 0.60 | 0.77 | 0.71 | 0.66 | 0.63 | 0.61 | 0.57 | 0.55 | 0.75 | 0.64 | 0.56 | 0.54 | 0.47 | 0.43 | 0.42 |
| | $T_\alpha^{CLX}$ | 0.76 | 0.72 | 0.68 | 0.62 | 0.60 | 0.53 | 0.50 | 0.73 | 0.68 | 0.62 | 0.58 | 0.51 | 0.50 | 0.47 | 0.63 | 0.58 | 0.50 | 0.46 | 0.39 | 0.32 | 0.31 |
| | | | | | | | | | | | | **Model 2** | | | | | | | | | | |
| | $T_\alpha^{\tau,plug}$ | 0.79 | 0.74 | 0.72 | 0.67 | 0.62 | 0.60 | 0.59 | 0.73 | 0.66 | 0.65 | 0.61 | 0.58 | 0.56 | 0.56 | 0.70 | 0.60 | 0.54 | 0.48 | 0.46 | 0.45 | 0.45 |
| | $T_\alpha^{\tau,jack}$ | 0.78 | 0.71 | 0.70 | 0.64 | 0.59 | 0.55 | 0.55 | 0.72 | 0.64 | 0.61 | 0.58 | 0.55 | 0.52 | 0.51 | 0.69 | 0.58 | 0.50 | 0.46 | 0.44 | 0.40 | 0.39 |
| | $T_\alpha^{\tau,ps}$ | 0.73 | 0.64 | 0.63 | 0.58 | 0.53 | 0.50 | 0.50 | 0.69 | 0.58 | 0.57 | 0.53 | 0.49 | 0.46 | 0.45 | 0.62 | 0.51 | 0.42 | 0.37 | 0.35 | 0.34 | 0.34 |
| | $T_\alpha^{CLX}$ | 0.07 | 0.04 | 0.02 | 0.01 | 0.01 | 0.01 | 0.00 | 0.07 | 0.04 | 0.02 | 0.01 | 0.00 | 0.00 | 0.00 | 0.03 | 0.02 | 0.01 | 0.00 | 0.00 | 0.00 | 0.00 |
| | | | | | | | | | | | | **Model 3** | | | | | | | | | | |
| | $T_\alpha^{\tau,plug}$ | 0.85 | 0.82 | 0.80 | 0.78 | 0.76 | 0.74 | 0.73 | 0.80 | 0.76 | 0.74 | 0.71 | 0.70 | 0.68 | 0.68 | 0.81 | 0.73 | 0.66 | 0.64 | 0.61 | 0.58 | 0.56 |
| | $T_\alpha^{\tau,jack}$ | 0.85 | 0.80 | 0.79 | 0.75 | 0.73 | 0.71 | 0.69 | 0.79 | 0.74 | 0.72 | 0.68 | 0.66 | 0.64 | 0.63 | 0.79 | 0.71 | 0.62 | 0.60 | 0.56 | 0.51 | 0.52 |
| | $T_\alpha^{\tau,ps}$ | 0.82 | 0.78 | 0.75 | 0.70 | 0.68 | 0.64 | 0.61 | 0.77 | 0.71 | 0.66 | 0.61 | 0.57 | 0.57 | 0.55 | 0.75 | 0.65 | 0.56 | 0.54 | 0.47 | 0.44 | 0.42 |
| | $T_\alpha^{CLX}$ | 0.00 | 0.00 | 0.03 | 0.07 | 0.16 | 0.22 | 0.37 | 0.00 | 0.01 | 0.02 | 0.06 | 0.15 | 0.26 | 0.37 | 0.00 | 0.00 | 0.00 | 0.06 | 0.16 | 0.26 | 0.37 |

- (Tridiagonal matrix $\mathbf{\Sigma}'$) Let $\mathbf{R}' = (r'_{ij}) \in \mathbb{R}^{d\times d}$ be a tridiagonal matrix with 1 on the main diagonal and 0.5 on the first diagonal. We then set $\mathbf{\Sigma}' = \mathbf{D}\mathbf{R}'\mathbf{D}$.
- (Multidiagonal matrix $\mathbf{\Sigma}^\star$) Let $\mathbf{R}^\star = (r^\star_{ij}) \in \mathbb{R}^{d\times d}$ with $r_{ij} = 0.8^{|i-j|}$ and $\mathbf{\Sigma}^\star = \mathbf{D}\mathbf{R}^\star\mathbf{D}$.

Under the null hypothesis, we sample $n_1+n_2$ data points from the following 3 models with $\mathbf{\Sigma} = \mathbf{\Sigma}^*, \mathbf{\Sigma}'$, and $\mathbf{\Sigma}^\star$.

- **Model 1** (Normal distribution) In this model, under the null hypothesis we generate $n_1+n_2$ random vectors from $N(\mathbf{0},\mathbf{\Sigma})$.
- **Model 2** (Multivariate $t$ distribution) We sample from $\boldsymbol{\mu} + \boldsymbol{Z}/\sqrt{W/\nu}$ with $W \sim \chi^2(\nu)$ and $\boldsymbol{Z} \sim N(\mathbf{0},\mathbf{\Sigma})$, where $W$ and $\boldsymbol{Z}$ are independent. Under the null hypothesis, we generate $n_1+n_2$ data points with $\boldsymbol{\mu} = \mathbf{0}$ and $\nu = 3$.
- **Model 3** (Marginal Cauchy distribution) Generate $n_1+n_2$ random vectors from $N(\mathbf{0},\mathbf{\Sigma})$. We then use a monotone function to transform each coordinate to follow the Cauchy distribution Cauchy$(\mu,s)$ whose density function is $s/\pi(s^2+(x-\mu)^2)$. In the simulation, we set $\mu = 0$ and $s = 1$.

Under above models, the two populations of $\boldsymbol{X}$ and $\boldsymbol{Y}$ have the same covariance matrices. We use them to show that our proposed methods can control the size correctly under the null hypothesis. For the power analysis, we introduce a random symmetric matrix $\mathbf{\Delta} = (\delta_{k\ell}) \in \mathbb{R}^{d\times d}$ with exactly 8 nonzero entries. Among the 8 entries, 4 entries are randomly selected from the upper triangle of $\mathbf{\Delta}$, with a magnitude generated from the uniform distribution on $(0, \zeta\sigma^2_{\max})$, where $\sigma^2_{\max}$ is the maximal value of $\mathbf{\Sigma}$'s main diagonal. Other 4 entries are determined by symmetry. We then set $\widetilde{\mathbf{\Sigma}}_1 = \mathbf{\Sigma} + \delta\mathbf{I}$ and $\widetilde{\mathbf{\Sigma}}_2 = \mathbf{\Sigma} + \mathbf{\Delta} + \delta\mathbf{I}$ with $\delta = |\min\{\lambda_{\min}(\mathbf{\Sigma}+\mathbf{\Delta}), \lambda_{\min}(\mathbf{\Sigma})\}| + 0.05$. In place of $\mathbf{\Sigma}$, we use the matrices $\widetilde{\mathbf{\Sigma}}_1$ and $\widetilde{\mathbf{\Sigma}}_2$ to generate samples for $\boldsymbol{X}$ and $\boldsymbol{Y}$ under the alternative hypothesis.

We set $n_1 = n_2 = n$ with $n = 200, 500$ and $d = 50, 100, 200, 300, 500, 700, 1000$. The nominal significance level $\alpha$ is 0.05. Table 1 presents empirical sizes. We see that $\mathrm{T}^{\tau,\mathrm{ps}}_\alpha$ always attains the desired size even for extremely large $d$. When $d$ is significantly larger than $n$, both $\mathrm{T}^{\tau,\mathrm{plug}}_\alpha$ and $\mathrm{T}^{\tau,\mathrm{jack}}_\alpha$ suffer from the size distortion. When $d$ approximates $n$, $\mathrm{T}^{\tau,\mathrm{jack}}_\alpha$ is still valid but $\mathrm{T}^{\tau,\mathrm{plug}}_\alpha$ fails. These size distortions decrease as $n$ increases. Although the theoretical limiting results are similar for all the proposed methods, the simulation results show that the estimation errors of variance heavily affect the proposed tests' finite sample performances, and $\mathrm{T}^{\tau,\mathrm{ps}}_\alpha$ benefits a lot from avoiding estimating the variance directly. Moreover, for heavy tail distributions such as multivariate $t$ and Cauchy distributions, we also see that $\mathrm{T}^{\mathrm{CLX}}_\alpha$ from Cai et al. (2013) becomes too conservative.

By examining the empirical powers in Table 2, for distributions with heavy tails or strong tail dependence, $\mathrm{T}^{\mathrm{CLX}}_\alpha$'s power decreases dramatically, making $\mathrm{T}^{\mathrm{CLX}}_\alpha$ inappropriate for such applications. These finite sample results also suggest that among three proposed methods $\mathrm{T}^{\tau,\mathrm{plug}}_\alpha$ is most aggressive and $\mathrm{T}^{\tau,\mathrm{ps}}_\alpha$ is most conservative.

These finite sample (with $n$ around several hundreds) results suggest that $\mathrm{T}^{\tau,\mathrm{plug}}_\alpha$ is useful only when $d$ is smaller than $n$. With $d$ approximates $n$, we recommend to use $\mathrm{T}^{\tau,\mathrm{jack}}_\alpha$ because it has averagely higher power. When $d$ is significantly larger than $n$, $\mathrm{T}^{\tau,\mathrm{ps}}_\alpha$ is recommended because of its good size control.

## 4.2. Real Data Example

In this section, we use proposed methods to analyze the dependence structure of brain activity. We use the resting-state functional magnetic resonance imaging (fMRI) data of normal children and diseased children with the disease attention deficit hyperactivity disorder (ADHD). Functional neuroimaging studies have revealed abnormalities in various brain regions of ADHD patients (Lou et al. (1990); Giedd et al. (2001); Shafritz et al. (2004); Yufeng et al. (2007); Zou et al. (2008)). As a marker of brain activity, amplitude of low-frequency fluctuation (ALFF) is a powerful tool to investigate this disorder. ALFF is the total power within the frequency range between 0.01 and 0.1 Hz of the fMRI time series. Generally speaking, it captures average slow fluctuations of brain activity. For the detailed definition of ALFF, we refer to Yufeng et al. (2007). Existing literature suggests the existence of significant differences in mean values of ALFF between the normal and diseased children (Zou et al. (2008)). By using our methods, we aim to test the dependence structure of ALFF between brain regions. Considering the nonlinear relationship and robustness, we use Kendall's tau matrix to measure the dependence structure.

We then introduce our data processing procedure. We use the standard methodology of software C-PAC[7] to correct body motion, brain heterogeneity, and many other kinds of measure errors. We then calculate voxel-wise ALFF of each person's fMRI images. As the voxel number is very large ($61 \times 73 \times 61$ for 3mm brain template), to limit the number of testing parameters, a common approach is to extract signals from specified regions of interest (ROIs) based on the anatomical structure of brain. In our experiments, we combine two kinds of brain areas including Brodmann (BA) and automated anatomical labeling (AAL) on the gray matter to build new 227 brain regions. In each brain region, we average obtained voxel-wise ALFF to get data points with the dimension $d = 227$.

After the introduction of data processing, we then describe the data set in detail. The resting-state fMRI data for ADHD is available on the Internet[8]. Considering the dimension (d=227) of data points, we use samples from Peking University and Kennedy Krieger Institute of Johns Hopkins University to build a sample with 119 ADHD patients and 200 control members.

**Table 3.** Region based two-sample tests of ALFF between ADHD patients and control members.

| | Mean vectcor | | | Kendall's tau matrix | | |
|---|---|---|---|---|---|---|
| | $T_\alpha^{\mathrm{Bai}}$ | $T_\alpha^{\mathrm{Sri}}$ | $T_\alpha^{\mathrm{Cai}}$ | $T_\alpha^{\mathrm{plug}}$ | $T_\alpha^{\mathrm{jack}}$ | $T_\alpha^{\mathrm{CLX}}$ |
| Test statistics | 2.8358 | 2.5978 | 22.7085 | 25.892 | 24.371 | 23.572 |
| P-values | 0.0023 | 0.0047 | 0.0006 | 0.0105 | 0.0223 | 0.0330 |

In the application, we test both mean vectors and Kendall's tau matrices between the diseased and normal groups and show the results in Table 3. In the context of high-dimensional mean tests, we use three existing methods: $T_\alpha^{\mathrm{Bai}}$ in Bai and Yin (1993), $T_\alpha^{\mathrm{Sri}}$ in Srivastava and Du (2008), and $\mathrm{T}_\alpha^{\mathrm{Cai}}$ in Cai et al. (2014). Except for the known mean differences, the results for Kendall's tau matrices also suggest that the dependence structure of brain activities for ADHD patients are also very different from normal children, which is worth investigating for related researchers.

[7]See the website http://fcp-indi.github.io/docs/user/.
[8]See the website http://fcon 1000.projects.nitrc.org/indi/adhd200/.

## 5. Summary and discussion

This paper considers the problem of testing the equality of high-dimensional U-statistic based matrices. We provide a lower bound for testing the equality of correlation matrices and prove the proposed methods' optimality. Based on thorough numerical comparisons, $\mathrm{T}_\alpha^{\mathrm{plug}}$ performs well only when $d$ is significantly smaller than $n$. When $d$ is very large, we recommend to use $\mathrm{T}_\alpha^{\mathrm{ps}}$ for correctly controlling the size . In addition, $\mathrm{T}_\alpha^{\mathrm{ps}}$ performs quite well for distributions with heavy tails or strong tail dependence. Therefore, $\mathrm{T}_\alpha^{\mathrm{ps}}$ is potentially more useful for financial applications in which heavy-tailness is a common phenomenon. There are many possible future directions of this work. For example, instead of two-sample problems, it is interesting to generalize the idea to $k$-sample testing problems ($k > 2$). This may require a nontrivial extension of theoretical analysis.

For testing Kendall's tau matrices, we show that the variance estimation error is a key factor influencing a test procedure's power. In fact, the test $\mathrm{T}_\alpha^{\mathrm{ps}}$, which exploits the exact value of variance under the uncorrelated condition ($\tau = 0$), achieves a better finite-sample performance especially when $d$ is very large. We can generalize such idea to many other applications. We also provide an upper bound of the Jackknife variance estimation error in the proof of Theorem 2.2. This result is also useful for other properties of U-statistics.

Next, we discuss the imposed assumptions. We note that the sparsity assumption **(A1)** plays a key role for obtaining the limiting extreme value distribution. It is not clear on whether this assumption is necessary, but it is satisfied in many high-dimensional applications. When **(A1)** is not satisfied, it is possible to exploit the bootstrap method to construct a test statistic. This is left as for future investigation. Regarding **(A2)**, we note that Cai et al. (2013) assume a stronger scaling assumption: $\log(d) = o(n^{1/5})$. We strengthen this scaling by assuming $\log(d) = O(n^{1/3-\epsilon})$ for an arbitrary $\epsilon > 0$. This is from the fact that U-statistics studied in this paper are assumed to have bounded kernels.

In the simulation studies, we use $\mathrm{T}_\alpha^{\mathrm{CLX}}$ as a comparison benchmark. In Supplement F of Supplement Material, we provide another heuristic test (denoted by $\mathrm{T}_\alpha^{\mathrm{R}}$) for testing the equality of Pearson's correlation matrices. The performances of $\mathrm{T}_\alpha^{\mathrm{CLX}}$ and $\mathrm{T}_\alpha^{\mathrm{R}}$ are similar for off diagonal disturbances.

### Supplementary Material

**Technical Proofs and More Simulation for "An Extreme-Value Approach for Testing the Equality of Large U-Statistic Based Correlation Matrices"**
(doi: COMPLETED BY THE TYPESETTER; supplement.pdf). We provide additional proof and simulation in Supplementary Material. Supplementary Material consists of 6 parts: Supplements A-F. Among them, Supplements A-D prove the theorems that are not proven in Appendix A. Supplement E introduces some useful definitions. Supplement F presents more simulation results.

# Acknowledgement

The authors are grateful for the support of NSF CAREER Award DMS1454377, NSF IIS1408910, NSF IIS1332109, NIH R01MH102339, NIH R01GM083084, and NIH R01HG06841. Fang's research is supported by a Google research fellowship and NSF DMS-1712536. Xinsheng's research is supported by the National Natural Science Foundation of China (Grant No. 11071045). Cheng's research is also supported by a fellowship from CSC (China Scholarship Council).

# Appendix A The Proof of main theorem

This appendix contains the proof of main theorem, i.e., Theorem 2.2. In the sequel, we use $C$, $C_1$, $C_2$, ..., to denote constants that do not depend on $n$, $d$, $q$ and they can vary from place to place.

***Proof.*** As explained in the sketch of proof, our analysis proceeds in three steps.

**Step (i)**. In this step, we prove that it is sufficient to establish (2.14) for proving the theorem. For this, we need to sharply characterize the estimation error of the Jackknife variance estimator of U-statistics. For this, we introduce the following lemma.

**Lemma A.1.** Let $\widehat{\sigma}^2(\widehat{u}_{a,ij})$ be the Jackknife estimator of $\widehat{u}_{a,ij}$ and $\sigma^2(\widehat{u}_{a,ij})$ be the variance of $\widehat{u}_{a,ij}$. Recalling the definition of $h_{ij}$ and $\zeta_{a,ij}$ in (2.9) and (2.10), $\zeta_{1,ij}$ and $\zeta_{2,ij}$ are the variances of $h_{ij}(\boldsymbol{X}_\ell)$ and $h_{ij}(\boldsymbol{Y}_\ell)$. We have that $m^2\zeta_{a,ij}$ is the limit of $n_a\sigma^2(\widehat{u}_{a,ij})$ as $n_a$ goes to infinity. We also have that $m^2\zeta_{a,ij}$ is the limit of $n_a\widehat{\sigma}^2(\widehat{u}_{a,ij})$ as $n_a$ goes to infinity. Moreover, under Assumption (**A2**), as $n$, $q \to \infty$ we have

$$\mathbb{P}\Big(\max_{1\le i,j\le q}\big|n_a\widehat{\sigma}^2(\widehat{u}_{a,ij}) - m^2\zeta_{a,ij}\big| \ge C\frac{\varepsilon_n}{\log q}\Big) = o(1), \tag{A.1}$$

where $\varepsilon_n = o(1)$ and $a = 1, 2$.

The detailed proof of Lemma A.1 is in Supplement B.1 of Supplementary Material. This Lemma presents an upper bound of Jackknife variance estimation error, which enables us to obtain the convergence rate of Jackknife variance estimator. To prove this lemma, we decompose $\widehat{\sigma}^2(\widehat{u}_{a,ij})$ into different pieces and bound each piece separately. The details of this decomposition are in Supplements B.1 and C.1. Both the result and the proof of Lemma A.1 are nontrivial and are of independent technical interest.

Lemma A.1 implies that both of the following two events

$$\mathcal{E}_1 := \Big\{ \max_{1\le i,j\le q} \big| n_1 \widehat{\sigma}^2(\widehat{u}_{1,ij}) - m^2\zeta_{1,ij} \big| < C\frac{\varepsilon_n}{\log q} \Big\},$$
$$\mathcal{E}_2 := \Big\{ \max_{1\le i,j\le q} \big| n_2 \widehat{\sigma}^2(\widehat{u}_{2,ij}) - m^2\zeta_{2,ij} \big| < C\frac{\varepsilon_n}{\log q} \Big\},$$

happen with probability going to one as $n, q \to \infty$. Under $\mathcal{E}_1$ and $\mathcal{E}_2$, by $\zeta_{a,ij} \ge r_a > 0$ (Assumption **(A2)**), we have

$$\Big| n_1\widehat{\sigma}^2(u_{1,ij})/(m^2\zeta_{1,ij}) - 1 \Big| < C\varepsilon_n/\log q \quad \text{and} \quad \Big| n_2\widehat{\sigma}^2(u_{2,ij})/(m^2\zeta_{2,ij}) - 1 \Big| < C\varepsilon_n/\log q.$$

We set $M_n := \max_{1\le i,j\le q} M_{ij}$ and $\widetilde{M}_n := \max_{1\le i,j\le q} \widetilde{M}_{ij}$. By the definition of $M_n$ and $\widetilde{M}_n$, we calculate the relative difference of $M_{ij}$ and $\widetilde{M}_{ij}$ as

$$\Big| \frac{M_{ij} - \widetilde{M}_{ij}}{\widetilde{M}_{ij}} \Big| \le \Big| \frac{\widehat{\sigma}^2(\widehat{u}_{1,ij}) - m^2\zeta_{1,ij}/n_1}{\widehat{\sigma}^2(\widehat{u}_{1,ij})} \Big| + \Big| \frac{\widehat{\sigma}^2(\widehat{u}_{2,ij}) - m^2\zeta_{2,ij}/n_2}{\widehat{\sigma}^2(\widehat{u}_{2,ij})} \Big| \le C\frac{\varepsilon_n}{\log q}. \tag{A.2}$$

Therefore, we have $|M_{ij} - \widetilde{M}_{ij}| \le C\varepsilon_n \widetilde{M}_{ij}/\log q$, which implies that

$$|M_n - \widetilde{M}_n| \le \max_{1\le i,j\le n_1} |M_{ij} - \widetilde{M}_{ij}| \le C\widetilde{M}_n\varepsilon_n/\log q. \tag{A.3}$$

Combining $\widetilde{M}_n/\log q = O_p(1)$ and $\varepsilon_n = o(1)$, to prove Theorem 2.2 it suffices to show that as $n, q \to \infty$, (2.14) holds for any $x \in \mathbb{R}$.

**Step (ii)**. In this step, we use the Hoeffding decomposition (Lemma D.4 in Supplementary Material) to decompose U-statistics. We then prove the residual term $\Delta_{n_a,ij}/\binom{n_a}{m}$ is negligible, i.e., to prove the theorem it is sufficient to prove (2.17) as $n$, $q \to \infty$.

For notational simplicity, we set

$$\widetilde{N}_{ij} := (\widehat{u}_{1,ij} - \widehat{u}_{2,ij})/\sqrt{m^2\zeta_{1,ij}/n_1 + m^2\zeta_{2,ij}/n_2}. \tag{A.4}$$

Recall that in (2.13) and (2.16) we define $\widetilde{M}_{ij}$ and $T_{ij}$ as

$$\widetilde{M}_{ij} := \frac{(\widehat{u}_{1,ij} - \widehat{u}_{2,ij})^2}{m^2\zeta_{1,ij}/n_1 + m^2\zeta_{2,ij}/n_2},\ T_{ij} := \frac{\sum_{\alpha=1}^{n_1} h_{ij}(\boldsymbol{X}_\alpha)/n_1 - \sum_{\alpha=1}^{n_2} h_{ij}(\boldsymbol{Y}_\alpha)/n_2}{\sqrt{\zeta_{1,ij}/n_1 + \zeta_{2,ij}/n_2}}. \tag{A.5}$$

By the definition of $\widetilde{M}_{ij}$, we have $\widetilde{M}_{ij} = (\widetilde{N}_{ij})^2$. Combining the definition of $T_{ij}$ and (2.15), we have

$$\widetilde{N}_{ij} = T_{ij} + \frac{\binom{n_1}{m}^{-1}\Delta_{n_1,ij} - \binom{n_2}{m}^{-1}\Delta_{n_2,ij}}{\sqrt{m^2\zeta_{1,ij}/n_1 + m^2\zeta_{2,ij}/n_2}}. \tag{A.6}$$

We then introduce the following lemma to analyze the difference of $\widetilde{N}_{ij}$ and $T_{ij}$.

**Lemma A.2.** As $n$, $q \to \infty$, we have

$$\Big| \max_{1\le i,j\le q} (\widetilde{N}_{ij})^2 - \max_{1\le i,j\le q} (T_{ij})^2 \Big| = o_p(1). \tag{A.7}$$

The detailed proof of Lemma A.2 is in Supplement B.2. This lemma illustrates that $\max_{1\le,i,j\le q} \widetilde{N}_{ij}$ and $T_n := \max_{1\le i,j\le q} T_{ij}$ have the same limiting distribution. Hence, to prove Theorem 2.2 it suffices to show (2.17) as $n$, $q \to \infty$.

**Step (iii)**. In this step, we aim to prove (2.17). In (2.17), $T_n$ is the maximum of $T_{ij}$ over $S := \{(i,j) : 1 \le i,j \le q\}$ and these $T_{ij}$'s are not independent of each other. Therefore, we cannot straightforwardly exploit the extreme value theorem under the independent setting to obtain the limiting distribution of $T_n$. To solve this problem, we construct normal approximation to obtain the extreme value distribution of $(T_{ij})_{1\le i,j\le q}$ under the setting that $T_{ij}$ can be dependent of each other. The construction of such normal approximation requires most correlations of different $T_{ij}$ to be small. Correlations between different $T_{ij}$'s are related to the correlations of entries of $\boldsymbol{X}$ and $\boldsymbol{Y}$. Assumption **(A1)** specifies sufficient conditions on the correlations of entries of $\boldsymbol{X}$ and $\boldsymbol{Y}$.

To obtain more insight of Assumption **(A1)**, we introduce the following notations. We use $S_0$ to denote pairs of $(i,j)$ such that $X_i$ and $X_j$ are highly correlated ($|u_{1,ij}| > (\log q)^{-1-\alpha_0}$) or $Y_i$ and $Y_j$ are highly correlated ($|u_{2,ij}| > (\log q)^{-1-\alpha_0}$). Recalling the formal definition of $S_0$ in (2.11), Assumption **(A1)** implies that the number of highly correlated ($|u_{a,ij}| > (\log q)^{-1-\alpha_0}$) entries of $\boldsymbol{X}$ and $\boldsymbol{Y}$ is small. More specifically, Assumption **(A1)** assumes $|S_0| = o(q^2)$.

We can prove that correlations between $T_{ij}$'s on $S \setminus S_0$ are all small. We then use the Bofferroni inequality (Lemma 1 of Cai et al. (2013)) and normal approximation to obtain the limiting distribution of $\max_{(i,j)\in S\setminus S_0}(T_{ij})^2$ so as to prove (2.17).

We then present the detailed proof of (2.17). Firstly, we prove that it suffices to take the maximum of $T_{ij}$ over $S \setminus S_0$ but instead of over $S$ as in (2.17). By setting $y_q = x + 4\log q - \log(\log q)$, we have

$$\Big|\mathbb{P}\Big(\max_{(i,j)\in S}(T_{ij})^2 \ge y_q\Big) - \mathbb{P}\Big(\max_{(i,j)\in S\setminus S_0}(T_{ij})^2 \ge y_q\Big)\Big| \le \mathbb{P}\Big(\max_{(i,j)\in S_0}(T_{ij})^2 \ge y_q\Big). \tag{A.8}$$

The next lemma implies that, as $n$, $q \to \infty$, we have $\mathbb{P}\big(\max_{(i,j)\in S_0}(T_{ij})^2 \ge y_q\big) \to 0$.

**Lemma A.3.** Under Assumptions **(A1)** and **(A2)**, as $n$, $q \to \infty$, we have

$$\mathbb{P}\Big(\max_{(i,j)\in S_0}(T_{ij})^2 \ge y_q\Big) \to 0.$$

The detailed proof of Lemma A.3 is in Supplement B.3 of Supplementary Material. By Lemma A.3, we have $\mathbb{P}\big(\max_{(i,j)\in S_0}(T_{ij})^2 \ge y_q\big) \to 0$ as $n, q \to \infty$. Moreover, by (A.8), we have that $\mathbb{P}\big(\max_{(i,j)\in S}(T_{ij})^2 \ge y_q\big)$ and $\mathbb{P}\big(\max_{(i,j)\in S\setminus S_0}(T_{ij})^2 \ge y_q\big)$ have the same limit value as $n$, $q \to \infty$. Therefore, to obtain (2.17), it suffices to prove

$$\mathbb{P}\Big(\max_{(i,j)\in S\setminus S_0}(T_{ij})^2 - 4\log q + \log(\log q) \le x\Big) \to \exp\Big(-\frac{1}{\sqrt{8\pi}}\exp(-\frac{x}{2})\Big), \tag{A.9}$$

as $n$, $q \to \infty$. The problem is then reduced to prove (A.9).

For simplicity, by rearranging the two-dimensional indices $\{(i,j) : (i,j) \in S \setminus S_0\}$ in any order, we set them as $\{(i_k, j_k) : 1 \le k \le h\}$ with $h = |S \setminus S_0|$. If we denote $T_k := T_{i_k j_k}$, (A.9) becomes

$$\mathbb{P}\Big( \max_{1 \le k \le h} (T_k)^2 - 4 \log q + \log(\log q) \le x \Big) \to \exp\Big( - \frac{1}{\sqrt{8\pi}} \exp(-\frac{x}{2}) \Big). \tag{A.10}$$

Secondly, we exploit normal approximation to obtain the limiting distribution of $\max_{1 \le k \le h}(T_k)^2$. This normal approximation is useful for getting the extreme value distribution of weakly dependent data. By excluding all the pairs in $S_0$, correlations between $T_k$'s are all small. Therefore, we can use this normal approximation to get the limiting distribution of $\max_{1 \le k \le h}(T_k)^2$. In detail, we first use the Boferroni inequality to obtain both lower and upper bounds of $\mathbb{P}\big(\max_{1 \le k \le h}(T_k)^2 \ge y_q\big)$. The obtained lower and upper bounds can then be shown to have the same limiting distribution, which is the extreme value distribution with the cumulative distribution function of $\exp\big( - (8\pi)^{-1/2} \exp(-x/2)\big)$.

To describe the procedure of normal approximation, we need some additional notations. We introduce

$$\begin{cases} \widehat{Z}_{\beta,ij} = n_2 h_{ij}(\boldsymbol{X}_\beta)/n_1 & \text{for} \quad 1 \le \beta \le n_1, \\ \widehat{Z}_{\beta,ij} = -h_{ij}(\boldsymbol{Y}_{\beta - n_1}) & \text{for} \quad n_1 + 1 \le \beta \le n_1 + n_2, \end{cases} \tag{A.11}$$

where $h_{ij}$ is defined in (2.9). Moreover, by the definition of $T_{ij}$ in (2.16), we have

$$T_k := T_{i_k j_k} = \sum_{\beta=1}^{n_1+n_2} \widehat{Z}_{\beta, i_k j_k} \Big/ \sqrt{n_2^2 \zeta_{1,i_k j_k}/n_1 + n_2 \zeta_{2,i_k j_k}}. \tag{A.12}$$

After introducing these notations, we explain how to use normal approximation to get the extreme value distribution of $\max_{1 \le k \le h}(T_k)^2$. Firstly, by the Boferroni inequality (Lemma 1 of Cai et al. (2013)), for any integer $M$ with $0 < M < [h/2]$, we have

$$\begin{aligned} \sum_{\ell=1}^{2M} (-1)^{\ell-1} \sum_{1 \le k_1 < \cdots < k_\ell \le h} \mathbb{P}\big( \bigcap_{j=1}^{\ell} E_{k_j} \big) \le & \mathbb{P}\big( \max_{1 \le k \le h} (T_k)^2 \ge y_q \big) \\ & \le \sum_{\ell=1}^{2M-1} (-1)^{\ell-1} \sum_{1 \le k_1 < \cdots < k_\ell \le h} \mathbb{P}\big( \bigcap_{j=1}^{\ell} E_{k_j} \big), \end{aligned} \tag{A.13}$$

where we set $E_{k_j} = \{(T_{k_j})^2 \ge y_q\}$. In next step, to simplify $\mathbb{P}\big( \bigcap_{j=1}^{\ell} E_{k_j} \big)$, we define

$$\widetilde{Z}_{\beta k} = \widehat{Z}_{\beta, i_k j_k} / (n_2 \zeta_{1,i_k j_k}/n_1 + \zeta_{2,i_k j_k})^{1/2} \quad \text{and} \quad \boldsymbol{W}_\beta = (\widetilde{Z}_{\beta k_1}, \ldots, \widetilde{Z}_{\beta k_\ell})^T, \tag{A.14}$$

for $1 \le k \le h$ and $1 \le \beta \le n_1+n_2$. Therefore, we have $T_{k_j} = (n_2)^{-1}\sum_{\beta=1}^{n_1+n_2}\widetilde{Z}_{\beta k_j}$. Define $\|\mathbf{v}\|_{\min} = \min_{1\le i\le \ell}|v_i|$ for vector $\mathbf{v}\in\mathbb{R}^\ell$. With these notations, we rewrite $\mathbb{P}\big(\bigcap_{j=1}^{\ell} E_{k_j}\big)$ as

$$\mathbb{P}\Big(\bigcap_{j=1}^{\ell} E_{k_j}\Big) = \mathbb{P}\Big(\Big\|n_2^{-1/2}\sum_{\beta=1}^{n_1+n_2}\boldsymbol{W}_\beta\Big\|_{\min} \ge y_q^{1/2}\Big).$$

Secondly, we use a normal vector $\boldsymbol{N}_\ell$ to approximate $n_2^{-1/2}\sum_{\beta=1}^{n_1+n_2}\boldsymbol{W}_\beta$. In detail, we set $\boldsymbol{N}_\ell$ as a normal vector with the same mean vector and the same covariance matrix as $n_2^{-1/2}\sum_{\beta=1}^{n_1+n_2}\boldsymbol{W}_\beta$. More specifically, we have

$$\boldsymbol{N}_\ell := (N_{k_1},\dots,N_{k_\ell})^T \text{ with } \mathbb{E}[\boldsymbol{N}_\ell]=0,\ \mathrm{Var}(\boldsymbol{N}_\ell)=n_1\mathrm{Var}(\boldsymbol{W}_1)/n_2+\mathrm{Var}(\boldsymbol{W}_{n_1+1}). \tag{A.15}$$

The following lemma uses $\boldsymbol{N}_\ell$ to rewrite the the upper and lower bounds in (A.13).

**Lemma A.4.** Under Assumption **(A2)**, as $n,\ q\to\infty$, we have

$$\mathbb{P}\Big(\max_{1\le k\le h}(T_k)^2 \ge y_q\Big) \le \sum_{\ell=1}^{2M-1}(-1)^{\ell-1}\sum_{1\le k_1<\cdots<k_\ell\le h}\mathbb{P}\Big(\|\boldsymbol{N}_\ell\|_{\min} \ge y_q^{1/2}-\epsilon_n(\log q)^{-1/2}\Big)+o(1), \tag{A.16}$$

$$\mathbb{P}\Big(\max_{1\le k\le h}(T_k)^2 \ge y_q\Big) \ge \sum_{\ell=1}^{2M}(-1)^{\ell-1}\sum_{1\le k_1<\cdots<k_\ell\le h}\mathbb{P}\Big(\|\boldsymbol{N}_\ell\|_{\min} \ge y_q^{1/2}+\epsilon_n(\log q)^{-1/2}\Big)-o(1). \tag{A.17}$$

The detailed proof of Lemma A.4 is in Supplement B.4 of Supplementary Material. At last, to complete the proof, we need to prove that the right hand sides of (A.16) and (A.17) have the same limit value $1-\exp\big(-(\sqrt{8\pi})^{-1}\exp(-x/2)\big)$ as $n,\ q\to\infty$. To calculate the limit value, we need the following lemma.

**Lemma A.5.** Under Assumption **(A3)**, for any integer $\ell \ge 1$ and $x\in\mathbb{R}$, we have

$$\sum_{1\le k_1<\ldots<k_\ell\le h}\mathbb{P}\Big(\|\boldsymbol{N}_\ell\|_{\min} \ge y_q^{1/2}\pm\epsilon_n(\log q)^{-1/2}\Big) = \frac{1}{\ell!}\Big(\frac{1}{\sqrt{8\pi}}\exp(-\frac{x}{2})\Big)^\ell(1+o(1)). \tag{A.18}$$

The detailed proof of Lemma A.5 is in Supplement B.5. By plugging (A.18) into (A.16) and (A.17), we construct the following inequities:

$$\limsup_{n,q\to\infty}\mathbb{P}\Big(\max_{1\le k\le h}(T_k)^2 \ge y_q\Big) \le \sum_{\ell=1}^{2M-1}(-1)^{\ell-1}\frac{1}{\ell!}\Big(\frac{1}{\sqrt{8\pi}}\exp(-\frac{x}{2})\Big)^\ell,$$
$$\liminf_{n,q\to\infty}\mathbb{P}\Big(\max_{1\le k\le h}(T_k)^2 \ge y_q\Big) \ge \sum_{\ell=1}^{2M}(-1)^{\ell-1}\frac{1}{\ell!}\Big(\frac{1}{\sqrt{8\pi}}\exp(-\frac{x}{2})\Big)^\ell,$$

for any positive integer $M$. Letting $M\to\infty$, we prove (A.10). Therefore, we finish the proof of Theorem 2.2. □

# Supplementary Material: Technical Proofs and More Simulation for "An Extreme-Value Approach for Testing the Equality of Large U-Statistic Based Correlation Matrices"

Cheng Zhou*, Fang Han**, Xinsheng Zhang†, and Han Liu‡

This document contains the additional details of the paper "An Extreme-Value Approach for Testing the Equality of Large U-Statistic Based Correlation Matrices" authored by Cheng Zhou, Fang Han, Xinsheng Zhang and Han Liu. It is organized as follows. Supplement A contains detailed proofs of theorems that are not proven in Appendix A, including Theorems 2.4, 2.6, 3.2, 3.3, 3.4, 3.5, 3.6 and Remark 3.8. Supplements B and C prove the introduced lemmas during the proof. Supplement D contains additional technical lemmas. Supplement E contains the definition of meta-elliptical distribution and its properties. Supplement F contains more simulation results.

## Supplement A: Detailed Proofs of Theorems

This appendix contains detailed proofs of theorems that are not proven in Appendix A, including Theorems 2.4, 2.6, 3.2, 3.3, 3.4, 3.5, 3.6 and Remark 3.8.

### A.1. Proof of Theorem 2.4

***Proof.*** In Theorem 2.4, we aim to prove that, under Assumption **(A2)**, as $n, q \to \infty$, we have

$$\inf_{(\mathbf{U}_1,\mathbf{U}_2)\in\mathbb{A}(4)} \mathbb{P}(M_n \geq G^{-}(\alpha) + 4\log q - \log(\log q)) \to 1, \tag{A.1}$$

where $G^{-}(\alpha) := -\log(8\pi) - 2\log\big(-\log(1-\alpha)\big)$. In (A.1), we set $\mathbf{U}_a$ and $\mathbb{A}(C)$ as $\mathbf{U}_a = (u_{a,ij}) \in \mathbb{R}^{q\times q}$ and

$$\mathbb{A}(C) = \left\{(\mathbf{U}_1,\mathbf{U}_2) : \max_{(i,j)\in S} \frac{|u_{1,ij} - u_{2,ij}|}{\sqrt{m^2\zeta_{1,ij}/n_1 + m^2\zeta_{2,ij}/n_2}} \geq C\sqrt{\log q}\right\},$$

where $S = \{(i,j) : 1 \leq i,j \leq q\}$. We use $\mathbb{A}(C)$ to characterize the set of alternative hypotheses for analyzing the power of $\mathrm{T}_\alpha$ in (2.6).

If one entry of $\mathbf{U}_1 - \mathbf{U}_2$ has a magnitude large enough (for example $\gg \sqrt{\log q/n}$), we have $(\mathbf{U}_1, \mathbf{U}_2) \in \mathbb{A}(C)$. Therefore, large perturbations on a few entries compared to the null hypothesis $\mathbf{U}_1 = \mathbf{U}_2$ can easily make the pair $(\mathbf{U}_1, \mathbf{U}_2)$ belong to $\mathbb{A}(C)$ for some constant $C$. In Theorem 2.4, we require that $(\mathbf{U}_1, \mathbf{U}_2) \in \mathbb{A}(4)$, which implies

$$\max_{(i,j)\in S} \frac{(u_{1,ij} - u_{2,ij})^2}{m^2\zeta_{1,ij}/n_1 + m^2\zeta_{2,ij}/n_2} \geq 16 \log q. \tag{A.2}$$

Under the alternative hypothesis, $u_{1,ij} = u_{2,ij}$ cannot hold for all $(i,j) \in S$. This motivates us to define

$$\begin{aligned} M_n^1 &:= \max_{(i,j)\in S} \frac{(\widehat{u}_{1,ij} - \widehat{u}_{2,ij} - u_{1,ij} + u_{2,ij})^2}{\widehat{\sigma}^2(\widehat{u}_{1,ij}) + \widehat{\sigma}^2(\widehat{u}_{2,ij})} \\ M_n &:= \max_{(i,j)\in S} M_{ij} = \max_{(i,j)\in S} \frac{(\widehat{u}_{1,ij} - \widehat{u}_{2,ij})^2}{\widehat{\sigma}^2(\widehat{u}_{1,ij}) + \widehat{\sigma}^2(\widehat{u}_{2,ij})}, \end{aligned} \tag{A.3}$$

because $M_n^1$ and $M_n$ are different under the alternative hypothesis.

To prove (2.17), using the inequality $(a \pm b)^2 \leq 2a^2 + 2b^2$, we get

$$(u_{1,ij} - u_{2,ij})^2 \leq 2(\widehat{u}_{1,ij} - \widehat{u}_{2,ij} - u_{1,ij} + u_{2,ij})^2 + 2(\widehat{u}_{1,ij} - \widehat{u}_{2,ij})^2.$$

Therefore, by the the definitions of $M_n$ and $M_n^1$ in (A.3), we get

$$\max_{(i,j)\in S} \frac{(u_{1,ij} - u_{2,ij})^2}{\widehat{\sigma}^2(\widehat{u}_{1,ij}) + \widehat{\sigma}^2(\widehat{u}_{2,ij})} \leq 2M_n^1 + 2M_n. \tag{A.4}$$

Under $(\mathbf{U}_1, \mathbf{U}_2) \in \mathbb{A}(4)$, to prove $\mathbb{P}\big(M_n \geq G^-(\alpha) + 4\log q - \log(\log q)\big) \to 1$, we need the following two lemmas.

**Lemma A.1.** Under Assumption **(A2)**, as $n, q \to \infty$, we have

$$\mathbb{P}\Big(M_n^1 \leq 4\log q - \frac{1}{2}\log(\log q)\Big) \to 1. \tag{A.5}$$

The detailed proof of Lemma A.1 is in Supplement B.6.

**Lemma A.2.** Under Assumption **(A2)**, as $n$, $q \to \infty$, we have

$$\mathbb{P}\Big(\max_{(i,j)\in S} \frac{(u_{1,ij} - u_{2,ij})^2}{\widehat{\sigma}^2(\widehat{u}_{1,ij}) + \widehat{\sigma}^2(\widehat{u}_{2,ij})} \geq 16\log q\Big) \to 1. \tag{A.6}$$

uniformly for all $(\mathbf{U}_1, \mathbf{U}_2) \in \mathbb{A}(4)$.

The detailed proof of Lemma A.2 is in Supplement B.7.

Combining (A.4), (A.5) and (A.6), as $n, q \to \infty$, with probability going to one, we have

$$M_n \geq \frac{1}{2} \max_{(i,j)\in S} \frac{(u_{1,ij} - u_{2,ij})^2}{\widehat{\sigma}^2(\widehat{u}_{1,ij}) + \widehat{\sigma}^2(\widehat{u}_{2,ij})} - M_n^1 \geq 4\log q + \frac{1}{2}\log(\log q).$$

Therefore, as $n, q \to \infty$, we have

$$1 \geq \mathbb{P}\big(M_n \geq G^-(\alpha) + 4\log q - \log(\log q)\big) \geq \mathbb{P}\Big(M_n \geq 4\log q + \frac{1}{2}\log(\log q)\Big) \to 1.$$

Hence, we finish the proof of Theorem 2.4. □

## A.2. Proof of Theorem 3.2

***Proof.*** In Theorem 3.2, we aim to prove (3.10) and (3.11). (3.10) is an application of Theorem 2.2 for Kendall's tau matrices. Therefore, we only need to prove (3.11). (3.11) is the same as (3.10) except that in (3.11) we use $\widehat{\sigma}^2_{\text{plug}}(\widehat{\tau}_{a,ij})$ (defined in (3.5)) to replace the Jackknife variance estimator $\widehat{\sigma}^2(\widehat{\tau}_{a,ij})$ in (3.10).

We prove (3.11) in the same three steps as those for Theorem 2.2 except that we need to verify whether $\widehat{\sigma}^2_{\text{plug}}(\widehat{\tau}_{a,ij})$ satisfies the following equation:

$$\mathbb{P}\Big(\max_{1\le i<j\le d}|n_a\widehat{\sigma}^2_{\text{plug}}(\widehat{\tau}_{a,ij}) - 4\zeta_{a.ij}| > C\frac{\varepsilon_n}{\log d}\Big) = o(1), \tag{A.7}$$

where $\varepsilon_n = o(1)$. (A.7) is the same as (A.1) except that in (A.7) we use $\widehat{\sigma}^2_{\text{plug}}(\widehat{\tau}_{a,ij})$ to replace $\widehat{\sigma}^2(\widehat{\tau}_{a,ij})$ in (A.1). (A.7) presents an upper bound of plug-in variance estimation error, which enables us to obtain the convergence rate of plug-in variance estimator. Both the result and the proof of (A.7) are nontrivial and are of independent technical interest. Once (A.7) is obtained, by the same argument in the first step of the proof of Theorem 2.2, we obtain $|M_n^{\tau,\text{plug}} - \max_{(i,j)\in S}\widetilde{M}_{ij}| = o_p(1)$, where $M_n^{\tau,\text{plug}}$ and $\widetilde{M}_{ij}$ are defined in (3.6) and (2.13). Following the steps of the proof for Theorem 2.2, we obtain (3.11).

The proof of (A.7) proceeds in two steps. We define $\Pi_{cc,ij}$ and $\Pi_{c,ij}$ in (3.4) and (3.2). In the first step, we rewrite (A.7) as

$$\mathbb{P}\Big(\max_{1\le i<j\le d}\Big|\Big((\widehat{\Pi}_{cc,ij} - (\widehat{\Pi}_{c,ij})^2) - (\Pi_{cc,ij} - (\Pi_{c,ij})^2)\Big)\Big| > C\frac{\varepsilon_n}{\log d}\Big) = o(1), \tag{A.8}$$

where $\widehat{\Pi}_{cc,ij}$ and $\widehat{\Pi}_{c,ij}$ are U-statistics and they estimate $\Pi_{cc,ij}$ and $\Pi_{c,ij}$. In the second step, we obtain upper bounds of $\max_{1\le i<j\le d}|\widehat{\Pi}_{cc,ij} - \Pi_{cc,ij}|$ and $\max_{1\le i<j\le d}|(\widehat{\Pi}_{c,ij})^2 - (\Pi_{c,ij})^2|$. Using the obtained upper bounds, we prove (A.8).

**Step (i)**. For notational simplicity, we ignore the subscript $a$ ($a = 1$ or 2) for $\boldsymbol{X}$ and $\boldsymbol{Y}$. To rewrite (A.7) as (A.8), we need $4\zeta_{ij} = 16(\Pi_{cc,ij} - (\Pi_{c,ij})^2)$ and $\widehat{\sigma}^2_{\text{plug}}(\widehat{\tau}_{ij}) = 16(\widehat{\Pi}_{cc,ij} - \widehat{\Pi}^2_{c,ij})/n$. By the definition of $\widehat{\sigma}^2_{\text{plug}}(\widehat{\tau}_{ij})$ in (3.5), we have $\widehat{\sigma}^2_{\text{plug}}(\widehat{\tau}_{ij}) = 16(\widehat{\Pi}_{cc,ij} - \widehat{\Pi}^2_{c,ij})/n$. We then need to prove $4\zeta_{ij} = 16(\Pi_{cc,ij} - (\Pi_{c,ij})^2)$. By (3.3), we have that the variance of $\widehat{\tau}_{ij}$ is

$$\frac{8}{n(n-1)}\Pi_{c,ij}(1-\Pi_{c,ij}) + 16\frac{1}{n}\frac{n-2}{n-1}(\Pi_{cc,ij} - \Pi^2_{c,ij}).$$

Therefore, $16(\Pi_{cc,ij} - (\Pi_{c,ij})^2)$ is the asymptotic variance of $\sqrt{n}\widehat{\tau}_{ij}$ as $n\to\infty$. Lemma A.1 implies that $4\zeta_{ij}$ is also the asymptotic variance of $\sqrt{n}\widehat{\tau}_{ij}$ as $n\to\infty$. Hence, we have

$$16(\Pi_{cc,ij} - (\Pi_{c,ij})^2) = 4\zeta_{ij}, \tag{A.9}$$

because both sides of (A.9) are the asymptotic variance of $\sqrt{n}\widehat{\tau}_{ij}$. Therefore, combining $\widehat{\sigma}^2_{\text{plug}}(\widehat{\tau}_{ij}) = 16(\widehat{\Pi}_{cc,ij} - \widehat{\Pi}^2_{c,ij})/n$ and (A.9), to prove (A.7), it suffices to prove (A.8) as $n,\ q\to\infty$.

**Step (ii)**. We aim to prove (A.8), in which $\widehat{\Pi}_{cc,ij}$ is a U-statistic estimator of $\Pi_{cc,ij}$ (defined in (3.4)). Therefore, we use the exponential inequality (Lemma D.2) to obtain

$$\mathbb{P}\Big(\max_{1\le i<j\le d}|\widehat{\Pi}_{cc,ij}-\Pi_{cc,ij}|>t\Big)\le C_1d^2\exp(-C_2nt^2). \tag{A.10}$$

$\widehat{\Pi}_{c,ij}$ is also a U-statistic which estimates $\Pi_{c,ij}$ (defined in (3.2)). Similarly to (A.10), we use the exponential inequality (Lemma D.2) to obtain

$$\mathbb{P}\Big(\max_{1\le i<j\le d}|\widehat{\Pi}_{c,ij}-\Pi_{c,ij}|>t\Big)\le C_1d^2\exp(-C_2nt^2). \tag{A.11}$$

By the definitions of $\widehat{\Pi}_{c,ij}$ and $\Pi_{c,ij}$, we have $|\widehat{\Pi}_{c,ij}|\le 2$ and $|\Pi_{c,ij}|\le 2$. Therefore, we have $|\widehat{\Pi}_{c,ij}+\Pi_{c,ij}|\le 4$ which combined with (A.11), implies that

$$\begin{aligned}\mathbb{P}\Big(\max_{1\le i<j\le d}|(\widehat{\Pi}_{c,ij})^2-(\Pi_{c,ij})^2|>t\Big)&\le\mathbb{P}\Big(\max_{1\le i<j\le d}4|\widehat{\Pi}_{c,ij}-\Pi_{c,ij}|>t\Big)\\&\le C_1d^2\exp(-C_2nt^2).\end{aligned} \tag{A.12}$$

Considering $\log d=O(n^{1/3-\epsilon})$ (see Assumption **(A2)**), (A.10) and (A.12), by setting $\varepsilon_n=1/(\log d)^{\kappa_0}$ with $\kappa_0>0$ sufficiently small, we have

$$\begin{aligned}\mathbb{P}\Big(\max_{1\le i<j\le d}\big|\widehat{\Pi}_{cc,ij}-\Pi_{cc,ij}\big|>C\frac{\varepsilon_n}{\log d}\Big)&=o(1),\\\mathbb{P}\Big(\max_{1\le i<j\le d}\big|(\widehat{\Pi}_{c,ij})^2-(\Pi_{c,ij})^2\big|>C\frac{\varepsilon_n}{\log d}\Big)&=o(1).\end{aligned} \tag{A.13}$$

By the triangle inequality, we have

$$\begin{aligned}&\mathbb{P}\Big(\max_{1\le i<j\le d}\Big|\Big(\big(\widehat{\Pi}_{cc,ij}-(\widehat{\Pi}_{c,ij})^2\big)-\big(\Pi_{cc,ij}-(\Pi_{c,ij})^2\big)\Big)\Big|>C\tfrac{\varepsilon_n}{\log d}\Big)\\&\le\ \mathbb{P}\Big(\max_{1\le i<j\le d}\big|\widehat{\Pi}_{cc,ij}-\Pi_{cc,ij}\big|>0.5C\tfrac{\varepsilon_n}{\log d}\Big)\\&\quad+\mathbb{P}\Big(\max_{1\le i<j\le d}\big|(\widehat{\Pi}_{c,ij})^2-(\Pi_{c,ij})^2\big|>0.5C\tfrac{\varepsilon_n}{\log d}\Big).\end{aligned} \tag{A.14}$$

Combining (A.13) and (A.14), we prove (A.8). This completes the proof. □

## A.3. Proof of Theorem 3.3

***Proof.*** By setting

$$M_{ij}^{\tau,\mathrm{ps}}:=\frac{(\widehat{\tau}_{1,ij}-\widehat{\tau}_{2,ij})^2}{\sigma_{1,\mathrm{ps}}^2/n_1+\sigma_{2,\mathrm{ps}}^2/n_2}\qquad\text{and}\qquad M_n^{\tau,\mathrm{ps}}:=\max_{(i,j)\in S}M_{ij}^{\tau,\mathrm{ps}}, \tag{A.15}$$

where $S=\{(i,j):1\le i<j\le d\}$, we aim to prove that as $n,\ d\to\infty$, we have

$$\mathbb{P}\big(M_n^{\tau,\mathrm{ps}}-4\log d+\log(\log d)\le x\big)\to\exp\Big(-\frac{1}{\sqrt{8\pi}}\exp(-\frac{x}{2})\Big). \tag{A.16}$$

(A.16) is the same as (3.10) except that in $M_n^{\tau,\mathrm{ps}}$ we use $\sigma^2_{a,\mathrm{ps}}/n_a$ to replace $\widehat{\sigma}^2(\widehat{\tau}_{a,ij})$ in (3.10).

We list a proof sketch of Theorem 3.3 and then present the detailed proof. For notational simplicity, we introduce

$$\begin{aligned} g^\tau_{ij}(\boldsymbol{X}_\alpha) &:= \mathbb{E}[\mathrm{sign}(X_{\alpha i} - X_{\beta i})\,\mathrm{sign}(X_{\alpha j} - X_{\beta j})|\boldsymbol{X}_\alpha], \\ g^\tau_{ij}(\boldsymbol{Y}_\alpha) &:= \mathbb{E}[\mathrm{sign}(Y_{\alpha i} - Y_{\beta i})\,\mathrm{sign}(Y_{\alpha j} - Y_{\beta j})|\boldsymbol{Y}_\alpha], \end{aligned} \tag{A.17}$$

with $\beta \neq \alpha$. The proof proceeds in three steps. In the first step, we use the Hoeffding method (Hoeffding, 1948) to decompose Kendall's tau as

$$\begin{aligned} \widehat{\tau}_{1,ij} &= \frac{2}{n_1}\sum_{\alpha=1}^{n_1} h^\tau_{ij}(\boldsymbol{X}_\alpha) + \frac{2}{n_1(n_1-1)}\Delta^\tau_{n_1,ij}, \\ \widehat{\tau}_{2,ij} &= \frac{2}{n_2}\sum_{\alpha=1}^{n_2} h^\tau_{ij}(\boldsymbol{Y}_\alpha) + \frac{2}{n_2(n_2-1)}\Delta^\tau_{n_2,ij}, \end{aligned} \tag{A.18}$$

where $h^\tau_{ij}(\boldsymbol{X}_\alpha) := g^\tau_{ij}(\boldsymbol{X}_\alpha) - \tau_{1,ij}$, $h^\tau_{ij}(\boldsymbol{Y}_\alpha) := g^\tau_{ij}(\boldsymbol{Y}_\alpha) - \tau_{2,ij}$ and

$$\begin{aligned} \Delta^\tau_{n_1,ij} &= \sum_{1\le k<\ell\le n_1} \big(\mathrm{sign}(X_{ki} - X_{\ell i})\,\mathrm{sign}(X_{kj} - X_{\ell j}) - \tau_{1,ij} - h^\tau_{ij}(\boldsymbol{X}_k) - h^\tau_{ij}(\boldsymbol{X}_\ell)\big), \\ \Delta^\tau_{n_2,ij} &= \sum_{1\le k<\ell\le n_1} \big(\mathrm{sign}(Y_{ki} - Y_{\ell i})\,\mathrm{sign}(Y_{kj} - Y_{\ell j}) - \tau_{2,ij} - h^\tau_{ij}(\boldsymbol{Y}_k) - h^\tau_{ij}(\boldsymbol{Y}_\ell)\big). \end{aligned}$$

$2\sum_{\alpha=1}^{n_1} h^\tau_{ij}(\boldsymbol{X}_\alpha)/n_1$ and $2\sum_{\alpha=1}^{n_1} h^\tau_{ij}(\boldsymbol{Y}_\alpha)/n_2$ are terms for the sum of i.i.d random variables and $2\Delta^\tau_{n_a,ij}/n_a(n_a-1)$ is the residual term. Similar to (A.5), we define $\widehat{T}_{ij}$ as

$$\widehat{T}_{ij} := \Big(\sum_{\alpha=1}^{n_1} h^\tau_{ij}(\boldsymbol{X}_\alpha)/n_1 - \sum_{\alpha=1}^{n_2} h^\tau_{ij}(\boldsymbol{Y}_\alpha)/n_2\Big)/\sqrt{\sigma^2_{1,\mathrm{ps}}/4n_1 + \sigma^2_{2,\mathrm{ps}}/4n_2}.$$

We then prove that the residual term $\Delta^\tau_{n_a,ij}/n_a(n_a-1)$ is negligible, i.e., to obtain Theorem 3.3, it suffices to prove that as $n, d \to \infty$, we have

$$\mathbb{P}\Big(\max_{(i,j)\in S}(\widehat{T}_{ij})^2 - 4\log d + \log(\log d) \le x\Big) \to \exp\Big(-\frac{1}{\sqrt{8\pi}}\exp(-\frac{x}{2})\Big). \tag{A.19}$$

In the second step, we prove that it is sufficient to take the maximum of $(\widehat{T}_{ij})^2$ over $S \setminus S_0$ but instead of over $S$ as in (A.19), i.e., to obtain Theorem 3.3, it suffices to prove that as $n, d \to \infty$, we have

$$\mathbb{P}\Big(\max_{(i,j)\in S\setminus S_0}(\widehat{T}_{ij})^2 - 4\log d + \log(\log d) \le x\Big) \to \exp\Big(-\frac{1}{\sqrt{8\pi}}\exp(-\frac{x}{2})\Big). \tag{A.20}$$

In the last step, we prove that it is sufficient to replace $\widehat{T}_{ij}$ in (A.20) with $T_{ij}$, i.e., to obtain Theorem 3.3, it suffices to prove that as $n$, $d \to \infty$, we have

$$\mathbb{P}\Big(\max_{(i,j)\in S\setminus S_0}(T_{ij})^2 - 4\log d + \log(\log d) \le x\Big) \to \exp\Big(-\frac{1}{\sqrt{8\pi}}\exp(-\frac{x}{2})\Big), \tag{A.21}$$

where $T_{ij}$ and $S_0$ are defined in (A.5) and (2.11). To prove sufficiency of (A.21), we need to prove $|\max_{(i,j)\in S\setminus S_0}(\widehat{T}_{ij})^2 - \max_{(i,j)\in S\setminus S_0}(T_{ij})^2| = o_p(1)$. $T_{ij}$ is the same as $\widehat{T}_{ij}$ except that in $T_{ij}$ we use $\sigma^2_{a,\mathrm{ps}}/4$ to replace $\zeta_{a,ij}$ in $\widehat{T}_{ij}$. Indeed, we prove that for $(i,j) \in S \setminus S_0$, we have $|4\zeta_{a,ij} - \sigma^2_{a,\mathrm{ps}}| \le C(\log d)^{-1-\alpha_0}$. The proof utilises the special structure of meta-elliptical distribution and is one of technical contributions.

After proving sufficiency of (A.21), to obtain Theorem 3.3, we need to prove that (A.21) is correct. Because (A.21) and (A.9) are the same, we prove (A.21) by following the same proof of (A.9). Hence, we complete the proof of Theorem 3.3.

The detailed proof is the following.

**Step (i)**. For notational simplicity, we introduce

$$N^{\tau,\mathrm{ps}}_{ij} := \frac{\widehat{\tau}_{1,ij} - \widehat{\tau}_{2,ij}}{\sqrt{\sigma^2_{1,\mathrm{ps}}/n_1 + \sigma^2_{2,\mathrm{ps}}/n_2}} = \widehat{T}_{ij} + \widehat{W}_{ij}, \;\; \widehat{W}_{ij} := \frac{\dfrac{\Delta^{\tau}_{n_1,ij}}{n_1(n_1-1)} - \dfrac{\Delta^{\tau}_{n_2,ij}}{n_2(n_2-1)}}{\sqrt{\sigma^2_{1,\mathrm{ps}}/4n_1 + \sigma^2_{2,\mathrm{ps}}/4n_2}}. \tag{A.22}$$

By (A.15) and (A.22), we have $M^{\tau,\mathrm{ps}}_{ij} := (N^{\tau,\mathrm{ps}}_{ij})^2$. We then introduce the following lemma to finish the proof of this step.

**Lemma A.3.** As $n$, $d \to \infty$, under Assumption **(A2)** we have

$$\big|\max_{1\le i<j\le d}(N^{\tau,\mathrm{ps}}_{ij})^2 - \max_{1\le i\le j\le d}(\widehat{T}_{ij})^2\big| = o_p(1). \tag{A.23}$$

The detailed proof of Lemma A.3 is in Supplement B.8. By Lemma A.3, we obtain that $\max_{1\le i<j\le d}(N^{\tau,\mathrm{ps}}_{ij})^2$ and $\max_{1\le i\le j\le d}(\widehat{T}_{ij})^2$ have the same limiting distribution as $n$, $d \to \infty$. Therefore, for obtaining Theorem 3.3, it suffices to prove (A.19) as $n$, $d \to \infty$.

**Step (ii)**. In the second step, we aim to prove that (A.20) is enough to prove the theorem. By setting $y_d = x + 4\log d - \log(\log d)$, we have

$$\Big|\mathbb{P}\Big(\max_{(i,j)\in S}(\widehat{T}_{ij})^2 > y_d\Big) - \mathbb{P}\Big(\max_{(i,j)\in S\setminus S_0}(\widehat{T}_{ij})^2 > y_d\Big)\Big| \le \mathbb{P}\Big(\max_{(i,j)\in S_0}(\widehat{T}_{ij})^2 > y_d\Big).$$

We then introduce an additional lemma.

**Lemma A.4.** Under Assumptions **(A1)** and **(A2)**, as $n$, $d \to \infty$, we have

$$\mathbb{P}\Big(\max_{(i,j)\in S_0}(\widehat{T}_{ij})^2 > y_d\Big) = o(1). \tag{A.24}$$

The detailed proof of Lemma A.4 is in Supplement B.9. By Lemma A.4, we obtain that as $n$, $d \to \infty$, $\mathbb{P}\Big(\max_{(i,j)\in S}(\widehat{T}_{ij})^2 > y_d\Big)$ and $\mathbb{P}\Big(\max_{(i,j)\in S\setminus S_0}(\widehat{T}_{ij})^2 > y_d\Big)$ have the same limit. Hence, we complete the proof of the second step.

**Step (iii)**. In this step, we aim to prove that it is sufficient to have (A.21) as $n$, $d \to \infty$. For this, we need the following lemma.

**Lemma A.5.** If $\boldsymbol{X}$ and $\boldsymbol{Y}$ follow the meta-elliptical distribution, we have

$$|4\zeta_{a,ij} - \sigma_{a,\mathrm{ps}}^2| \le C|\tau_{a,ij}|,$$

where $C$ is a constant irrelevant to $i$ or $j$. The definitions of $\zeta_{a,ij}$ and $\sigma_{a,\mathrm{ps}}$ are in (2.10) and (3.7).

The detailed proof of Lemma A.5 is in Supplement B.10. By the definition of $S_0$ in (2.11), for any $(i,j) \in S \setminus S_0$, we have $|\tau_{a,ij}| \le (\log d)^{-1-\alpha_0}$. By Lemma A.5, we have $|4\zeta_{a,ij} - \sigma_{a,\mathrm{ps}}^2| \le C|\tau_{a,ij}|$. Hence, for $(i,j) \in S \setminus S_0$, we have $|4\zeta_{a,ij} - \sigma_{a,\mathrm{ps}}^2| \le C|\tau_{a,ij}| \le C(\log d)^{-1-\alpha_0}$. Considering $\zeta_{a,ij} \ge r_a > 0$ (see Assumption **(A2)**), for $(i,j) \in S \setminus S_0$, we obtain

$$\left|\sigma_{a,\mathrm{ps}}^2/4\zeta_{a,ij} - 1\right| \le C(\log d)^{-1-\alpha_0} \quad \text{and} \quad \left|4\zeta_{a,ij}/\sigma_{a,\mathrm{ps}}^2 - 1\right| \le C(\log d)^{-1-\alpha_0}. \tag{A.25}$$

To prove the sufficiency of (A.21), we need $|\max_{(i,j)\in S\setminus S_0}(\widehat{T}_{ij})^2 - \max_{(i,j)\in S\setminus S_0}(T_{ij})^2| = o_p(1)$. For this, we calculate the relevant difference of $(\widehat{T}_{ij})^2$ and $(T_{ij})^2$ as

$$\left|\frac{(\widehat{T}_{ij})^2 - (T_{ij})^2}{(T_{ij})^2}\right| \le \left|\frac{\sigma_{1,\mathrm{ps}}^2 - 4\zeta_{1,ij}}{\sigma_{1,\mathrm{ps}}^2}\right| + \left|\frac{\sigma_{2,\mathrm{ps}}^2 - 4\zeta_{2,ij}}{\sigma_{2,\mathrm{ps}}^2}\right|.$$

By (A.25), we then have $|(\widehat{T}_{ij})^2 - (T_{ij})^2| \le C(T_{ij})^2(\log d)^{-1-\alpha_0}$ for $(i,j) \in S \setminus S_0$. Hence, we have

$$\begin{aligned}\left|\max_{(i,j)\in S\setminus S_0}(\widehat{T}_{ij})^2 - \max_{(i,j)\in S\setminus S_0}(T_{ij})^2\right| &\le \max_{(i,j)\in S\setminus S_0}\left|(\widehat{T}_{ij})^2 - (T_{ij})^2\right| \\ &\le C(\log d)^{-1-\alpha_0}\max_{(i,j)\in S\setminus S_0}(T_{ij})^2.\end{aligned}$$

Considering $\max_{(i,j)\in S\setminus S_0}(T_{ij})^2 = O_p(\log d)$, $|\max_{(i,j)\in S\setminus S_0}(\widehat{T}_{ij})^2 - \max_{(i,j)\in S\setminus S_0}(T_{ij})^2| = o_p(1)$ holds. Therefore, it suffices to prove (A.21) as $n,\ d \to \infty$. Because (A.21) and (A.9) are the same, we prove (A.21) by following the same proof of (A.9). Hence, we prove Theorem 3.3. □

## A.4. Proof of Theorem 3.4

***Proof.*** To obtain Theorem 3.4, we aim to prove (3.13), (3.14) and (3.15). The proof of (3.13) is the same as Theorem 2.4, as Kendall's tau is a special kind of U-statistic with a bounded kernel. Hence, we only need to prove (3.14) and (3.15).

First, we prove (3.14). Under the alternative hypothesis, $\tau_{1,ij} = \tau_{2,ij}$ cannot hold for all $(i,j) \in S$, where $S = \{(i,j) : 1 \le i < j \le d\}$. This motivates us to define

$$\begin{aligned} M_n^{1,\mathrm{plug}} &:= \max_{(i,j)\in S}\frac{(\widehat{\tau}_{1,ij} - \widehat{\tau}_{2,ij} - \tau_{1,ij} + \tau_{2,ij})^2}{\widehat{\sigma}_{\mathrm{plug}}^2(\widehat{\tau}_{1,ij}) + \widehat{\sigma}_{\mathrm{plug}}^2(\widehat{\tau}_{2,ij})}, \\ M_n^{\tau,\mathrm{plug}} &:= \max_{(i,j)\in S}\frac{(\widehat{\tau}_{1,ij} - \widehat{\tau}_{2,ij})^2}{\widehat{\sigma}_{\mathrm{plug}}^2(\widehat{\tau}_{1,ij}) + \widehat{\sigma}_{\mathrm{plug}}^2(\widehat{\tau}_{2,ij})}, \end{aligned} \tag{A.26}$$

because $M_n^{1,\mathrm{plug}}$ and $M_n^{\tau,\mathrm{plug}}$ are different under the alternative hypothesis. Similarly to Theorem 2.4, to prove (3.14), using the inequality $(a \pm b)^2 \le 2a^2 + 2b^2$, we obtain

$$(\tau_{1,ij} - \tau_{2,ij})^2 \le 2(\widehat{\tau}_{1,ij} - \widehat{\tau}_{2,ij} - \tau_{1,ij} + \tau_{2,ij})^2 + 2(\widehat{\tau}_{1,ij} - \widehat{\tau}_{2,ij})^2.$$

By the the definition of $M_n^{\tau,\mathrm{plug}}$ and $M_n^{1,\mathrm{plug}}$ in (A.26), we get

$$\max_{(i,j)\in S} \frac{(\tau_{1,ij} - \tau_{2,ij})^2}{\widehat{\sigma}^2_{\mathrm{plug}}(\widehat{u}_{1,ij}) + \widehat{\sigma}^2_{\mathrm{plug}}(\widehat{\tau}_{2,ij})} \le 2M_n^{1,\mathrm{plug}} + 2M_n^{\tau,\mathrm{plug}}. \tag{A.27}$$

To complete the proof of (3.14), we need two additional lemmas.

**Lemma A.6.** Under Assumption **(A2)**, as $n$, $d \to \infty$, we have

$$\mathbb{P}\big(M_n^{1,\mathrm{plug}} \le 4\log d - \frac{1}{2}\log(\log d)\big) \to 1. \tag{A.28}$$

The detailed proof of Lemma A.6 is in Supplement B.11.

**Lemma A.7.** Under Assumption **(A2)**, as $n$, $d \to \infty$, we have

$$\mathbb{P}\Big(\max_{(i,j)\in S} \frac{(\tau_{1,ij} - \tau_{2,ij})^2}{\widehat{\sigma}^2_{\mathrm{plug}}(\widehat{\tau}_{1,ij}) + \widehat{\sigma}^2_{\mathrm{plug}}(\widehat{\tau}_{2,ij})} \ge 16\log d\Big) \to 1, \tag{A.29}$$

uniformly for $(\mathbf{U}_1^\tau, \mathbf{U}_2^\tau) \in \mathbb{U}(4)$.

The detailed proof of Lemma A.7 is in Supplement B.12. Combining (A.27), Lemmas A.6 and A.7, under $(\mathbf{U}_1^\tau, \mathbf{U}_2^\tau) \in \mathbb{U}(4)$, as $n, d \to \infty$, with probability going to one, we have

$$M_n^{\tau,\mathrm{plug}} \ge \frac{1}{2} \max_{(i,j)\in S} \frac{(\tau_{1,ij} - \tau_{2,ij})^2}{\widehat{\sigma}^2_{\mathrm{plug}}(\widehat{\tau}_{1,ij}) + \widehat{\sigma}^2_{\mathrm{plug}}(\widehat{\tau}_{2,ij})} - M_n^{1,\mathrm{plug}} \ge 4\log d + \frac{1}{2}\log(\log d).$$

Therefore, under $(\mathbf{U}_1^\tau, \mathbf{U}_2^\tau) \in \mathbb{U}(4)$, as $n, d \to \infty$, we have

$$1 \ge \mathbb{P}\big(M_n^{\tau,\mathrm{plug}} \ge G^-(\alpha) + 4\log d - \log(\log d)\big) \ge \mathbb{P}\big(M_n^{\tau,\mathrm{plug}} \ge 4\log d + \frac{1}{2}\log(\log d)\big) \to 1,$$

where $G^-(\alpha) := -\log(8\pi) - 2\log\big(-\log(1-\alpha)\big)$. Hence, we prove (3.14).

Secondly, we aim to prove (3.15). Under the alternative hypothesis, $\tau_{1,ij} = \tau_{2,ij}$ cannot hold for all $(i,j) \in S$. This motivates us to define

$$M_n^{1,\mathrm{ps}} := \max_{(i,j)\in S} \frac{(\widehat{\tau}_{1,ij} - \widehat{\tau}_{2,ij} - \tau_{1,ij} + \tau_{2,ij})^2}{\sigma^2_{\mathrm{ps}}/n_1 + \sigma^2_{\mathrm{ps}}/n_2},\ M_n^{\tau,\mathrm{ps}} := \max_{(i,j)\in S} \frac{(\widehat{\tau}_{1,ij} - \widehat{\tau}_{2,ij})^2}{\sigma^2_{\mathrm{ps}}/n_1 + \sigma^2_{\mathrm{ps}}/n_2}, \tag{A.30}$$

as $M_n^{1,\mathrm{ps}}$ and $M_n^{\tau,\mathrm{ps}}$ are different under the alternative hypothesis.

We then introduce an additional lemma.

**Lemma A.8.** Under Assumptions **(A1)** and **(A2)**, as $n$, $d \to \infty$, we have

$$\mathbb{P}\big(M_n^{1,\mathrm{ps}} \leq 4\log d - \frac{1}{2}\log(\log d)\big) \to 1. \tag{A.31}$$

The detailed proof of Lemma A.8 is in Supplement B.13. To prove (3.15), using $(a \pm b)^2 \leq 2a^2 + 2b^2$, we obtain

$$(\tau_{1,ij} - \tau_{2,ij})^2 \leq 2(\widehat{\tau}_{1,ij} - \widehat{\tau}_{2,ij} - \tau_{1,ij} + \tau_{2,ij})^2 + 2(\widehat{\tau}_{1,ij} - \widehat{\tau}_{2,ij})^2.$$

By definitions of $M_n^{\tau,\mathrm{ps}}$ and $M_n^{1,\mathrm{ps}}$ in (A.30), we have

$$\max_{(i,j)\in S} \frac{(\tau_{1,ij} - \tau_{2,ij})^2}{\sigma_{1,\mathrm{ps}}^2/n_1 + \sigma_{2,\mathrm{ps}}^2/n_2} \leq 2M_n^{1,\mathrm{ps}} + 2M_n^{\tau,\mathrm{ps}}. \tag{A.32}$$

Under $(\mathbf{U}_1^\tau, \mathbf{U}_2^\tau) \in \mathbb{V}(4)$, we have

$$\max_{(i,j)\in S} \frac{(\tau_{1,ij} - \tau_{2,ij})^2}{\sigma_{1,\mathrm{ps}}^2/n_1 + \sigma_{2,\mathrm{ps}}^2/n_2} \geq 16\log d. \tag{A.33}$$

Combining Lemma A.8, (A.32) and (A.33), under $(\mathbf{U}_1^\tau, \mathbf{U}_2^\tau) \in \mathbb{V}(4)$, as $n$, $d \to \infty$, with probability going to one, we have

$$M_n^{\tau,\mathrm{ps}} \geq \frac{1}{2}\max_{(i,j)\in S} \frac{(\tau_{1,ij} - \tau_{2,ij})^2}{\sigma_{1,\mathrm{ps}}^2/n_1 + \sigma_{2,\mathrm{ps}}^2/n_2} - M_n^{1,\mathrm{plug}} \geq 4\log d + \frac{1}{2}\log(\log d).$$

Therefore, under $(\mathbf{U}_1^\tau, \mathbf{U}_2^\tau) \in \mathbb{V}(4)$, as $n, d \to \infty$, we have

$$1 \geq \mathbb{P}\big(M_n^{\tau,\mathrm{ps}} \geq G^-(\alpha) + 4\log d - \log(\log d)\big) \geq \mathbb{P}\big(M_n^{\tau,\mathrm{ps}} \geq 4\log d + \frac{1}{2}\log(\log d)\big) \to 1,$$

where $G^-(\alpha) := -\log(8\pi) - 2\log\big(-\log(1-\alpha)\big)$. Hence, we prove (3.15). This completes the proof of Theorem 3.4. □

## A.5. Proof of Theorem 3.5

***Proof.*** It suffices to take $\mathcal{T}_\alpha$ to be the set of level $\alpha$ tests over the normal distributions with covariance matrix $\mathbf{\Sigma}$, where $\mathrm{Diag}(\mathbf{\Sigma}) = \mathbf{I}_d$, since it contains all the $\alpha$-level tests over the collection of the assumed distributions. For these normal distributions, we define

$$\mathbb{H}(c') = \big\{(\mathbf{\Sigma}_1, \mathbf{\Sigma}_2) : \sigma_{1,ii} = 1, \sigma_{2,ii} = 1, \|\mathbf{\Sigma}_1 - \mathbf{\Sigma}_2\|_{\max} \geq c'\sqrt{\log d/n}\big\}.$$

By $\sigma_{ij} = \sin(\tau_{ij}\pi/2)$ (Theorem E.3 in Supplement E), we have $|\sigma_{1,ij} - \sigma_{2,ij}| \leq |\tau_{1,ij} - \tau_{2,ij}|\pi/2$. Therefore, for any $c'$, there is a $c_0$ such that $\mathbb{H}(c') \subset \mathbb{U}(c_0)$. For simplicity, we set $n_1 = n_2 = n$.

Let's consider the Gaussian setting and define

$$\mathcal{F}(\rho) = \{\mathbf{\Sigma} = \mathbf{I}_d + \rho \mathbf{e}_1 \mathbf{e}_j^T + \rho \mathbf{e}_j \mathbf{e}_1^T,\ \mathbf{e}_k = (\underbrace{0, \ldots, 0}_{k-1}, 1, 0, \ldots, 0) \text{ for } 1 \le k \le d,\ j = 2, \ldots, d\},$$

where $\rho = c' \sqrt{\log d / n}$ and $\mu_\rho$ is the uniform measure on $\mathcal{F}(\rho)$. For simplicity, under the null hypothesis, we set $\mathbf{\Sigma}_1 = \mathbf{\Sigma}_2 = \mathbf{I}_d$. Under the alternative hypothesis, we set $\mathbf{\Sigma}_1 = \mathbf{\Sigma} \sim \mu_\rho$ and $\mathbf{\Sigma}_2 = \mathbf{I}_d$. Therefore, under the alternative hypothesis, there is a $c_0$ such that we have $(\mathbf{\Sigma}_1, \mathbf{\Sigma}_2) \in \mathbb{H}(c') \subset \mathbb{U}(c_0)$.

Let $\mathbb{P}_{\mathbf{\Sigma}}$ denote the probability measure of samples for $\boldsymbol{X} \sim N(\mathbf{0}, \mathbf{\Sigma}_1)$ and $\boldsymbol{Y} \sim N(\mathbf{0}, \mathbf{\Sigma}_2)$ with $\mathbf{\Sigma}_1 = \mathbf{\Sigma}$ and $\mathbf{\Sigma}_2 = \mathbf{I}_d$. We set $\mathbb{P}_{\mu_\rho} = \int \mathbb{P}_{\mathbf{\Sigma}} d\mu_\rho(\mathbf{\Sigma})$. In particular, let $\mathbb{P}_0$ denote the probability measure of samples for $\boldsymbol{X} \sim N(\mathbf{0}, \mathbf{\Sigma}_1)$ and $\boldsymbol{Y} \sim N(\mathbf{0}, \mathbf{\Sigma}_2)$ with $\mathbf{\Sigma}_1 = \mathbf{\Sigma}_2 = \mathbf{I}_d$. We then have

$$\inf_{T_\alpha \in \mathcal{T}_\alpha} \sup_{\mathbf{\Sigma} \in \mathcal{F}(\rho)} \mathbb{P}_{\mathbf{\Sigma}}(T_\alpha = 0) \ge 1 - \alpha - \sup_{A: \mathbb{P}_0(A) \le \alpha} |\mathbb{P}_{\mu_\rho}(A) - \mathbb{P}_0(A)| \ge 1 - \alpha - \frac{1}{2} ||\mathbb{P}_{\mu_\rho} - \mathbb{P}_0||_{TV},$$

where $||\cdot||_{TV}$ denotes the total variation norm. By setting $L_{\mu_\rho}(\mathbf{z}) := \frac{d\mathbb{P}_{\mu_\rho}}{d\mathbb{P}_0}(\mathbf{z})$, considering the Jensen's inequality, we have

$$||\mathbb{P}_{\mu_\rho} - \mathbb{P}_0||_{TV} = \int |L_{\mu_\rho}(\mathbf{z}) - 1| d\mathbb{P}_0(\mathbf{z}) = \mathbb{E}_{\mathbb{P}_0} |L_{\mu_\rho}(\boldsymbol{Z}) - 1| \le (|\mathbb{E}_{\mathbb{P}_0} L^2_{\mu_\rho}(\boldsymbol{Z}) - 1|)^{1/2}.$$

Hence, as long as $\mathbb{E}_{\mathbb{P}_0} L^2_{\mu_\rho}(\boldsymbol{Z}) = 1 + o(1)$, we have

$$\inf_{T_\alpha \in \mathcal{T}_\alpha} \sup_{\mathbf{\Sigma} \in \mathcal{F}(\rho)} \mathbb{P}_{\mathbf{\Sigma}}(T_\alpha = 0) \ge 1 - \alpha - o(1) > 0,$$

which is the desired result. We then aim to prove $\mathbb{E}_{\mathbb{P}_0} L^2_{\mu_\rho}(\boldsymbol{Z}) = 1 + o(1)$. By construction, we have

$$L_{\mu_\rho} = \frac{1}{d-1} \sum_{\mathbf{\Sigma} \in \mathcal{F}(\rho)} \Big( \prod_{i=1}^n \frac{1}{|\mathbf{\Sigma}|^{1/2}} \exp(-\frac{1}{2} \boldsymbol{Z}_i^T (\mathbf{\Omega} - \mathbf{I}_d) \boldsymbol{Z}_i) \Big),$$

where $\mathbf{\Omega} := \mathbf{\Sigma}^{-1}$ and $\{\boldsymbol{Z}_i\}$ are independent random vectors with $\boldsymbol{Z}_i = (Z_{i1}, Z_{i2}, \ldots, Z_{id})^T \sim N(\mathbf{0}, \mathbf{I}_d)$. Therefore, we have

$$\mathbb{E}_{\mathbb{P}_0} L^2_{\mu_\rho} = \frac{1}{(d-1)^2} \sum_{\mathbf{\Sigma}_1, \mathbf{\Sigma}_2 \in \mathcal{F}(\rho)} \mathbb{E} \left( \prod_{i=1}^n \frac{1}{|\mathbf{\Sigma}_1|^{1/2} |\mathbf{\Sigma}_2|^{1/2}} \exp\Big( -\frac{1}{2} \boldsymbol{Z}_i^T (\mathbf{\Omega}_1 + \mathbf{\Omega}_2 - 2\mathbf{I}_d) \boldsymbol{Z}_i \Big) \right), \quad \text{(A.34)}$$

where $\mathbf{\Omega}_i = \mathbf{\Sigma}_i^{-1}$ for $i = 1, 2$. By setting

$$\mathbf{A} = \frac{\rho}{1-\rho^2} \begin{pmatrix} 2\rho & -1 & -1 \\ -1 & \rho & 0 \\ -1 & 0 & \rho \end{pmatrix}, \qquad \mathbf{B} = \frac{2\rho}{1-\rho^2} \begin{pmatrix} \rho & -1 \\ -1 & \rho \end{pmatrix},$$

$\boldsymbol{Z}_{i,\{1,2\}} = (Z_{i1}, Z_{i2})^T$ and $\boldsymbol{Z}_{i,\{1,2,3\}} = (Z_{i1}, Z_{i2}, Z_{i3})^T$, we have

$$\mathbb{E}_{\mathbb{P}_0} L^2_{\mu_\rho} = \underbrace{\frac{d-2}{d-1} \prod_{i=1}^n \left( \frac{1}{1-\rho^2} \mathbb{E} \exp\left( -\frac{1}{2} \boldsymbol{Z}^T_{i,\{1,2,3\}} \mathbf{A} \boldsymbol{Z}_{i,\{1,2,3\}} \right) \right)}_{A_3}$$
$$+ \underbrace{\frac{1}{d-1} \prod_{i=1}^n \left( \frac{1}{1-\rho^2} \mathbb{E} \exp\left( -\frac{1}{2} \boldsymbol{Z}^T_{i,\{1,2\}} \mathbf{B} \boldsymbol{Z}_{i,\{1,2\}} \right) \right)}_{A_4},$$

where $A_3$ represents the sum of terms with $\boldsymbol{\Sigma}_1 \neq \boldsymbol{\Sigma}_2$ in (A.34) and $A_4$ represents the sum of terms with $\boldsymbol{\Sigma}_1 = \boldsymbol{\Sigma}_2$ in (A.34). For $A_3$, by the standard argument in calculating moment generating function of the Gaussian quadratic form (Baldessari, 1967), we have

$$A_3 = \frac{d-2}{d-1} \cdot \frac{1}{(1-\rho^2)^n} \left( (1+\lambda_1(\mathbf{A}))(1+\lambda_2(\mathbf{A}))(1+\lambda_3(\mathbf{A})) \right)^{-n/2}.$$

Moreover, we have $(1+\lambda_1(\mathbf{A}))(1+\lambda_2(\mathbf{A}))(1+\lambda_3(\mathbf{A})) = |\mathbf{A}+\mathbf{I}_3| = (1-\rho^2)^{-2}$. Therefore, we have $A_3 = (d-2)/(d-1) = 1 + o(1)$. For $A_4$, it is easy to get $\lambda_1(\mathbf{B}) = 2\rho/(1-\rho)$ and $\lambda_2(\mathbf{B}) = -2\rho/(1+\rho)$. Accordingly, similarly to the calculation of $A_3$, we have $A_4 = (d-1)^{-1} \cdot (1-\rho^2)^{-n}$. Considering $\rho = c'\sqrt{\log d/n}$, as long as $c' < 1$, we have

$$A_4 \leq \frac{1}{d-1} \cdot (1 - c'^2 \log d/n)^{-n} = (d-1)^{-1} \exp(c'^2 \log d)(1+o(1)) = o(1), \qquad (A.35)$$

as $n, d \to \infty$. Combining $A_3 = 1 + o(1)$ and (A.35), we prove $\mathbb{E}_{\mathbb{P}_0} L^2_{\mu_\rho}(\boldsymbol{Z}) = 1 + o(1)$. This completes the proof. □

## A.6. Proof of Theorem 3.6

***Proof.*** To prove Theorem 3.6, we introduce the set $B_0$. We set it as

$$C_0 = \{(i,j) : i \in \Omega(r) \bigcup \Gamma\} \bigcup \{(i,j) : j \in \Omega(r) \bigcup \Gamma\} \qquad \text{and} \qquad B_0 = S_0 \bigcup C_0,$$

where $\Omega(r)$, $\Gamma$ and $S_0$ are defined in (3.17), Assumption **(A4)** and (2.11). For Kendall's tau matrices, we set $S = \{(i,j) : 1 \leq i < j \leq d\}$ and $q = d$. By the definition of $B_0$, we have $|\tau_{a,k\ell}| \leq r < 1$ for any $k \neq \ell \in \{i_1, j_1, i_2, j_2\}$, where $(i_1, j_1) \in S \setminus B_0$ and $(i_2, j_2) \in S \setminus B_0$. For any $(i,j) \in S \setminus B_0$, we also have $|\tau_{a,ij}| \leq (\log d)^{-1-\alpha_0}$, as $S_0$ is a subset of $B_0$. By Assumptions **(A1)** and **(A4)**, we have $|B_0| = o(d^2)$.

We then prove Theorem 3.6 similarly to the proofs of Theorems 3.2 and 3.3 except that we replace $S_0$ with $B_0$. However, as we don't require Assumption **(A3)** in Theorem 3.6, we don't have Lemma A.5. Therefore, we only need to prove (A.18) under Assumptions **(A1)**, **(A2)** and **(A4)**.

We then begin to prove (A.18) under Assumptions **(A1)**, **(A2)** and **(A4)**. As we use $B_0$ to replace $S_0$ in the proofs of Theorem 3.2 and 3.3, we need to redefine some notations. After rearranging the two-dimensional indices $\{(i,j) : (i,j) \in S \setminus B_0\}$ in any order, we set them as $\{(i_k, j_k) : 1 \le k \le h\}$ with $h = |S \setminus B_0|$. We denote $T_k = T_{i_k j_k}$, where the definition of $T_{ij}$ is the same as (A.5) except that we use Kendall's tau as the U-statistic in (A.5). We only need to check whether (A.18) is correct for $\{(i_k, j_k) : 1 \le k \le h\}$ under Assumptions **(A1)**, **(A2)** and **(A4)**. The definition of $\boldsymbol{N}_\ell$ is same as (A.15). By $\boldsymbol{N}_\ell$'s definition in (A.15), the entry in the $a$-th row, $b$-th column of $\mathrm{Var}(\boldsymbol{N}_\ell)$ is

$$\frac{n_2 \operatorname{Cov}\big(h_{i_{k_a} j_{k_a}}(\boldsymbol{X}_1), h_{i_{k_b} j_{k_b}}(\boldsymbol{X}_1)\big)/n_1 + \operatorname{Cov}\big(h_{i_{k_a} j_{k_a}}(\mathbf{Y}_1), h_{i_{k_b} j_{k_b}}(\mathbf{Y}_1)\big)}{\sqrt{n_2 \zeta_{1,i_{k_a} j_{k_a}}/n_1 + \zeta_{a,i_{k_a} j_{k_a}}}\sqrt{n_2 \zeta_{1,i_{k_b} j_{k_b}}/n_1 + \zeta_{a,i_{k_b} j_{k_b}}}}, \tag{A.36}$$

where $h_{ij}$ is defined in (2.9). Apparently, when $a = b$, (A.36) equals one.

The proof of (A.18) proceeds in three steps. In the first step, for any $(i,j)$, $(k,\ell) \in S \setminus B_0$, we define a condition $(\star)$. This condition specifies some $(i,j)$, $(k,\ell) \in S \setminus B_0$, for which there exits an $i_1 \in \{i,j,k,\ell\}$ such that with any $j_1 \in \{i,j,k,\ell\} \setminus i_1$, we have $|\tau_{a,i_1 j_1}| = O((\log d)^{-1-\alpha_0})$. In the second step, we prove that for $(i,j)$, $(k,\ell) \in S \setminus B_0$ satisfying $(\star)$, as $n$, $d \to \infty$, we have

$$\frac{n_2 \operatorname{Cov}\big(h_{ij}(\boldsymbol{X}_1), h_{k\ell}(\boldsymbol{X}_1)\big)/n_1 + \operatorname{Cov}\big(h_{ij}(\mathbf{Y}_1), h_{k\ell}(\mathbf{Y}_1)\big)}{\sqrt{n_2 \zeta_{1,ij}/n_1 + \zeta_{2,ij}}\sqrt{n_2 \zeta_{1,k\ell}/n_1 + \zeta_{2,k\ell}}} = O((\log d)^{-1-\alpha_0}). \tag{A.37}$$

In the third step, we prove that for $(i,j)$, $(k,\ell) \in S \setminus B_0$ dissatisfying $(\star)$, we have

$$\left|\frac{n_2 \operatorname{Cov}\big(h_{ij}(\boldsymbol{X}_1), h_{k\ell}(\boldsymbol{X}_1)\big)/n_1 + \operatorname{Cov}\big(h_{ij}(\mathbf{Y}_1), h_{k\ell}(\mathbf{Y}_1)\big)}{\sqrt{n_2 \zeta_{1,ij}/n_1 + \zeta_{2,ij}}\sqrt{n_2 \zeta_{1,k\ell}/n_1 + \zeta_{2,k\ell}}}\right| \le C < 1. \tag{A.38}$$

as $n$, $d \to \infty$. By the proof of Lemma 5 in Cai et al. (2013), (A.37) and (A.38) are sufficient conditions for obtaining (A.18). Therefore, by proving (A.18), we finish the proof of Theorem 3.6.

Before presenting the detailed proof of (A.18) , we first introduce two additional lemmas. These lemmas simplify the proof process.

**Lemma A.9.** Suppose $(X_1, X_2, X_3, X_4)^T \sim N(\mathbf{0}, \boldsymbol{\Sigma}_{\text{full}})$ is Gaussian distributed with

$$\boldsymbol{\Sigma}_{\text{full}} = \begin{bmatrix} 1 & \boldsymbol{\varrho}^T \\ \boldsymbol{\varrho} & \boldsymbol{\Sigma} \end{bmatrix} \quad \text{and} \quad \operatorname{Diag}(\boldsymbol{\Sigma}) = \mathbf{I}_3, \tag{A.39}$$

where $\boldsymbol{\varrho} = (\rho_1, \rho_2, \rho_3)^T$. We have

$$\begin{aligned} &\big|\mathbb{E}(\Phi(X_1) - 1/2)(\Phi(X_2) - 1/2)(\Phi(X_3) - 1/2)(\Phi(X_4) - 1/2)\big| \\ &\qquad \le \left(\frac{1}{8\pi} + \frac{1}{4\sqrt{2\pi}}\right)\sqrt{\boldsymbol{\varrho}^T \boldsymbol{\Sigma}^{-1} \boldsymbol{\varrho}}, \end{aligned} \tag{A.40}$$

where $\Phi(\cdot)$ is the cumulative distribution function of a standard normal distribution.

The detailed proof of Lemma A.9 is in Supplement B.14.

**Lemma A.10.** Suppose $(X_1, X_2, X_3, X_4)^T \sim N(\mathbf{0}, \mathbf{\Sigma}_{\text{full}})$ is Gaussian distributed with

$$\mathbf{\Sigma}_{\text{full}} = \begin{bmatrix} 1 & \rho_1 & a_1 & a_2 \\ \rho_1 & 1 & a_3 & a_4 \\ a_1 & a_3 & 1 & \rho_2 \\ a_2 & a_4 & \rho_2 & 1 \end{bmatrix}.$$

Then, when $|\rho_1|, |\rho_2|, |a_1|, \ldots, |a_4| \le r < 1$, we have

$$C_r := \sup_{\substack{|\rho_1|,|\rho_2|,\\ |a_1|,\ldots,|a_4|\le r}} |\text{Corr}\{(\Phi(X_1) - 1/2)(\Phi(X_2) - 1/2), (\Phi(X_3) - 1/2)(\Phi(X_4) - 1/2)\}| < 1.$$

Moreover, we have $C_r = 1$ only when $r = 1$ and the set $\{\rho_1, \rho_2, a_1, \ldots, a_4\}$ attains the boundary.

The detailed proof of Lemma A.10 is in Supplement B.15. After introducing the two lemmas, we then prove (A.18) in detail.

**Step (i)**. In this step, we define the condition $(\star)$. For this, we set graph $G_{ijk\ell} = (V_{ijk\ell}, E_{ijk\ell})$, where $V_{ijk\ell} = \{i,\ j,\ k,\ \ell\}$ is the set of vertices and $E_{ijk\ell}$ is the set of edges. We say that there is an edge between $a \neq b \in \{i,\ j,\ k,\ \ell\}$ if and only if $|\tau_{ab}| \ge (\log d)^{-1-\alpha_0}$. If the number of different vertices in $V_{ijk\ell}$ is 3, we say that $G_{ijk\ell}$ is a three vertices graph $(3 - G)$. If the number of different vertices in $V_{ijk\ell}$ is 4, we say that $G_{ijk\ell}$ is a four vertices graph $(4 - G)$. If there is no edge connected to a vertex in $G_{ijk\ell}$, we say that it is isolated. For any $1 \le k_a \neq k_b \le h$, we have that $G_{i_{k_a} j_{k_a} i_{k_b} j_{k_b}}$ is either $3 - G$ or $4 - G$. For any $1 \le k_a \neq k_b \le h$, we say $\mathcal{G} := G_{i_{k_a} j_{k_a} i_{k_b} j_{k_b}}$ satisfies $(\star)$, if a graph $\mathcal{G}$ has the property:

$(\star)$ If $G_{i_{k_a} j_{k_a} i_{k_b} j_{k_b}}$ is $4 - G$, there is at least one isolated vertex in $G_{i_{k_a} j_{k_a} i_{k_b} j_{k_b}}$; otherwise $G_{i_{k_a} j_{k_a} i_{k_b} j_{k_b}}$ is $3 - G$ and $E_{i_{k_a} j_{k_a} i_{k_b} j_{k_b}} = \emptyset$.

**Step (ii)**. In this step, we prove that for $(i,j)$, $(k,\ell) \in S \setminus B_0$ satisfying the condition $(\star)$, (A.37) is correct. To prove (A.37), by $\zeta_{a,ij} \ge r_a > 0$ (see Assumption **(A2)**), it suffices to prove

$$\begin{aligned} \text{Cov}\big(h_{ij}(\boldsymbol{X}_1), h_{k\ell}(\boldsymbol{X}_1)\big) &= O\big((\log d)^{-1-\alpha_0}\big), \\ \text{Cov}\big(h_{ij}(\boldsymbol{Y}_1), h_{k\ell}(\boldsymbol{Y}_1)\big) &= O\big((\log d)^{-1-\alpha_0}\big). \end{aligned} \tag{A.41}$$

For simplicity, we only show the proof of $\boldsymbol{X}$. We treat Kendall's tau as a special kind of U-statistics. By (2.9), we have $h_{ij}(\boldsymbol{X}_1) = g_{ij}(\boldsymbol{X}_1) - \tau_{1,ij}$. Therefore, we have

$$\text{Cov}\big(h_{ij}(\boldsymbol{X}_1), h_{k\ell}(\boldsymbol{X}_1)\big) = \mathbb{E}[g_{ij}(\boldsymbol{X}_1) g_{k\ell}(\boldsymbol{X}_1)] - \tau_{1,ij}\tau_{1,k\ell}. \tag{A.42}$$

By $(i,j)$, $(k,\ell) \in S \setminus B_0$, we have

$$\tau_{1,ij} = O((\log d)^{-1-\alpha_0}) \qquad \text{and} \qquad \tau_{1,k\ell} = O((\log d)^{-1-\alpha_0}). \tag{A.43}$$

Combining (A.42) and (A.43), for getting (A.41), it suffices to show

$$\mathbb{E}[g_{ij}(\boldsymbol{X}_1) g_{k\ell}(\boldsymbol{X}_1)] = O\big((\log d)^{-1-\alpha_0}\big). \tag{A.44}$$

To prove (A.44), by the definition of $g_{ij}$ in (A.17), we know for Kendall's tau, we have

$$\begin{aligned} g_{ij}(\boldsymbol{X}_1) =& \mathbb{P}\big((X_{1i}-X_{2i})(X_{1j}-X_{2j})>0 | X_{1i}, X_{1j}\big) \\ &- \mathbb{P}\big((X_{1i}-X_{2i})(X_{1j}-X_{2j})<0 | X_{1i}, X_{1j}\big). \end{aligned} \tag{A.45}$$

Genz and Bretz (2009) show that for a two-dimensional normal vector $(N_1, N_2)^T \sim N(0, \widetilde{\boldsymbol{\Sigma}})$ with

$$\widetilde{\boldsymbol{\Sigma}} = \begin{bmatrix} 1 & \rho \\ \rho & 1 \end{bmatrix},$$

we have

$$\mathbb{P}(N_1 > a_1, N_2 > a_2) = \Phi(-a_1)\Phi(-a_2) + \frac{1}{2\pi}\int_0^{\sin^{-1}(\rho)} \exp\big(-\frac{a_1^2 - 2a_1a_2\sin(\theta) + a_2^2}{2\cos^2(\theta)}\big)d\theta. \tag{A.46}$$

Noticing that $\boldsymbol{X}_1$ and $\boldsymbol{X}_2$ follow the Gaussian copula distribution, by (A.45), (A.46) and Theorem E.3, we have

$$g_{ij}(\boldsymbol{X}_1) = (2F_i(X_{1i}) - 1)(2F_j(X_{1j}) - 1) + O(\tau_{1,ij}), \tag{A.47}$$

where $F_t$ is the cumulative distribution function of $X_{1t}$. Combining (A.43) and (A.47), to obtain (A.44), we only need to prove

$$\begin{aligned} 2^4\mathbb{E}\big[\big(F_i(X_{1i}) - 1/2\big)\big(F_j(X_{1j}) - 1/2\big)\big(F_k(X_{1k}) - 1/2\big)\big(F_\ell(X_{1l}) - 1/2\big)\big] \\ = O\big((\log d)^{-1-\alpha_0}\big). \end{aligned} \tag{A.48}$$

As $\boldsymbol{X}$ belongs to the Gaussian copula family,we assume $(X_{1i}, X_{1j}, X_{1k}, X_{1\ell})^T \sim N(\mathbf{0}, \boldsymbol{\Sigma}_{\text{full}})$ is Gaussian distributed with

$$\text{Diag}(\boldsymbol{\Sigma}) = \mathbf{I}_3 \quad \text{and} \quad \boldsymbol{\Sigma}_{\text{full}} = \begin{bmatrix} 1 & \boldsymbol{\varrho}^T \\ \boldsymbol{\varrho} & \boldsymbol{\Sigma} \end{bmatrix},$$

for simplicity. Accordingly, we have $F_t(\cdot) = \Phi(\cdot)$ for $t = 1, \ldots, d$, where $\Phi(\cdot)$ is the cumulative distribution function of a standard normal distribution. Hence, to obtain (A.48), it suffices to prove that as $n, d \to \infty$, we have

$$\begin{aligned} \mathbb{E}\big[\big(\Phi(X_{1i}) - 1/2\big)\big(\Phi(X_{1j}) - 1/2\big)\big(\Phi(X_{1k}) - 1/2\big)\big(\Phi(X_{1\ell}) - 1/2\big)\big] \\ = O\big((\log d)^{-1-\alpha_0}\big). \end{aligned} \tag{A.49}$$

If $G_{ijk\ell}$ is $4-G$ and satisfies $(\star)$, without loss of generality, we suppose that $X_{1i}$ is the isolated vertex. Hence, by Theorem E.3, we have $\|\boldsymbol{\varrho}\|_2 = O\big((\log d)^{-1-\alpha_0}\big)$. Noticing that the smallest eigenvalue of any 4 by 4 principal sub-matrix of $\boldsymbol{\Sigma}_a$ is bounded away from 0, by Lemma A.9, we have (A.49). Hence (A.37) is correct. If $G_{ijk\ell}$ is $3-G$ and satisfies $(\star)$, we have (A.49) similarly. This completes the proof of **Step (ii)**.

**Step (iii)**. For $(i, j)$, $(k, \ell) \in S \setminus B_0$, to get (A.38) as $n, d \to \infty$, we only need to prove

$$\big|\text{Corr}\big(h_{ij}(\boldsymbol{X}_1), h_{k\ell}(\boldsymbol{X}_1)\big)\big| \le C < 1. \tag{A.50}$$

By (A.47) and (A.43), to prove (A.50), we only need to prove that as $n,\ d \to \infty$, we have

$$\big|\mathrm{Corr}\big((F_i(X_{1i}) - 1/2)(F_j(X_{1j}) - 1/2), (F_k(X_{1k}) - 1/2)(F_\ell(X_{1\ell}) - 1/2)\big)\big| \le C < 1.$$

For simplicity, we suppose $(X_{1i}, X_{1j}, X_{1k}, X_{1\ell})^T \sim N(\mathbf{0}, \boldsymbol{\Sigma}_{\text{full}})$ is Gaussian distributed with

$$\mathrm{Diag}(\boldsymbol{\Sigma}) = \mathbf{I}_3 \quad \text{and} \quad \boldsymbol{\Sigma}_{\text{full}} = \begin{bmatrix} 1 & \boldsymbol{\varrho}^T \\ \boldsymbol{\varrho} & \boldsymbol{\Sigma} \end{bmatrix}.$$

We then need to prove that as $n,\ d \to \infty$, we have

$$\big|\mathrm{Corr}\big((\Phi(X_{1i}) - 1/2)(\Phi(X_{1j}) - 1/2), (\Phi(X_{1k}) - 1/2)(\Phi(X_{1\ell}) - 1/2)\big)\big| \le C < 1. \quad \text{(A.51)}$$

By Lemma A.10, we have (A.51). Therefore, we prove (A.38). By following the proof of Lemma 5 in Cai et al. (2013), we prove (A.18) under Assumptions **(A1)**, **(A2)** and **(A4)**. This completes the proof of Theorem 3.6.

□

## A.7. Proof of Theorems 2.6 and Remark 3.8

***Proof.*** The proof of Theorem 2.6 is similar and simpler than the proof of Theorem 2.2. The proof of Remark 3.8 is similar to the proofs of Theorems 3.2, 3.3 and 3.6. Due to the close similarity, the proof is omitted. □

# Supplement B: Proofs of Lemmas in Appendix A and Supplement A

In this appendix, we present proofs of lemmas in Appendix A and Supplement A. In the sequel, we use $C$, $C_1$, $C_2$, …, to denote constants that do not depend on $n$, $d$, $q$ and can vary from place to place.

## B.1. Proof of Lemma A.1 of Appendix A

***Proof.*** By Lemma D.3, we have that $m^2\zeta_{a,ij}$ is the limiting variance of $\sqrt{n_a}\widehat{u}_{a,ij}$ as $n$, $q \to \infty$. Combining Lemma D.3 and the definition of $h_{ij}$ in (2.9), we also have that $\zeta_{1,ij}$ and $\zeta_{2,ij}$ are variances of $h_{ij}(\boldsymbol{X}_\ell)$ and $h_{ij}(\mathbf{Y}_\ell)$.

In Lemma A.1, we aim to prove (A.1). For simplicity, we only provide the proof of $\boldsymbol{X}$. By the definition of Jackknife variance estimator $\widehat{\sigma}^2(\widehat{u}_{a,ij})$ in (2.5), we rewrite (A.1) as

$$\mathbb{P}\Big(\max_{1\le i,j\le q}\Big|\frac{m^2(n_1-1)}{(n_1-m)^2}\sum_{\alpha=1}^{n_1}(q_{1\alpha,ij} - \widehat{u}_{1,ij})^2 - m^2\zeta_{1,ij}\Big| \ge C\frac{\varepsilon_n}{\log q}\Big) = o(1), \quad \text{(B.1)}$$

where $q_{1\alpha,ij}$ is defined as

$$q_{1\alpha,ij} := \binom{n_1-1}{m-1}^{-1} \sum_{\substack{1\le \ell_1<\cdots<\ell_{m-1}\le n_1\\ \ell_k\neq\alpha, k=1,\cdots,m-1}} \Phi_{ij}(\boldsymbol{X}_\alpha, \boldsymbol{X}_{\ell_1},\ldots,\boldsymbol{X}_{\ell_{m-1}}).$$

To prove Lemma A.1, we also need the centralized version of $q_{1\alpha,ij}$. By $\Psi_{ij}(\boldsymbol{X}_{\ell_1},\ldots,\boldsymbol{X}_{\ell_m}) = \Phi_{ij}(\boldsymbol{X}_{\ell_1},\ldots,\boldsymbol{X}_{\ell_m}) - u_{1,ij}$, we define the centralized version of $q_{1\alpha,ij}$ as

$$\widetilde{q}_{1\alpha,ij} := \binom{n_1-1}{m-1}^{-1} \sum_{\substack{1\le \ell_1<\cdots<\ell_{m-1}\le n_1\\ \ell_k\neq\alpha, k=1,\cdots,m-1}} \Psi_{ij}(\boldsymbol{X}_\alpha, \boldsymbol{X}_{\ell_1},\ldots,\boldsymbol{X}_{\ell_{m-1}}). \tag{B.2}$$

To bound the left hand side of (B.1), it is easy to obtain

$$\begin{aligned}&\mathbb{P}\Big(\max_{1\le i,j\le q}\Big|\frac{m^2(n_1-1)}{(n_1-m)^2}\sum_{\alpha=1}^{n_1}(q_{1\alpha,ij}-\widehat{u}_{1,ij})^2 - m^2\zeta_{1,ij}\Big| \ge C\frac{\varepsilon_n}{\log q}\Big)\\ &\quad\le q^2\max_{1\le i,j\le q}\mathbb{P}\Big(\Big|\frac{m^2(n_1-1)}{(n_1-m)^2}\sum_{\alpha=1}^{n_1}(q_{1\alpha,ij}-\widehat{u}_{1,ij})^2 - m^2\zeta_{1,ij}\Big| \ge C\frac{\varepsilon_n}{\log q}\Big).\end{aligned}$$

We then replace $q_{1\alpha,ij}$ and $\widehat{u}_{1,ij}$ with their centralized versions $\widetilde{q}_{1\alpha,ij}$ and $\widetilde{u}_{1,ij}$ to obtain

$$\begin{aligned}&q^2\max_{1\le i,j\le q}\mathbb{P}\Big(\Big|\frac{m^2(n_1-1)}{(n_1-m)^2}\sum_{\alpha=1}^{n}(q_{1\alpha,ij}-\widehat{u}_{1,ij})^2 - m^2\zeta_{1,ij}\Big| \ge C\frac{\varepsilon_n}{\log q}\Big)\\ &\quad= q^2\max_{1\le i,j\le q}\mathbb{P}\Big(\Big|\frac{m^2(n_1-1)}{(n_1-m)^2}\sum_{\alpha=1}^{n}(\widetilde{q}_{1\alpha,ij}-\widetilde{u}_{1,ij})^2 - m^2\zeta_{1,ij}\Big| \ge C\frac{\varepsilon_n}{\log q}\Big).\end{aligned}$$

Considering that $\zeta_{1,ij}$ is the variance of $h_{ij}(\boldsymbol{X}_\ell)$, by setting $\bar{h}_{1,ij} := \sum_{\alpha=1}^{n_1} h_{ij}(\boldsymbol{X}_\alpha)/n_1$, we use $\sum_{\alpha=1}^{n_1}(h_{ij}(\boldsymbol{X}_\alpha)-\bar{h}_{1,ij})^2/n_1$ to approximate $\zeta_{1,ij}$. Therefore, we insert the term $\sum_{\alpha=1}^{n_1}(h_{ij}(\boldsymbol{X}_\alpha)-\bar{h}_{1,ij})^2/n_1$ and use the triangle inequality to obtain

$$\begin{aligned}&q^2\max_{1\le i,j\le q}\mathbb{P}\Big(\Big|\frac{m^2(n_1-1)}{(n_1-m)^2}\sum_{\alpha=1}^{n}(\widetilde{q}_{1\alpha,ij}-\widetilde{u}_{1,ij})^2 - m^2\zeta_{1,ij}\Big| \ge C\frac{\varepsilon_n}{\log q}\Big)\\ &\le \underbrace{q^2\max_{1\le i,j\le q}\mathbb{P}\Big(\Big|\frac{m^2(n_1-1)}{(n_1-m)^2}\Big(\sum_{\alpha=1}^{n}(\widetilde{q}_{1\alpha,ij}-\widetilde{u}_{1,ij})^2 - \sum_{\alpha=1}^{n_1}(h_{ij}(\boldsymbol{X}_\alpha)-\bar{h}_{1,ij})^2\Big)\Big| \ge \frac{1}{2}C\frac{\varepsilon_n}{\log q}\Big)}_{A_1}\\ &\quad+\underbrace{q^2\max_{1\le i,j\le q}\mathbb{P}\Big(\Big|\frac{m^2(n_1-1)}{(n_1-m)^2}\sum_{\alpha=1}^{n_1}(h_{ij}(\boldsymbol{X}_\alpha)-\bar{h}_{1,ij})^2 - m^2\zeta_{1,ij}\Big| \ge \frac{1}{2}C\frac{\varepsilon_n}{\log q}\Big)}_{A_2}.\end{aligned}$$

We then introduce an additional lemma to complete the proof.

**Lemma B.1.** Under Assumption **(A2)**, as $n$, $q \to \infty$, we have

$$q^2 \max_{1\le i,j\le q} \mathbb{P}\Big(\Big|\frac{m^2(n_1-1)}{(n_1-m)^2}\big(\sum_{\alpha=1}^{n}(\widetilde{q}_{1\alpha,ij}-\widetilde{u}_{1,ij})^2-\sum_{\alpha=1}^{n_1}(h_{ij}(\boldsymbol{X}_\alpha)-\bar{h}_{1,ij})^2\big)\Big|\ge C\frac{\varepsilon_n}{\log q}\Big)=o(1), \quad \text{(B.3)}$$

$$q^2 \max_{1\le i,j\le q} \mathbb{P}\Big(\Big|\frac{m^2(n_1-1)}{(n_1-m)^2}\sum_{\alpha=1}^{n_1}(h_{ij}(\boldsymbol{X}_\alpha)-\bar{h}_{1,ij})^2-m^2\zeta_{1,ij}\Big|\ge C\frac{\varepsilon_n}{\log q}\Big)=o(1), \quad \text{(B.4)}$$

where $\varepsilon_n = o(1)$. Results for $\boldsymbol{Y}$ are the same.

The detailed proof of Lemma B.1 is in Supplement C.1. By Lemma B.1, we have $A_1 = o(1)$ and $A_2 = o(1)$. Hence, we prove Lemma A.1. □

## B.2. Proof of Lemma A.2 of Appendix A

***Proof.*** In Lemma A.2, we aim to prove (A.7). For notational simplicity, we set

$$W_{ij} = \frac{\binom{n_1}{m}^{-1}\Delta_{n_1,ij} - \binom{n_2}{m}^{-1}\Delta_{n_2,ij}}{\sqrt{m^2\zeta_{1,ij}/n_1 + m^2\zeta_{2,ij}/n_2}} \qquad \text{and} \qquad W_{a,ij} = \sqrt{n_a}\binom{n_a}{m}^{-1}\Delta_{n_a,ij}. \quad \text{(B.5)}$$

By (A.6), we have $W_{ij} = \widetilde{N}_{ij} - T_{ij}$. We also have $W_{ij} = c_1 W_{1,ij} + c_2 W_{2,ij}$, where

$$c_1 = 1/\sqrt{m^2\zeta_{1,ij} + m^2\zeta_{2,ij}n_1/n_2} \qquad \text{and} \qquad c_2 = -1/\sqrt{m^2\zeta_{1,ij}n_2/n_1 + m^2\zeta_{2,ij}}.$$

To prove (A.6), by setting $L_2 := |\max_{1\le i,j\le q}(\widetilde{N}_{ij})^2 - \max_{1\le i,j\le q}(T_{ij})^2|$, it suffices to prove that as $n$, $q \to \infty$ we have $L_2 = o_p(1)$. By the definitions of $\widetilde{N}_{ij}$ and $W_{ij}$ in (A.4) and (B.5), we obtain

$$L_2 \le \max_{1\le i,j\le q}|\widetilde{N}_{ij} - T_{ij}| \max_{1\le i,j\le q}|\widetilde{N}_{ij} + T_{ij}| \le \max_{1\le i,j\le q}|W_{ij}|(\max_{1\le i,j\le q} 2|T_{ij}| + \max_{1\le i,j\le q}|W_{ij}|).$$

Considering $\max_{1\le i,j\le q} 2|T_{ij}| + \max_{1\le i,j\le q}|W_{ij}| = O_p(\log q)$, to obtain $L_2 = o_p(1)$, we only need to prove

$$\mathbb{P}\Big(\max_{1\le i,j\le q}|W_{ij}| > 1/\log q\Big) \le q^2 \max_{1\le i,j\le q}\mathbb{P}\big(|W_{ij}| > 1/\log q\big) \to 0, \quad \text{(B.6)}$$

as $n$, $q \to \infty$. By Proposition 2.3(c) of Arcones and Gine (1993) and $W_{ij} = \widetilde{N}_{ij} - T_{ij}$, we have

$$\mathbb{P}\Big(|W_{ij}| > \frac{1}{\log q}\Big) \le C\exp\Bigg(-\frac{C_1(n^{m/2-1/2}/\log q)^{2/m}}{\sigma^{2/m} + C_2\big((n^{m/2-1/2}/\log q)^{1/m}n^{-1/2}\big)^{2/(m+1)}}\Bigg). \quad \text{(B.7)}$$

Combining $\log q = O(n^{1/3-\epsilon})$ (see Assumption **(A2)**) and (B.7), we obtain (B.6). We then have $L_2 = o_p(1)$. This completes the proof. □

## B.3. Proof of Lemma A.3 of Appendix A

***Proof.*** In Lemma A.3, we aim to prove that as $n,\ q \to \infty$, we have $\mathbb{P}\big(\max_{(i,j)\in S_0}(T_{ij})^2 \geq y_q\big) \to 0$, where $T_{ij}$ is defined in (A.5) and $y_q = x + 4\log q - \log(\log q)$. For notational simplicity, we introduce

$$g_{ij}(\boldsymbol{X}_\alpha) := \mathbb{E}[\Phi_{ij}(\boldsymbol{X}_{\alpha_1},\ldots,\boldsymbol{X}_{\alpha_m})|\boldsymbol{X}_\alpha] \quad \text{and} \quad g_{ij}(\boldsymbol{Y}_\alpha) := \mathbb{E}[\Phi_{ij}(\boldsymbol{Y}_{\alpha_1},\ldots,\boldsymbol{Y}_{\alpha_m})|\boldsymbol{Y}_\alpha],$$

where $\Phi_{ij}$ is the kernel function of the U-statistic $\widehat{u}_{a,ij}$. By Lemma D.3, we have $\zeta_{1,ij} = \mathrm{Var}(g_{ij}(\boldsymbol{X}_\alpha))$ and $\zeta_{2,ij} = \mathrm{Var}(g_{ij}(\boldsymbol{Y}_\alpha))$. We then introduce an additional lemma.

**Lemma B.2.** If $g_{ij}$ is bounded, under Assumption **(A2)**, we have

$$\mathbb{P}\bigg(\max_{(i,j)\in\Lambda}\frac{(\frac{1}{n_1}\sum_{\alpha=1}^{n_1}g_{ij}(\boldsymbol{X}_\alpha)-\frac{1}{n_2}\sum_{\alpha=1}^{n_2}g_{ij}(\boldsymbol{Y}_\alpha)-u_{1,ij}+u_{2,ij})^2}{\zeta_{1,ij}/n_1+\zeta_{2,ij}/n_2} > t^2\bigg) \leq C_1|\Lambda|(1-\Phi(t)),$$

uniformly for any $\Lambda \subseteq S$ and $t \in [0, O(n^{1/6-\varepsilon})]$, where $\varepsilon$ is an arbitrary positive number.

The detailed proof of Lemma B.2 is in Supplement C.2. Considering $\log q = O(n^{1/3-\epsilon})$ (see Assumption **(A2)**), by choosing $\varepsilon < \epsilon/2$ with $\varepsilon$ small enough, we have $\sqrt{y_q} \in [0, O(n^{1/6-\varepsilon})]$. By Lemma D.5, we have $1-\Phi(t) \leq C\exp(-C_1t^2)$ as $\Phi$ is the cumulative distribution function of a standard normal distribution. Considering $|S_0| = o(q^2)$ (see Assumption **(A1)**) and $y_q = x + 4\log q - \log(\log q)$, by setting $\Lambda$ in Lemma B.2 as $S_0$, we obtain

$$\mathbb{P}\Big(\max_{(i,j)\in S_0}(T_{ij})^2 \geq y_q\Big) \leq o(q^2)\exp(-C_1y_q) = o(1).$$

This completes the proof. □

## B.4. Proof of Lemma A.4 of Appendix A

***Proof.*** To prove Lemma A.4, we obtain both upper and lower bounds of $\mathbb{P}\big(\max_{1\leq k\leq h}(T_k)^2 \geq y_q\big)$ by using the Boferroni inequality and normal approximation. By the Boferroni inequality (Lemma 1 of Cai et al. (2013)), for any integer $M$ with $0 < M < [h/2]$, we have

$$\sum_{\ell=1}^{2M}(-1)^{\ell-1}\sum_{1\leq k_1<\cdots<k_\ell\leq h}\mathbb{P}\Big(\bigcap_{j=1}^{\ell}E_{k_j}\Big)\leq\mathbb{P}\Big(\max_{1\leq k\leq h}(T_k)^2 \geq y_q\Big)\leq\sum_{\ell=1}^{2M-1}(-1)^{\ell-1}\sum_{1\leq k_1<\cdots<k_\ell\leq h}\mathbb{P}\Big(\bigcap_{j=1}^{\ell}E_{k_j}\Big), \quad \text{(B.8)}$$

where we set $E_{k_j} = \big\{(T_{k_j})^2 \geq y_q\big\}$. By the definition of $T_k$ in (A.12), we have

$$T_k = T_{i_kj_k} = \sum_{\beta=1}^{n_1+n_2}\widehat{Z}_{\beta,i_kj_k}\Big/\sqrt{n_2^2\zeta_{1,i_kj_k}/n_1 + n_2\zeta_{2,i_kj_k}}.$$

For any $\ell$-tuple $k_1,\ldots,k_\ell$ satisfying $1 \leq k_1 < \ldots < k_\ell \leq h$, we have

$$\widetilde{Z}_{\beta k} = \widehat{Z}_{\beta,i_kj_k}/(n_2\zeta_{1,i_kj_k}/n_1 + \zeta_{2,i_kj_k})^{1/2} \qquad \text{and} \qquad \boldsymbol{W}_\beta = (\widetilde{Z}_{\beta k_1},\ldots,\widetilde{Z}_{\beta k_\ell})^T,$$

for $1 \leq k \leq h$ and $1 \leq \beta \leq n_1 + n_2$. Therefore, $T_k = (n_2)^{-1/2} \sum_{\beta=1}^{n_1+n_2} \widetilde{Z}_{\beta k}$. Define $\|\mathbf{v}\|_{\min} = \min_{1\leq i\leq \ell} |v_i|$ for any vector $\mathbf{v} \in \mathbb{R}^\ell$. By $E_{k_j} = \{(T_{k_j})^2 \geq y_q\}$, we have

$$\mathbb{P}\Big(\bigcap_{j=1}^{\ell} E_{k_j}\Big) = \mathbb{P}\Big(\|n_2^{-1/2} \sum_{\beta=1}^{n_1+n_2} \boldsymbol{W}_\beta\|_{\min} \geq y_q^{1/2}\Big).$$

We set $\boldsymbol{N}_\ell$ as a normal vector with the same mean vector and covariance matrix as $n_2^{-1/2} \sum_{\beta=1}^{n_1+n_2} \boldsymbol{W}_\beta$. More specifically, $\boldsymbol{N}_\ell := (N_{k_1}, \dots, N_{k_\ell})^T$ is a normal vector with $\mathbb{E}[\boldsymbol{N}_\ell] = 0$ and $\mathrm{Var}(\boldsymbol{N}_\ell) = n_1 \mathrm{var}(\boldsymbol{W}_1)/n_2 + \mathrm{Var}(\boldsymbol{W}_{n_1+1})$. Therefore, $\boldsymbol{N}_\ell$ is the normal approximation of $n_2^{-1/2} \sum_{\beta=1}^{n_1+n_2} \boldsymbol{W}_\beta$. We then aim to use $\boldsymbol{N}_\ell$ to rewrite obtained lower and upper bounds in (B.8). As $|\widetilde{Z}_{\beta k}|$ is bounded, we set $|\widetilde{Z}_{\beta k}| \leq K$, where $K$ is a constant. By Theorem 1 of Zaitsev (1987), we have

$$\begin{aligned}
&\mathbb{P}\Big(\|n_2^{-1/2} \sum_{\beta=1}^{n_1+n_2} \boldsymbol{W}_\beta\|_{\min} \geq y_q^{1/2}\Big) \\
&\quad \leq \mathbb{P}\Big(\|\boldsymbol{N}_\ell\|_{\min} \geq y_q^{1/2} - \epsilon_n (\log q)^{-1/2}\Big) + C_1 \ell^{5/2} \exp\Big(-\frac{n^{1/2}\epsilon_n}{C_2 \ell^3 K (\log q)^{1/2}}\Big),
\end{aligned}$$

where $\epsilon_n \to 0$, which will be specified later. Considering that $\ell$ is a fixed integer that does not depend on $n$ and $q$, by $\log q = O(n^{1/3-\epsilon})$ (see Assumption **(A2)**), we can let $\epsilon_n \to 0$ sufficiently slow such that we have

$$c_1 \ell^{5/2} \exp\Big(-\frac{n^{1/2}\epsilon_n}{c_2 \ell^3 K (\log q)^{1/2}}\Big) = O(q^{-J}),$$

for any large $J > 0$. Hence, we prove (A.16). The proof of (A.17) is similar. This completes the proof of Lemma A.4. □

## B.5. Proof of Lemma A.5 of Appendix A

***Proof.*** We aim to prove (A.18) under Assumption **(A3)**. $g_{ij}$ is defined in (2.9). Assumption **(A3)** requires

$$u_{1,ijk\ell} = \mathbb{E}[g_{ij}(\boldsymbol{X}_\alpha) g_{k\ell}(\boldsymbol{X}_\alpha)] = O((\log q)^{-1-\alpha_0}), \tag{B.9}$$

$$u_{2,ijk\ell} = \mathbb{E}[g_{ij}(\mathbf{Y}_\alpha) g_{k\ell}(\mathbf{Y}_\alpha)] = O((\log q)^{-1-\alpha_0}), \tag{B.10}$$

for $(i,j),(k,\ell) \in S \setminus S_0$, where $S = \{(i,j) : 1 \leq i,j \leq q\}$ and $S_0$ is defined in (2.11). By definitions of $S$ and $S_0$, for $(i,j) \in S \setminus S_0$, we have $|u_{a,ij}| \leq O((\log q)^{-1-\alpha_0})$. $\boldsymbol{N}_\ell$ is an $\ell$-dimensional normal random vector with zero mean and covariance matrix $\mathbf{V}_\ell$. Combining (2.9), (A.36), (B.9) and (B.10), we obtain $\|\mathbf{V}_\ell - \mathbf{I}_\ell\|_2 = O\big((\log q)^{-1-\alpha_0}\big)$ uniformly for all $k_1, k_2, \ldots, k_\ell$, where $k_1, \ldots, k_\ell$ is any $\ell-$tuple of positive integers satisfying $1 \le k_1 < \ldots < k_\ell \le h$. By setting $w_q = y_q^{1/2} \pm \epsilon_n (\log q)^{-1/2}$, we have

$$\begin{aligned}\mathbb{P}\big(\|\boldsymbol{N}_\ell\|_{\min} \ge y_q^{1/2} \pm \epsilon_n(\log q)^{-1/2}\big) =& \mathbb{P}\big(|N_{k_1}| \ge w_q, |N_{k_2}| \ge w_q, \ldots, |N_{k_\ell}| \ge w_q\big) \\ =& \frac{1}{(2\pi)^{\ell/2}|\mathbf{V}_\ell|^{1/2}} \int_{\|\mathbf{x}\|_{\min} \ge w_q} \exp\big(-\frac{1}{2}\mathbf{x}^T\mathbf{V}_\ell^{-1}\mathbf{x}\big) d\mathbf{x}.\end{aligned}$$

We then calculate the integral to obtain

$$\begin{aligned}&\frac{1}{(2\pi)^{\ell/2}|\mathbf{V}_\ell|^{1/2}} \int_{\|\mathbf{x}\|_{\min} \ge w_q} \exp\Big(-\frac{1}{2}\mathbf{x}^T\mathbf{V}_\ell^{-1}\mathbf{x}\Big) d\mathbf{x} \\ =& \frac{1}{(2\pi)^{\ell/2}|\mathbf{V}_\ell|^{1/2}} \int_{\substack{\|\mathbf{x}\|_{\min} \ge w_q, \\ \|\mathbf{x}\|_{\max} \le (\log q)^{1/2+\alpha_0/4}}} \exp\Big(-\frac{1}{2}\mathbf{x}^T\mathbf{V}_\ell^{-1}\mathbf{x}\Big) d\mathbf{x} + O(\exp(-(\log q)^{1+\alpha_0/2}/4)).\end{aligned}$$

Considering that $\|\mathbf{V}_\ell - \mathbf{I}_\ell\|_2 = O\big((\log q)^{-1-\alpha_0}\big)$ holds uniformly for all $(k_1, k_2, \ldots, k_\ell)$, we obtain

$$\begin{aligned}&\frac{1}{(2\pi)^{\ell/2}|\mathbf{V}_\ell|^{1/2}} \int_{\substack{\|\mathbf{x}\|_{\min} \ge w_q, \\ \|\mathbf{x}\|_{\max} \le (\log q)^{1/2+\alpha_0/4}}} \exp\Big(-\frac{1}{2}\mathbf{x}^T\mathbf{V}_\ell^{-1}\mathbf{x}\Big) d\mathbf{x} + O(\exp(-(\log q)^{1+\alpha_0/2}/4)) \\ =& \frac{1 + O((\log q)^{-\alpha_0/2})}{(2\pi)^{\ell/2}} \int_{\substack{\|\mathbf{x}\|_{\min} \ge w_q, \\ \|\mathbf{x}\|_{\max} \le (\log q)^{1/2+\alpha_0/4}}} \exp\Big(-\frac{1}{2}\mathbf{x}^T\mathbf{x}\Big) d\mathbf{x} + O\Big(\exp\big(-(\log q)^{1+\alpha_0/2}/4\big)\Big) \\ =& \frac{1 + O((\log q)^{-\alpha_0/2})}{(2\pi)^{\ell/2}} \int_{\|\mathbf{x}\|_{\min} \ge w_q} \exp\Big(-\frac{1}{2}\mathbf{x}^T\mathbf{x}\Big) d\mathbf{x} + O\Big(\exp\big(-(\log q)^{1+\alpha_0/2}/4\big)\Big).\end{aligned}$$

We then calculate the integral to get

$$\begin{aligned}&\frac{1 + O((\log q)^{-\alpha_0/2})}{(2\pi)^{\ell/2}} \int_{\|\mathbf{x}\|_{\min} \ge w_q} \exp\Big(-\frac{1}{2}\mathbf{x}^T\mathbf{x}\Big) d\mathbf{x} + O\Big(\exp\big(-(\log q)^{1+\alpha_0/2}/4\big)\Big) \\ &= \big(1 + o(1)\big)\Big(\frac{2}{\sqrt{8\pi}}\exp(-\frac{x}{2})\Big)^\ell q^{-2\ell}.\end{aligned}$$

Therefore, as $n, q \to \infty$, we have

$$\mathbb{P}\big(\|\boldsymbol{N}_\ell\|_{\min} \ge y_q^{1/2} \pm \epsilon_n(\log q)^{-1/2}\big) = \big(1 + o(1)\big)\Big(\frac{2}{\sqrt{8\pi}}\exp(-\frac{x}{2})\Big)^\ell q^{-2\ell}. \quad \text{(B.11)}$$

Considering $C_h^\ell = h!/(\ell!(h-\ell)!)$ and $2h/(q^2) = 1 + o(1)$, by (B.11), we prove (A.18). Hence, we complete the proof. □

## B.6. Proof of Lemma A.1 of Supplementary Material

***Proof.*** By Theorem 2.2, under Assumptions **(A1)**, **(A2)** and **(A3)**, we have

$$\mathbb{P}\Big(M_n^1 - 4\log q + \log(\log q) \le x\Big) \to \exp\Big(-\frac{1}{\sqrt{8\pi}}\exp(-\frac{x}{2})\Big). \quad \text{(B.12)}$$

However, in Lemma A.1, we don't assume **(A1)** and **(A3)**, so that we cannot get (A.5) from (B.12).

In Lemma A.1, we aim to prove (B.12) under Assumption **(A2)**. By (A.3) and (A.7), we have

$$\left| M_n^1 - \max_{(i,j)\in S} \frac{(\frac{1}{n_1}\sum_{\alpha=1}^{n_1} g_{ij}(\boldsymbol{X}_\alpha) - \frac{1}{n_2}\sum_{\alpha=1}^{n_2} g_{ij}(\boldsymbol{Y}_\alpha) - u_{1,ij} + u_{2,ij})^2}{\zeta_{1,ij}/n_1 + \zeta_{2,ij}/n_2} \right| = o_p(1). \quad \text{(B.13)}$$

By Lemma B.2, we have

$$\mathbb{P}\left( \max_{(i,j)\in S} \frac{(\frac{1}{n_1}\sum_{\alpha=1}^{n_1} g_{ij}(\boldsymbol{X}_\alpha) - \frac{1}{n_2}\sum_{\alpha=1}^{n_2} g_{ij}(\boldsymbol{Y}_\alpha) - u_{1,ij} + u_{2,ij})^2}{\zeta_{1,ij}/n_1 + \zeta_{2,ij}/n_2} > t^2 \right) \le C_1|S|(1-\Phi(t)), \quad \text{(B.14)}$$

where $\Phi(t)$ is the cumulative distribution function of a standard normal distribution. By Lemma D.5, we have $1-\Phi(t) \le \exp(-t^2/2)/(t\sqrt{2\pi})$ for $t > 0$. Considering $|S| = q^2$ and (B.14), by setting $t^2 = 4\log q - 0.5\log(\log q)$, under Assumption **(A2)**, we get

$$\mathbb{P}\left( \max_{(i,j)\in S} \frac{(\frac{1}{n_1}\sum_{\alpha=1}^{n_1} g_{ij}(\boldsymbol{X}_\alpha) - \frac{1}{n_2}\sum_{\alpha=1}^{n_2} g_{ij}(\boldsymbol{Y}_\alpha) - u_{1,ij} + u_{2,ij})^2}{\zeta_{1,ij}/n_1 + \zeta_{2,ij}/n_2} > 4\log q - \frac{1}{2}\log(\log q) \right) = o(1). \quad \text{(B.15)}$$

Combining (B.13) and (B.15), we prove (A.5). This completes the proof. □

## B.7. Proof of Lemma A.2 of Supplementary Material

***Proof.*** By setting $S = \{(i,j) : 1 \le i,j \le q\}$, Lemma A.1 implies that both of the following two events

$$\mathcal{E}_1 := \left\{ \max_{(i,j)\in S} \left| n_1\widehat{\sigma}^2(\widehat{u}_{1,ij}) - m^2\zeta_{1,ij} \right| < C\frac{\varepsilon_n}{\log q} \right\},$$
$$\mathcal{E}_2 := \left\{ \max_{(i,j)\in S} \left| n_2\widehat{\sigma}^2(\widehat{u}_{2,ij}) - m^2\zeta_{2,ij} \right| < C\frac{\varepsilon_n}{\log q} \right\},$$

happen with probability going to one as $n, q \to \infty$. Under $\mathcal{E}_1$ and $\mathcal{E}_2$, by $\zeta_{a,ij} \ge r_a > 0$ (see Assumption **(A2)**), we have

$$\left| n_1\widehat{\sigma}^2(u_{1,ij})/(m^2\zeta_{1,ij}) - 1 \right| < C\varepsilon_n/\log q \quad \text{and} \quad \left| n_2\widehat{\sigma}^2(u_{2,ij})/(m^2\zeta_{2,ij}) - 1 \right| < C\varepsilon_n/\log q.$$

Therefore, under $\mathcal{E}_1$ and $\mathcal{E}_2$, we have

$$\begin{aligned} &\left| \max_{(i,j)\in S} \frac{(u_{1,ij} - u_{2,ij})^2}{\widehat{\sigma}^2(\widehat{u}_{1,ij}) + \widehat{\sigma}^2(\widehat{u}_{2,ij})} - \max_{(i,j)\in S} \frac{(u_{1,ij} - u_{2,ij})^2}{m^2\zeta_{1,ij}/n_1 + m^2\zeta_{2,ij}/n_2} \right| \\ &\quad \le \max_{(i,j)\in S} \left| \frac{(u_{1,ij} - u_{2,ij})^2}{\widehat{\sigma}^2(\widehat{u}_{1,ij}) + \widehat{\sigma}^2(\widehat{u}_{2,ij})} - \frac{(u_{1,ij} - u_{2,ij})^2}{m^2\zeta_{1,ij}/n_1 + m^2\zeta_{2,ij}/n_2} \right| \le C\frac{\varepsilon_n}{\log q}. \end{aligned}$$

Hence, as $n,\ q \to \infty$, we have

$$\Big| \max_{(i,j)\in S} \frac{(u_{1,ij}-u_{2,ij})^2}{\widehat{\sigma}^2(\widehat{u}_{1,ij})+\widehat{\sigma}^2(\widehat{u}_{2,ij})} - \max_{(i,j)\in S} \frac{(u_{1,ij}-u_{2,ij})^2}{m^2\zeta_{1,ij}/n_1+m^2\zeta_{2,ij}/n_2}\Big| = o_p(1). \qquad \text{(B.16)}$$

Considering $(\mathbf{U}_1^\tau, \mathbf{U}_2^\tau) \in \mathbb{A}(4)$, by (B.16), we obtain (A.6). This completes the proof. □

## B.8. Proof of Lemma A.3 of of Supplementary Material

***Proof.*** By the definition of $N_{ij}^{\tau,\mathrm{ps}}$ in (A.3) and (3.8), we have $M_n^{\tau,\mathrm{ps}} = \max_{(i,j)\in S}(N_{ij}^{\tau,\mathrm{ps}})^2$. By setting $S = \{(i,j) : 1 \le i < j \le d\}$, we aim to prove that as $n$ , $d \to \infty$, we have

$$\Big| \max_{(i,j)\in S}(N_{ij}^{\tau,\mathrm{ps}})^2 - \max_{(i,j)\in S}(\widehat{T}_{ij})^2\Big| = o_p(1), \qquad \text{(B.17)}$$

where we define $\widehat{T}_{ij}$ in (A.22). By setting $L_2^{\tau,\mathrm{ps}} := |\max_{(i,j)\in S}(N_{ij}^{\tau,\mathrm{ps}})^2 - \max_{(i,j)\in S}(\widehat{T}_{ij})^2|$, we aim to prove $L_2^{\tau,\mathrm{ps}} = o_p(1)$. To prove $L_2^{\tau,\mathrm{ps}} = o_p(1)$, we construct an upper bound of $L_2^{\tau,\mathrm{ps}}$ as

$$L_2^{\tau,\mathrm{ps}} \le \max_{(i,j)\in S}|N_{ij}^{\tau,\mathrm{ps}} - \widehat{T}_{ij}| \max_{(i,j)\in S}|N_{ij}^{\tau,\mathrm{ps}} + \widehat{T}_{ij}| \le \max_{(i,j)\in S}|\widehat{W}_{ij}|(\max_{(i,j)\in S} 2|\widehat{T}_{ij}| + \max_{(i,j)\in S}|\widehat{W}_{ij}|).$$

Considering $\max_{(i,j)\in S} 2|\widehat{T}_{ij}| + \max_{(i,j)\in S}|\widehat{W}_{ij}| = O_p(\log d)$, to show $L_2^{\tau,\mathrm{ps}} = o_p(1)$, we only need to prove

$$\mathbb{P}\Big(\max_{(i,j)\in S}|\widehat{W}_{ij}| > 1/\log d\Big) \le d^2 \max_{(i,j)\in S}\mathbb{P}\big(|\widehat{W}_{ij}| > 1/\log d\big) \to 0, \qquad \text{(B.18)}$$

as $n,\ d \to \infty$. By Proposition 2.3(c) of Arcones and Gine (1993), we have

$$\mathbb{P}\Big(|\widehat{W}_{ij}| > \frac{1}{\log d}\Big) \le C \exp\Big(-\frac{C_2\sqrt{n}/\log d}{\sigma + C_1(\log d\sqrt{n})^{-1/3}}\Big). \qquad \text{(B.19)}$$

Considering $\log d = O(n^{1/3-\epsilon})$ (see Assumption **(A2)**), by (B.18) and (B.19), we obtain

$$\mathbb{P}\Big(\max_{(i,j)\in S}|\widehat{W}_{ij}| > 1/\log d\Big) \to 0,$$

as $n,\ d \to \infty$. Hence, we complete the proof. □

## B.9. Proof of Lemma A.4 of Supplementary Material

***Proof.*** $\sigma_{a,\mathrm{ps}}^2$ and $g_{ij}^\tau$ are defined in (3.7) and (A.17). We aim to prove (A.24). For this, we introduce the following lemma.

**Lemma B.3.** By setting $S = \{(i,j) : 1 \le i < j \le d\}$, under Assumption **(A2)**, we have

$$\mathbb{P}\Big(\max_{(i,j)\in\Lambda} \frac{(\frac{1}{n_1}\sum_{\alpha=1}^{n_1} g_{ij}^{\tau}(\boldsymbol{X}_\alpha) - \frac{1}{n_2}\sum_{\alpha=1}^{n_2} g_{ij}^{\tau}(\boldsymbol{Y}_\alpha) - \tau_{1,ij} + \tau_{1,ij})^2}{\sigma_{1,\mathrm{ps}}^2/(4n_1) + \sigma_{2,\mathrm{ps}}^2/(4n_2)} > t^2\Big) \le C|\Lambda|(1-\Phi(C_1 t)),$$

uniformly for any $\Lambda \subseteq S$ and $t \in [0, O(n^{1/6-\varepsilon})]$, where $\varepsilon$ is an arbitrary positive number. $\Phi(t)$ is the cumulative distribution function of a standard normal distribution.

The detailed proof of Lemma B.3 is in Supplement C.3. By Lemma B.3, we have

$$\mathbb{P}\Big(\max_{(i,j)\in S_0} (\widehat{T}_{ij})^2 > y_d\Big) \le C|S_0|\big(1-\Phi(C_1\sqrt{y_d})\big), \tag{B.20}$$

where $y_d = x + 4\log d - \log(\log d)$. By Lemma D.5, we have

$$1-\Phi(t) \le \exp(-t^2/2)/(t\sqrt{2\pi}) \qquad \text{for} \quad t > 0. \tag{B.21}$$

Considering $|S_0| = o(d^2)$ (see Assumption **A1**), by (B.20) and (B.21), we have

$$\mathbb{P}\Big(\max_{(i,j)\in S_0} (\widehat{T}_{ij})^2 > y_d\Big) = o(1)$$

. This completes the proof. □

## B.10. Proof of Lemma A.5 of Supplementary Material

***Proof.*** In Lemma A.5, we aim to prove $|4\zeta_{a,ij} - \sigma_{\mathrm{ps}}^2| \le C|\tau_{a,ij}|$. For simplicity, we ignore the subscript $a$ ($a = 1$ or 2) for $\boldsymbol{X}$ and $\boldsymbol{Y}$. By (A.9), we have $4\zeta_{ij} = 16(\Pi_{cc,ij} - (\Pi_{c,ij})^2)$, where $\Pi_{cc,ij}$ and $\Pi_{c,ij}$ are defined in (3.4) and (3.2). By the definition of the meta-elliptical distribution in Definition E.2, we have

$$(X_{\alpha i}, X_{\alpha j}) \overset{d}{=} (R_\alpha \cos(\Theta_\alpha), R_\alpha \sin(\Theta_\alpha + u)), \tag{B.22}$$

where $R_\alpha$ is a positive radial random variable and $\Theta_\alpha$ is an independent, uniformly distributed angle on $[-\pi, \pi]$ with $\alpha = 1, 2, \cdots, n_1$ and $u = \tau_{ij}\pi/2$. By the definition of $\sigma_{\mathrm{ps}}^2$ in (3.7), we have $\sigma_{\mathrm{ps}}^2 = 16(\Pi_{cc,ij}^{u=0} - (\Pi_{c,ij}^{u=0})^2)$, where $\Pi_{cc,ij}^{u=0}$ and $\Pi_{c,ij}^{u=0}$ is the $\Pi_{cc,ij}$ and $\Pi_{c,ij}$ under the condition $u = \tau_{ij}\pi/2 = 0$.

By the definition of $\Pi_{cc,ij}$ in (3.4), $\Pi_{cc,ij}$ is sum of the probability of four parts:

$$\begin{aligned}
E_1 &:= \{X_{2i} > X_{1i}, X_{2j} > X_{1j}, X_{3i} > X_{1i}, X_{3j} > X_{1j}\},\\
E_2 &:= \{X_{2i} < X_{1i}, X_{2j} < X_{1j}, X_{3i} < X_{1i}, X_{3j} < X_{1j}\},\\
E_3 &:= \{X_{2i} > X_{1i}, X_{2j} > X_{1j}, X_{3i} < X_{1i}, X_{3j} < X_{1j}\},\\
E_4 &:= \{X_{2i} < X_{1i}, X_{2j} < X_{1j}, X_{3i} > X_{1i}, X_{3j} > X_{1j}\},
\end{aligned}$$

i.e., $\Pi_{cc,ij} = \mathbb{P}(E_1) + \mathbb{P}(E_2) + \mathbb{P}(E_3) + \mathbb{P}(E_4)$. By (B.22), we rewrite the event $E_1$ as

$$E_1 = \{R_2\cos(\Theta_2) > X_{1i}, R_2\sin(\Theta_2+u) > X_{1j}, R_3\cos(\Theta_3) > X_{1i}, R_3\sin(\Theta_3+u) > X_{1j}\}.$$

We then define $g_u(x_{1i}, x_{1j}, r_2, r_3)$ as

$$g_u(x_{1i}, x_{1j}, r_2, r_3) := \mathbb{P}[E_1 | X_{1i} = x_{1i}, X_{1j} = x_{1j}, R_2 = r_2, R_3 = r_3]. \tag{B.23}$$

By $\mathbb{P}(E_1) = \mathbb{E}\Big[\mathbb{P}\big[E_1 | X_{1i}, X_{1j}, R_2, R_3\big]\Big]$, we have $\mathbb{P}(E_1) = \mathbb{E}[g_u(X_{1i}, X_{1j}, R_2, R_3)]$. By (B.22) and (B.23), we have

$$|g_u(x_{1i}, x_{1j}, r_2, r_3) - g_{u=0}(x_{1i}, x_{1j}, r_2, r_3)| \le C\tau_{ij}, \tag{B.24}$$

where $C$ is a positive constant. By (B.24), we obtain

$$\Big|\mathbb{E}\big[g_u(X_{1i}, X_{1j}, R_2, R_3)\big] - \mathbb{E}\big[g_{u=0}(X_{1i}, X_{1j}, R_2, R_3)\big]\Big| \le C|\tau_{ij}|.$$

By similar proofs on $E_2$, $E_3$ and $E_4$, we have $|\Pi_{cc,ij} - \Pi_{cc,ij}^{u=0}| \le C|\tau_{ij}|$. Considering the boundedness of $\Pi_{c,ij}$, we also have $|(\Pi_{c,ij})^2 - (\Pi_{c,ij}^{u=0})^2| \le C|\tau_{ij}|$. We then complete the proof by $\sigma_{\mathrm{ps}}^2 = 16(\Pi_{cc,ij}^{u=0} - (\Pi_{c,ij}^{u=0})^2)$ and $4\zeta_{ij} = 16(\Pi_{cc,ij} - (\Pi_{c,ij})^2)$. □

## B.11. Proof of Lemma A.6

***Proof.*** To prove Lemma A.6, we set

$$M_n^{1,\mathrm{plug}} := \max_{(i,j)\in S} \frac{(\widehat{\tau}_{1,ij} - \widehat{\tau}_{2,ij} - \tau_{1,ij} + \tau_{2,ij})^2}{\widehat{\sigma}^2_{\mathrm{plug}}(\widehat{\tau}_{1,ij}) + \widehat{\sigma}^2_{\mathrm{plug}}(\widehat{\tau}_{2,ij})},$$
$$\widetilde{M}_n^{1,\mathrm{plug}} := \max_{(i,j)\in S} \frac{(\widehat{\tau}_{1,ij} - \widehat{\tau}_{2,ij} - \tau_{1,ij} + \tau_{2,ij})^2}{4\zeta_{1,ij}/n_1 + 4\zeta_{2,ij}/n_2},$$

where $S = \{(i,j) : 1 \le i < j \le d\}$. By Theorem 3.2, under Assumptions **(A1)**, **(A2)** and **(A3)**, we have

$$\mathbb{P}\Big(M_n^{1,\mathrm{plug}} - 4\log d + \log(\log d) \le x\Big) \to \exp\Big(-\frac{1}{\sqrt{8\pi}}\exp(-\frac{x}{2})\Big). \tag{B.25}$$

However, in Lemma A.6, we don't assume **(A1)** and **(A3)**. Therefore, we cannot obtain (A.28) from (B.25).

In Lemma A.6, we aim to prove (A.28) under Assumption **(A2)**. By (A.7) and (A.7), as $n,\ d \to \infty$, we obtain

$$\Bigg|M_n^{1,\mathrm{plug}} - \max_{(i,j)\in S} \frac{(\frac{1}{n_1}\sum_{\alpha=1}^{n_1} g_{ij}^{\tau}(\boldsymbol{X}_\alpha) - \frac{1}{n_2}\sum_{\alpha=1}^{n_2} g_{ij}^{\tau}(\mathbf{Y}_\alpha) - \tau_{1,ij} + \tau_{2,ij})^2}{\zeta_{1,ij}/n_1 + \zeta_{2,ij}/n_2}\Bigg| = o_p(1), \tag{B.26}$$

where $g_{ij}^{\tau}$ is defined in (A.17). By Lemma B.2, we have

$$\mathbb{P}\Bigg(\max_{(i,j)\in S} \frac{(\frac{1}{n_1}\sum_{\alpha=1}^{n_1} g_{ij}^{\tau}(\boldsymbol{X}_\alpha) - \frac{1}{n_2}\sum_{\alpha=1}^{n_2} g_{ij}^{\tau}(\mathbf{Y}_\alpha) - \tau_{1,ij} + \tau_{2,ij})^2}{\zeta_{1,ij}/n_1 + \zeta_{2,ij}/n_2} > t^2\Bigg) \le C_1|S|(1 - \Phi(t)), \tag{B.27}$$

where $\Phi(t)$ is the cumulative distribution function of a standard normal distribution. By Lemma D.5, we have $1-\Phi(t) \le \exp(-t^2/2)/(t\sqrt{2\pi})$ for $t>0$. Considering $|S| = O(d^2)$ and (B.27), by setting $t^2 = 4\log d - 0.5\log(\log d)$, under Assumption **(A2)**, we obtain

$$\mathbb{P}\Bigg(\max_{(i,j)\in S} \frac{(\frac{1}{n_1}\sum_{\alpha=1}^{n_1} g_{ij}^{\tau}(\boldsymbol{X}_\alpha) - \frac{1}{n_2}\sum_{\alpha=1}^{n_2} g_{ij}^{\tau}(\mathbf{Y}_\alpha) - \tau_{1,ij} + \tau_{2,ij})^2}{\zeta_{1,ij}/n_1 + \zeta_{2,ij}/n_2} > 4\log d - \frac{1}{2}\log(\log d)\Bigg) = o(1). \tag{B.28}$$

Combining (B.26) and (B.28), we prove (A.28). □

## B.12. Proof of Lemma A.7

***Proof.*** By setting $S = \{(i,j) : 1 \le i < j \le d\}$, (A.7) implies that both of the following two events

$$\mathcal{E}_3 := \Big\{ \max_{(i,j)\in S} \big| n_1 \widehat{\sigma}^2_{\text{plug}}(\widehat{u}_{1,ij}) - m^2 \zeta_{1,ij} \big| < C\frac{\varepsilon_n}{\log q} \Big\},$$
$$\mathcal{E}_4 := \Big\{ \max_{(i,j)\in S} \big| n_2 \widehat{\sigma}^2_{\text{plug}}(\widehat{u}_{2,ij}) - m^2 \zeta_{2,ij} \big| < C\frac{\varepsilon_n}{\log q} \Big\},$$

happen with probability going to one as $n,d\to\infty$. Under $\mathcal{E}_3$ and $\mathcal{E}_4$, by $\zeta_{a,ij} \ge r_a > 0$ (see Assumption **(A2)**), we have

$$\max_{1\le i<j\le d} \left| \frac{n_a \widehat{\sigma}^2_{\text{plug}}(\widehat{\tau}_{a,ij})}{m^2 \zeta_{a,ij}} - 1 \right| \le C\frac{\varepsilon_n}{\log q} \quad \text{and} \quad \max_{1\le i<j\le d} \left| \frac{m^2 \zeta_{a,ij}}{n_a \widehat{\sigma}^2_{\text{plug}}(\widehat{\tau}_{a,ij})} - 1 \right| \le C\frac{\varepsilon_n}{\log q}.$$

Therefore, under $\mathcal{E}_3$ and $\mathcal{E}_4$, as $n,d\to\infty$, we have

$$\begin{aligned}&\left| \max_{(i,j)\in S} \frac{(\tau_{1,ij}-\tau_{2,ij})^2}{\widehat{\sigma}^2_{\text{plug}}(\widehat{\tau}_{1,ij}) + \widehat{\sigma}^2_{\text{plug}}(\widehat{\tau}_{2,ij})} - \max_{(i,j)\in S} \frac{(\tau_{1,ij}-\tau_{2,ij})^2}{m^2\zeta_{1,ij}/n_1 + m^2\zeta_{2,ij}/n_2} \right| \\ &\le \max_{(i,j)\in S} \left| \frac{(\tau_{1,ij}-\tau_{2,ij})^2}{\widehat{\sigma}^2_{\text{plug}}(\widehat{\tau}_{1,ij}) + \widehat{\sigma}^2_{\text{plug}}(\widehat{\tau}_{2,ij})} - \frac{(\tau_{1,ij}-\tau_{2,ij})^2}{m^2\zeta_{1,ij}/n_1 + m^2\zeta_{2,ij}/n_2} \right| \le C\frac{\varepsilon_n}{\log d}.\end{aligned}$$

Hence, as $n,d\to\infty$, we have

$$\left| \max_{(i,j)\in S} \frac{(\tau_{1,ij}-\tau_{2,ij})^2}{\widehat{\sigma}^2_{\text{plug}}(\widehat{\tau}_{1,ij}) + \widehat{\sigma}^2_{\text{plug}}(\widehat{u}_{2,ij})} - \max_{(i,j)\in S} \frac{(\tau_{1,ij}-\tau_{2,ij})^2}{m^2\zeta_{1,ij}/n_1 + m^2\zeta_{2,ij}/n_2} \right| = o_p(1). \tag{B.29}$$

Considering $(\mathbf{U}_1^{\tau}, \mathbf{U}_2^{\tau}) \in \mathbb{U}(4)$, by (B.29), we obtain (A.29) to complete the proof. □

## B.13. Proof of Lemma A.8

***Proof.*** We aim to prove (A.31) under Assumptions **(A1)** and **(A2)**. For Kendall's tau, we have

$$\widehat{T}_{ij} = \frac{\frac{1}{n_1}\sum_{\alpha=1}^{n_1} g_{ij}^{\tau}(\boldsymbol{X}_\alpha) - \frac{1}{n_2}\sum_{\alpha=1}^{n_2} g_{ij}^{\tau}(\mathbf{Y}_\alpha) - \tau_{1,ij} + \tau_{1,ij}}{\sqrt{\sigma_{1,\text{ps}}^2/4n_1 + \sigma_{2,\text{ps}}^2/4n_2}}, \tag{B.30}$$

$$T_{ij} = \frac{\frac{1}{n_1}\sum_{\alpha=1}^{n_1} g_{ij}^{\tau}(\boldsymbol{X}_\alpha) - \frac{1}{n_2}\sum_{\alpha=1}^{n_2} g_{ij}^{\tau}(\mathbf{Y}_\alpha) - \tau_{1,ij} + \tau_{1,ij}}{\zeta_{1,ij}/n_1 + \zeta_{2,ij}/n_2}. \tag{B.31}$$

Considering the definition of $M_n^{1,\text{ps}}$ in (A.30), by (A.23), we have

$$\left| M_n^{1,\text{ps}} - \max_{(i,j)\in S} \frac{(\frac{1}{n_1}\sum_{\alpha=1}^{n_1} g_{ij}^{\tau}(\boldsymbol{X}_\alpha) - \frac{1}{n_2}\sum_{\alpha=1}^{n_2} g_{ij}^{\tau}(\mathbf{Y}_\alpha) - \tau_{1,ij} + \tau_{1,ij})^2}{\sigma_{1,\text{ps}}^2/4n_1 + \sigma_{2,\text{ps}}^2/4n_2} \right| = o_p(1),$$

where $S = \{(i,j) : 1 \le i < j \le d\}$. Therefore, to obtain (A.31), it suffices to prove as $n$, $d \to \infty$, we have $\mathbb{P}\big(\max_{(i,j)\in S}(\widehat{T}_{ij})^2 > \widetilde{y}_d\big) = o(1)$, where $\widetilde{y}_d := 4\log d - \log(\log d)/2$.

We then aim to prove that as $n, d \to \infty$, we have

$$\left|\mathbb{P}\Big(\max_{(i,j)\in S}(\widehat{T}_{ij})^2 > \widetilde{y}_d\Big) - \mathbb{P}\Big(\max_{(i,j)\in S\setminus S_0}(\widehat{T}_{ij})^2 > \widetilde{y}_d\Big)\right| = o(1), \tag{B.32}$$

where $S_0$ is defined in (2.11). To prove (B.32), by

$$\left|\mathbb{P}\Big(\max_{(i,j)\in S}(\widehat{T}_{ij})^2 > \widetilde{y}_d\Big) - \mathbb{P}\Big(\max_{(i,j)\in S\setminus S_0}(\widehat{T}_{ij})^2 > \widetilde{y}_d\Big)\right| \le \mathbb{P}\Big(\max_{(i,j)\in S_0}(\widehat{T}_{ij})^2 > \widetilde{y}_d\Big),$$

we need to prove $\mathbb{P}\Big(\max_{(i,j)\in S_0}(\widehat{T}_{ij})^2 > \widetilde{y}_d\Big) = o(1)$ as $n, d \to \infty$. By Lemma B.3, we have

$$\mathbb{P}\Big(\max_{(i,j)\in S_0}(\widehat{T}_{ij})^2 > \widetilde{y}_d\Big) \le C|S_0|\big(1 - \Phi(C_1\sqrt{\widetilde{y}_d})\big), \tag{B.33}$$

where $\Phi(t)$ is the cumulative distribution function of a standard normal distribution. Lemma D.5 claims $1 - \Phi(t) \le \exp(-t^2/2)/(t\sqrt{2\pi})$ for $t > 0$. Considering $|S_0| = o(d^2)$ (see Assumption **(A1)**), by (B.33), we have $\mathbb{P}\Big(\max_{(i,j)\in S_0}(\widehat{T}_{ij})^2 > \widetilde{y}_d\Big) = o(1)$. Therefore, we prove that as $n, d \to \infty$, we have (B.32). To obtain (A.31), by (B.32), it suffices to prove that as $n$, $d \to \infty$, we have $\mathbb{P}\Big(\max_{(i,j)\in S\setminus S_0}(\widehat{T}_{ij})^2 > \widetilde{y}_d\Big) = o(1)$.

Considering the definition of $S_0$ in (2.11), by Lemma A.5, for any $(i,j) \in S \setminus S_0$, we have $|4\zeta_{a,ij} - \sigma_{a,\text{ps}}^2| \le C|\tau_{a,ij}| \le C(\log d)^{-1-\alpha_0}$. By $\zeta_{a,ij} \ge r_a > 0$ (see Assumption **(A2)**), for $(i,j) \in S \setminus S_0$, we then have

$$\big|\sigma_{a,\text{ps}}^2/4\zeta_{a,ij} - 1\big| \le C(\log d)^{-1-\alpha_0} \quad \text{and} \quad \big|4\zeta_{a,ij}/\sigma_{a,\text{ps}}^2 - 1\big| \le C(\log d)^{-1-\alpha_0}. \tag{B.34}$$

Therefore, by (B.34), for $(i,j) \in S \setminus S_0$, we have

$$\Big|\frac{(\widehat{T}_{ij})^2 - (T_{ij})^2}{(T_{ij})^2}\Big| \le \Big|\frac{\sigma_{1,\text{ps}}^2 - 4\zeta_{1,ij}}{\sigma_{1,\text{ps}}^2}\Big| + \Big|\frac{\sigma_{2,\text{ps}}^2 - 4\zeta_{2,ij}}{\sigma_{2,\text{ps}}^2}\Big| \le C(\log d)^{-1-\alpha_0}.$$

Therefore, for $(i,j) \in S \setminus S_0$, we have $|(\widehat{T}_{ij})^2 - (T_{ij})^2| \le C(T_{ij})^2(\log d)^{-1-\alpha_0}$. We then have

$$\begin{aligned} \Big| \max_{(i,j)\in S\setminus S_0} (\widehat{T}_{ij})^2 - \max_{(i,j)\in S\setminus S_0} (T_{ij})^2 \Big| &\le \max_{(i,j)\in S\setminus S_0} \Big| (\widehat{T}_{ij})^2 - (T_{ij})^2 \Big| \\ &\le C(\log d)^{-1-\alpha_0} \max_{(i,j)\in S\setminus S_0} (T_{ij})^2. \end{aligned} \tag{B.35}$$

Considering $\max_{(i,j)\in S\setminus S_0}(T_{ij})^2 = O_p(\log d)$, by (B.35), we have

$$\Big| \max_{(i,j)\in S\setminus S_0} (\widehat{T}_{ij})^2 - \max_{(i,j)\in S\setminus S_0} (T_{ij})^2 \Big| = o_p(1).$$

Therefore, it suffices to prove

$$\mathbb{P}\Big( \max_{(i,j)\in S\setminus S_0} (T_{ij})^2 > \widetilde{y}_d \Big) = o(1), \tag{B.36}$$

as $n, d \to \infty$. By Lemma B.2, we have

$$\mathbb{P}\Bigg( \max_{(i,j)\in S\setminus S_0} \frac{(\frac{1}{n_1}\sum_{\alpha=1}^{n_1} g_{ij}^{\tau}(\boldsymbol{X}_\alpha) - \frac{1}{n_2}\sum_{\alpha=1}^{n_2} g_{ij}^{\tau}(\boldsymbol{Y}_\alpha) - \tau_{1,ij} + \tau_{1,ij})^2}{\zeta_{1,ij}/n_1 + \zeta_{2,ij}/n_2} > t^2 \Bigg) \le C_1 |S\setminus S_0| (1 - \Phi(t)). \tag{B.37}$$

By Lemma D.5, we have $1 - \Phi(t) \le \exp(-t^2/2)/(t\sqrt{2\pi})$ for $t > 0$. Considering $|S \setminus S_0| = O(d^2)$ (see Assumption **(A1)**) and (B.37), by setting $t^2 = 4\log d - 0.5\log(\log d)$, as $n, d \to \infty$, we obtain

$$\mathbb{P}\Bigg( \max_{(i,j)\in S\setminus S_0} \frac{(\frac{1}{n_1}\sum_{\alpha=1}^{n_1} g_{ij}^{\tau}(\boldsymbol{X}_\alpha) - \frac{1}{n_2}\sum_{\alpha=1}^{n_2} g_{ij}^{\tau}(\boldsymbol{Y}_\alpha) - \tau_{1,ij} + \tau_{1,ij})^2}{\zeta_{1,ij}/n_1 + \zeta_{2,ij}/n_2} > \widetilde{y}_d \Bigg) = o(1).$$

Therefore, by (B.31), we prove (B.36). Hence, we complete the proof. ☐

## B.14. Proof of Lemma A.9

***Proof.*** We write $X_1 = Z_0 + Z_1$ with $Z_1 = \boldsymbol{\varrho}^T \boldsymbol{\Sigma}^{-1}(X_2, X_3, X_4)^T$ and $Z_0 = X_1 - Z_1$. We then obtain that $Z_0$ is independent of $(X_2, X_3, X_4)$, $Z_1 \sim N(0, \boldsymbol{\varrho}^T \boldsymbol{\Sigma}^{-1} \boldsymbol{\varrho})$ and $Z_0 \sim N(0, 1 - \boldsymbol{\varrho}^T \boldsymbol{\Sigma}^{-1} \boldsymbol{\varrho})$. We set $\sigma^2 := \boldsymbol{\varrho}^T \boldsymbol{\Sigma}^{-1} \boldsymbol{\varrho}$. Accordingly, we write

$$\begin{aligned} &\mathbb{E}\big[(\Phi(X_1) - 1/2)(\Phi(X_2) - 1/2)(\Phi(X_3) - 1/2)(\Phi(X_4) - 1/2)\big] \\ &\quad = \mathbb{E}\Big[\mathbb{E}\big[\Phi(X_1) - 1/2 \,|\, Z_1\big] \mathbb{E}\big[(\Phi(X_2) - 1/2)(\Phi(X_3) - 1/2)(\Phi(X_4) - 1/2) \,|\, Z_1\big]\Big]. \end{aligned} \tag{B.38}$$

We first focus on bounding $\big|\mathbb{E}\big[\Phi(X_1) - 1/2 \,|\, Z_1\big]\big|$. We use $\ell_3$ to denote it. Given $Z_1$, we have $Z_0 + Z_1 \,|\, Z_1 \sim N(Z_1, 1 - \sigma^2)$. Hence, we have $\ell_3 = \mathbb{E}\big[\Phi(\xi) - \mathbb{E}\Phi(\eta) \,|\, Z_1\big]$ with $\xi \sim N(Z_1, 1 - \sigma^2)$ and $\eta \sim N(0, 1)$. Because the equation holds for any joint distribution

of $N(Z_1, 1-\sigma^2)$ and $N(0,1)$, we have $\ell_3 = \big|\mathbb{E}\big[\Phi(Z_1+\sqrt{1-\sigma^2}Y_0) - \mathbb{E}\Phi(Y_0) \mid Z_1\big]\big|$, where $Y_0 \sim N(0,1)$ is independent of $Z_1$. Then we have

$$\ell_3 \leq \mathbb{E}\big[\big|\Phi(Z_1+\sqrt{1-\sigma^2}Y_0) - \Phi(Y_0)\big| \mid Z_1\big] \leq \frac{1}{\sqrt{2\pi}}\mathbb{E}\big[\big|Z_1+\sqrt{1-\sigma^2}Y_0 - Y_0\big| \mid Z_1\big].$$

Given $Z_1$, because $Z_1+\sqrt{1-\sigma^2}Y_0 - Y_0 \sim N(Z_1, (1-\sqrt{1-\sigma^2})^2)$, denoting $\sigma_2^2 = (1-\sqrt{1-\sigma^2})^2$, we then have

$$\Big|\mathbb{E}\big[\Phi(X_1) - \frac{1}{2} \mid Z_1\big]\Big| \leq \frac{1}{\sqrt{2\pi}}\Big[\sigma_2\sqrt{2/\pi}\exp(-Z_1^2/2\sigma_2^2) + Z_1\big[1-2\Phi(-Z_1/\sigma_2)\big]\Big]. \quad \text{(B.39)}$$

Moreover, it is easy to obtain $|\Phi(\cdot) - 1/2| \leq 1/2$. Therefore, we have

$$\mathbb{E}\Big[(\Phi(X_2)-1/2)(\Phi(X_3)-1/2)(\Phi(X_4)-1/2) \mid Z_1\Big] \leq 1/8. \quad \text{(B.40)}$$

Combining (B.38), (B.39) and (B.40), we then have

$$\begin{aligned}&\Big|\mathbb{E}\big[\Phi(X_1)-1/2)(\Phi(X_2)-1/2)(\Phi(X_3)-1/2)(\Phi(X_4)-1/2)\big]\Big| \\ &\quad \leq \frac{\sigma_2}{8\pi} + \frac{1}{8\sqrt{2\pi}}\mathbb{E}\big[Z_1\big(1-2\Phi(-Z_1/\sigma_2)\big)\big].\end{aligned} \quad \text{(B.41)}$$

Considering $\mathbb{E}[Z_1]^2 = \sigma^2$ and $\mathbb{E}[\Phi^2(Z_1/\sigma_2)] \leq 1$, we have

$$\mathbb{E}\big[Z_1\big(1-2\Phi(-Z_1/\sigma_2)\big)\big] = 2\mathbb{E}\big[Z_1\Phi(Z_1/\sigma_2)\big] \leq 2\mathbb{E}[Z_1^2]^{1/2}\mathbb{E}[\Phi^2(Z_1/\sigma_2)]^{1/2} \leq 2\sigma. \quad \text{(B.42)}$$

Noticing $\sigma_2 < \sigma$, combining (B.41) and (B.42), we obtain

$$\big|\mathbb{E}(\Phi(X_1)-1/2)(\Phi(X_2)-1/2)(\Phi(X_3)-1/2)(\Phi(X_4)-1/2)\big| \leq \left(\frac{1}{8\pi} + \frac{1}{4\sqrt{2\pi}}\right)\sigma.$$

This completes the proof. □

## B.15. Proof of Lemma A.10

***Proof.*** First, we prove that only when $r = 1$ and the set $\{\rho_1, \rho_2, a_1, \ldots, a_4\}$ attains the boundary, we have $C_r = 1$. When $C_r = 1$, we have

$$(\Phi(X_1)-1/2)(\Phi(X_2)-1/2) = a \cdot (\Phi(X_3)-1/2)(\Phi(X_4)-1/2),$$

for some constant $a$. This implies that

$$X_1 = \Phi^{-1}\left(\frac{a\cdot(\Phi(X_3)-1/2)(\Phi(X_4)-1/2)}{\Phi(X_2)-1/2} + 1/2\right).$$

We have $X_1 \sim N(0,1)$ if and only if we have

$$\frac{a\cdot(\Phi(X_3)-1/2)(\Phi(X_4)-1/2)}{\Phi(X_2)-1/2} \sim \text{Unif}(-1/2,1/2).$$

When $X_2 \neq \pm X_3$ and $X_2 \neq \pm X_4$, there is always a possibility such that $X_2$ is very close to zero and both $X_3$ and $X_4$ are away from zero. Under this condition, $(a\cdot(\Phi(X_3)-1/2)(\Phi(X_4)-1/2))/(\Phi(X_2)-1/2)$ is very large and outside $[-1/2,1/2]$. Therefore, $X_2$ must equal $\pm X_3$ or $\pm X_4$. Equivalently, $\{\rho_1,\rho_2,a_1,\dots,a_4\}$ attains the boundary. This completes the proof of the first part.

Secondly, it is obvious that there is a one-to-one map between $r$ and $C_r$. Accordingly, as long as $r<1$, $C_r<1$ only depends on $r$. □

# Supplement C: Proofs of Lemmas in Supplement B

In Supplement C, we present proofs of three lemmas introduced in Supplement B.

## C.1. Proof of Lemma B.1

***Proof.*** To prove (B.3) and (B.4), we need to analyze the following two terms:

$$A_1 := \mathbb{P}\Big(\Big|\frac{m^2(n_1-1)}{(n_1-m)^2}\big(\sum_{\alpha=1}^{n}(\widetilde{q}_{1\alpha,ij}-\widetilde{u}_{1,ij})^2-\sum_{\alpha=1}^{n_1}(h_{ij}(\boldsymbol{X}_\alpha)-\bar{h}_{1,ij})^2\big)\Big|\geq t\Big), \tag{C.1}$$

$$A_2 := \mathbb{P}\Big(\Big|\frac{m^2(n_1-1)}{(n_1-m)^2}\sum_{\alpha=1}^{n_1}(h_{ij}(\boldsymbol{X}_\alpha)-\bar{h}_{1,ij})^2-m^2\zeta_{1,ij}\Big|\geq t\Big). \tag{C.2}$$

We bound $A_1$ and $A_2$ separately. $A_2$ represents $\sum_{\alpha=1}^{n_1}(h_{ij}(\boldsymbol{X}_\alpha)-\bar{h}_{1,ij})^2/n_1$'s approximation error for $\zeta_{1,ij}$. Considering the definition of $\zeta_{a,ij}$ in Lemma D.3, by setting $\bar{h}_{1,ij}=\sum_{\alpha=1}^{n_1}h_{ij}(\boldsymbol{X}_\alpha)$, we have

$$\sum_{\alpha=1}^{n_1}(h_{ij}(\boldsymbol{X}_\alpha)-\bar{h}_{1,ij})^2=\sum_{\alpha=1}^{n_1}h_{ij}^2(\boldsymbol{X}_\alpha)-n_1(\bar{h}_{1,ij})^2 \quad \text{and} \quad \zeta_{1,ij}=\mathbb{E}[h_{ij}^2(\boldsymbol{X}_\alpha)]. \tag{C.3}$$

Thus, by using the triangle inequality, we obtain

$$\begin{aligned}&\mathbb{P}\Big(\Big|\frac{m^2(n_1-1)}{(n_1-m)^2}\sum_{\alpha=1}^{n_1}(h_{ij}(\boldsymbol{X}_\alpha)-\bar{h}_{1,ij})^2-m^2\zeta_{1,ij}\Big|\geq t/2\Big)\\ &\quad\leq \mathbb{P}\Big(\Big|\frac{1}{n_1}\sum_{\alpha=1}^{n_1}h_{ij}^2(\boldsymbol{X}_\alpha)-\mathbb{E}[h_{ij}^2(\boldsymbol{X}_\alpha)]\Big|\geq C_1t\Big)+\mathbb{P}\Big((\bar{h}_{1,ij})^2\geq C_2t\Big).\end{aligned}$$

By Lemma D.1, we use the Hoeffding's inequality to obtain

$$\mathbb{P}\Big(\Big|\sum_{\alpha=1}^{n_1} h_{ij}^2(\boldsymbol{X}_\alpha)/n_1 - \mathbb{E}[h_{ij}^2(\boldsymbol{X}_\alpha)]\Big| \geq C_1 t\Big) \leq C_3 \exp(-C_4 n_1 t^2). \tag{C.4}$$

Considering $h_{ij}(\boldsymbol{X}_\alpha)$ is bounded, by $\mathbb{E}[h_{ij}(\boldsymbol{X}_\alpha)] = 0$, we use Hoeffding inequality (Lemma D.1) to obtain

$$\mathbb{P}\Big((\bar{h}_{1,ij})^2 \geq C_2 t\Big) = \mathbb{P}\Big(|\bar{h}_{1,ij}| \geq \sqrt{C_2 t}\Big) \leq C_5 \exp(-C_6 n_1 t). \tag{C.5}$$

Hence, combing (C.4) and (C.5), we have

$$A_2 \leq C_3 \exp(-C_4 n_1 t^2) + C_5 \exp(-C_6 n_1 t). \tag{C.6}$$

By (C.6) and $\log q = O(n^{1/3-\epsilon})$ (Assumption **(A2)**), we have (B.4) by setting $\varepsilon_n = 1/(\log q)^{\kappa_0}$, where $\kappa_0 > 0$ is sufficiently small.

For $A_1$, we use (C.3) to rewrite it as

$$\mathbb{P}\Big(\Big|\frac{m^2(n_1-1)}{(n_1-m)^2}\big(\big(\sum_{\alpha=1}^{n_1}(\widetilde{q}_{1\alpha,ij})^2 - n_1(\widetilde{u}_{1,ij})^2\big) - \big(\sum_{\alpha=1}^{n_1}(h_{ij}(\boldsymbol{X}_\alpha))^2 - n_1(\bar{h}_{1,ij})^2\big)\big)\Big| \geq t/2\Big). \tag{C.7}$$

Similarly to $A_2$, we rearrange terms in $A_1$ and use the triangle inequality to obtain

$$\begin{aligned}
&\mathbb{P}\Big(\Big|\frac{m^2(n_1-1)}{(n_1-m)^2}\big(\big(\sum_{\alpha=1}^{n_1}(\widetilde{q}_{1\alpha,ij})^2 - n_1(\widetilde{u}_{1,ij})^2\big) - \big(\sum_{\alpha=1}^{n_1}(h_{ij}(\boldsymbol{X}_\alpha))^2 - n_1(\bar{h}_{1,ij})^2\big)\big)\Big| \geq t/2\Big) \\
&\leq \underbrace{\mathbb{P}\Big(\Big|\frac{m^2(n_1-1)}{(n_1-m)^2}\big(\sum_{\alpha=1}^{n_1}(\widetilde{q}_{1\alpha,ij})^2 - \sum_{\alpha=1}^{n_1}\big(h_{ij}(\boldsymbol{X}_\alpha)\big)^2\big)\Big| \geq t/4\Big)}_{I_1} \\
&\quad + \underbrace{\mathbb{P}\Big(\Big|\frac{m^2 n_1(n_1-1)}{(n_1-m)^2}\big((\widetilde{u}_{1,ij})^2 - (\bar{h}_{1,ij})^2\big)\Big| \geq t/4\Big)}_{I_2}.
\end{aligned} \tag{C.8}$$

Therefore, to bound $A_1$, we need to bound $I_1$ and $I_2$ separately. By

$$(\widetilde{u}_{1,ij})^2 - (\bar{h}_{1,ij})^2 = (\widetilde{u}_{1,ij} + \bar{h}_{1,ij})(\widetilde{u}_{1,ij} - \bar{h}_{1,ij}) \quad \text{and} \quad |\widetilde{u}_{1,ij} + \bar{h}_{1,ij}| \leq C,$$

we use the triangle inequality to bound $I_2$ and get

$$I_2 \leq \mathbb{P}\Big(\big|\widetilde{u}_{1,ij} - \bar{h}_{1,ij}\big| \geq Ct\Big) \leq \mathbb{P}\Big(\big|\widetilde{u}_{1,ij}\big| \geq Ct/2\Big) + \mathbb{P}\Big(\big|\bar{h}_{1,ij}\big| \geq Ct/2\Big).$$

We know that $\widetilde{u}_{1,ij}$ and $\bar{h}_{1,ij}$'s kernel functions are bounded. Noticing $\mathbb{E}[\widetilde{u}_{1,ij}] = \mathbb{E}[\bar{h}_{1,ij}] = 0$, we use the exponential inequalities (Lemmas D.1 and D.2) to obtain

$$I_2 \leq C_1 \exp(-C_2 n_1 t^2).$$

Considering $\log q = o(n^{1/3-\epsilon})$ in Assumption **(A2)**, by setting $\varepsilon_n = 1/(\log q)^{\kappa_0}$ with $\kappa_0 > 0$ sufficiently small, we have

$$q^2 \max_{1\le i,j\le q} \mathbb{P}\Big(\Big|\frac{m^2 n_1(n_1-1)}{(n_1-m)^2}\big((\widetilde{u}_{1,ij})^2 - (\bar{h}_{1,ij})^2\big)\Big| \ge C\frac{\varepsilon_n}{\log q}\Big) = o(1). \tag{C.9}$$

To prove (B.3), by (C.8) and (C.9), we only need to prove

$$q^2 \max_{1\le i,j\le q} \mathbb{P}\Big(\Big|\frac{m^2(n_1-1)}{(n_1-m)^2}\big(\sum_{\alpha=1}^{n_1}(\widetilde{q}_{1\alpha,ij})^2 - \sum_{\alpha=1}^{n_1}\big(h_{ij}(\boldsymbol{X}_\alpha)\big)^2\big)\Big| \ge C\frac{\varepsilon_n}{\log q}\Big) = o(1). \tag{C.10}$$

For this, we need to bound $I_1$. To bound $I_1$, we introduce the following decomposition of $\widetilde{q}_{1\alpha,ij}$. We set

$$\sum_{\substack{1\le \ell_1<\ldots<\ell_{m-1}\le n_1\\ \ell_k\neq\alpha,k=1,\ldots,m-1}}^{n_1} \Psi_{ij}(\boldsymbol{X}_\alpha,\boldsymbol{X}_{\ell_1},\ldots,\boldsymbol{X}_{\ell_{m-1}}) = Ah_{ij}(\boldsymbol{X}_\alpha) + BS_{1ij} + \Upsilon_{1ij}^{(\alpha)}, \tag{C.11}$$

with $A = \binom{n_1-1}{m_1-1} - \binom{n_1-2}{m-2}$, $B = \binom{n_1-2}{m-2}$ , $S_{1ij} := \sum_{\beta=1}^{n_1} h_{ij}(\boldsymbol{X}_\beta)$ and

$$\Upsilon_{1ij}^{(\alpha)} := \sum_{\substack{1\le \ell_1<\ldots<\ell_{m-1}\le n_1\\ \ell_j\neq\alpha,j=1,\ldots,m-1}}^{n_1} \Big(\Psi_{ij}(\boldsymbol{X}_\alpha,\boldsymbol{X}_{\ell_1}\ldots,\boldsymbol{X}_{\ell_m-1}) - \big(h_{ij}(\boldsymbol{X}_\alpha) + \sum_{k=1}^{m-1} h_{ij}(\boldsymbol{X}_{\ell_k})\big)\Big). \tag{C.12}$$

Furthermore, we set

$$V_{1ij}^2 := \sum_{\alpha=1}^{n_1}(h_{ij}(\boldsymbol{X}_\alpha))^2, \quad \Lambda_{1ij}^2 := \sum_{\alpha=1}^{n_1}(\Upsilon_{1ij}^{(\alpha)})^2 \quad \text{and} \quad D = \binom{n_1-1}{m-1}. \tag{C.13}$$

By $\widetilde{q}_{1\alpha,ij}$'s definition in (B.2), we have $\widetilde{q}_{1\alpha,ij} = (Ah_{ij}(\boldsymbol{X}_\alpha) + BS_{1ij} + \Upsilon_{1ij}^{(\alpha)})/D$. By setting $L_1$ as

$$L_1 := \Big|\frac{1}{n_1}\big(\frac{1}{D^2}\sum_{\alpha=1}^{n_1}\big(Ah_{ij}(\boldsymbol{X}_\alpha) + BS_{1ij} + \Upsilon_{1ij}^{(\alpha)}\big)^2 - \sum_{\alpha=1}^{n_1}\big(h_{ij}(\boldsymbol{X}_\alpha)\big)^2\big)\Big|, \tag{C.14}$$

bounding $I_1$ is equivalent to bounding $\mathbb{P}(L_1 \ge Ct)$. By expanding $\sum_{\alpha=1}^{n_1}(\widetilde{q}_{1\alpha,ij})^2$, we get

$$\sum_{\alpha=1}^{n_1}(\widetilde{q}_{1\alpha,ij})^2 = \frac{A^2V_{1ij}^2 + \Lambda_{1ij}^2 + (n_1B^2 + 2AB)(S_{1ij})^2 + 2A\sum_{\alpha=1}^{n_1} h_{ij}(\boldsymbol{X}_\alpha)\Upsilon_{1ij}^{(\alpha)} + 2BS_{1ij}\sum_{\alpha=1}^{n_1}\Upsilon_{1ij}^{(\alpha)}}{D^2}.$$

To bound $L_1$, we introduce

$$J_1 := \big|(A^2-D^2)V_{1ij}^2/(D^2n_1)\big|, \qquad J_2 := |\Lambda_{1ij}^2/(D^2n_1)|,$$
$$J_3 := |(n_1B^2 + 2AB)(S_{1ij})^2/(D^2n_1)|, \quad J_4 := 2AV_{1ij}\Lambda_{1ij}/(D^2n_1).$$

We use the Cauchy-Swartz inequality on $\sum_{\alpha=1}^{n_1} \Upsilon_{1ij}^{(\alpha)}$ to get

$$2BS_{1ij}\sum_{\alpha=1}^{n_1} \Upsilon_{1ij}^{(\alpha)}/(D^2 n_1) \leq 2B|S_{1ij}|\sqrt{n_1}\Lambda_{1ij}/(D^2 n_1).$$

This motivates us to set $J_5 := 2B|S_{1ij}|\sqrt{n_1}\Lambda_{1ij}/(D^2 n_1)$. By using the triangle inequality on $L_1$, we obtain $L_1 \leq J_1 + J_2 + J_3 + J_4 + J_5$. For $J_4$ and $J_5$, we use the Cauchy-Swartz inequality on $S_{1ij} = \sum_{\beta=1}^{n_1} h_{ij}(\boldsymbol{X}_\beta)$ to yield

$$J_4 + J_5 \leq \Big|\frac{2A}{n_1 D^2} V_{1ij}\Lambda_{1ij}\Big| + \Big|\frac{2Bn_1 V_{1ij}\Lambda_{1ij}}{n_1 D^2}\Big| \leq \Big|\frac{2A + 2Bn_1}{n_1 D^2} V_{1ij}\Lambda_{1ij}\Big|.$$

This motivates us to define $J_6 := |2(A + Bn_1)V_{1ij}\Lambda_{1ij}/n_1 D^2|$. Therefore, to bound $L_1$, we only need to bound $J_1$, $J_2$, $J_3$ and $J_6$ separately.

By the definitions of $A$ and $D$, we obtain

$$A = O(n_1^{m-1}), \qquad D = O(n_1^{m-1}) \qquad \text{and} \qquad D - A = \binom{n_1 - 2}{m - 2} = O(n_1^{m-2}).$$

Thus, for $J_1$, by the definition of $V_{1,ij}$ in (C.13), we use the Hoeffding inequality (Lemma D.1) to have

$$J_1 = \mathbb{P}\Big(\frac{1}{n_1^2} V_{1ij}^2 \geq Ct\Big) \leq C_1 \exp(-C_2 n_1^3 t^2 + C_3 n_1^2 t). \tag{C.15}$$

Similarly, for $J_3$, we use the Hoeffding inequality (Lemma D.1) to have

$$J_3 = \mathbb{P}\Big(\frac{(S_{1ij})^2}{n_1^2} \geq Ct\Big) \leq C_1 \exp(-C_2 n_1 t). \tag{C.16}$$

To bound $J_2$, by $\Lambda_{1ij}^2 := \sum_{\alpha=1}^{n_1} (\Upsilon_{1ij}^{(\alpha)})^2$, we have

$$\mathbb{P}\Big(\frac{\Lambda_{1ij}^2}{n_1^{2m-1}} \geq Ct\Big) = \mathbb{P}\Big(\sum_{\alpha=1}^{n_1} (\Upsilon_{1ij}^{(\alpha)})^2 \geq Cn_1^{2m-1} t\Big) \leq \sum_{\alpha=1}^{n_1} \mathbb{P}\Big((\Upsilon_{1ij}^{(\alpha)})^2 \geq Cn_1^{2m-2} t\Big).$$

By the definition of $\Upsilon_{1ij}^{(\alpha)}$ in (C.12), given $\boldsymbol{X}_\alpha$, we can treat

$$\Psi_{ij}(\boldsymbol{X}_\alpha, \boldsymbol{X}_{\ell_1} \dots, \boldsymbol{X}_{\ell_m - 1}) - \big(h_{ij}(\boldsymbol{X}_\alpha) + \sum_{k=1}^{m-1} h_{ij}(\boldsymbol{X}_{\ell_k})\big),$$

as a symmetric kernel function. Therefore, given $\boldsymbol{X}_\alpha$, $\Upsilon_{1ij}^{(\alpha)}/D$ is a U-statistic with a kernel function of zero mean and $m - 1$ order. Hence, by Lemma D.2, it follows that we have the following inequality:

$$\mathbb{P}\big((\Upsilon_{1ij}^{(\alpha)})^2 \geq Cn_1^{2m-2} t \big| \boldsymbol{X}_\alpha\big) = \mathbb{P}\big(|\Upsilon_{1ij}^{(\alpha)}| \geq Cn_1^{m-1}\sqrt{t} \big| \boldsymbol{X}_\alpha\big) \leq C_1 \exp(-C_2 n_1 t). \tag{C.17}$$

By taking expectation of (C.17), we have

$$\mathbb{P}\big((\Upsilon_{1ij}^{(\alpha)})^2 \geq Cn_1^2 t\big) \leq C_1 \exp(-C_2 n_1 t).$$

Hence, for $J_2$, we construct the following bound:

$$J_2 \leq \mathbb{P}\Big(\frac{\Lambda_{1ij}^2}{n_1^{2m-1}} \geq Ct\Big) \leq C_1 n_1 \exp(-C_2 n_1 t). \tag{C.18}$$

At last, by $A = O(n_1^{m-1})$, $B = O(n_1^{m-2})$ and $D = O(n_1^{m-1})$, bounding $J_6$ is equivalent to bounding $\mathbb{P}\big(V_{1ij}\Lambda_{1ij}/n_1^m \geq Ct\big)$. We then have

$$J_6 = \mathbb{P}\Big(\frac{V_{1ij}^2 \Lambda_{1ij}^2}{n_1^{2-2\kappa} n_1^{2m-2+2\kappa}} \geq Ct^2\Big) \leq \mathbb{P}\Big(\frac{V_{1ij}^2}{n_1^{2-2\kappa}} \geq Ct\Big) + \mathbb{P}\Big(\frac{\Lambda_{1ij}^2}{n_1^{2m-2+2\kappa}} \geq Ct\Big), \tag{C.19}$$

where $\kappa > 0$ will be given afterwards. Similarly to (C.15), we obtain

$$\mathbb{P}\Big(\frac{V_{1ij}^2}{n_1^{2-2\kappa}} \geq Ct\Big) \leq C_1 \exp(-C_2 n_1^{3-4\kappa} t^2 + C_3 n_1^{2-2\kappa} t). \tag{C.20}$$

For the other term in (C.19), we have

$$\mathbb{P}\Big(\frac{\Lambda_{1ij}^2}{n_1^{2m-2+2\kappa}} \geq Ct\Big) \leq \sum_{\alpha=1}^{n_1} \mathbb{P}\big((\Upsilon_{1ij}^{(\alpha)})^2 \geq Cn_1^{2m-3+2\kappa} t\big) = \sum_{\alpha=1}^{n_1} \mathbb{P}\Big(|\Upsilon_{1ij}^{(\alpha)}| \geq C\sqrt{n_1^{2m-3+2\kappa} t}\Big). \tag{C.21}$$

Similarly to (C.17), we have

$$\mathbb{P}\Big(|\Upsilon_{1ij}^{(\alpha)}| \geq C\sqrt{n_1^{2m-3+2\kappa} t}\Big) \leq C_4 \exp(-C_5 n_1^{2\kappa} t). \tag{C.22}$$

Therefore, Combining (C.19), (C.20), (C.21) and (C.22), we obtain

$$J_6 \leq C_1 \exp(-C_2 n_1^{3-4\kappa} t^2 + C_3 n_1^{2-2\kappa} t) + C_4 n_1 \exp(-C_5 n_1^{2\kappa} t). \tag{C.23}$$

Noticing $\log q = O(n^{1/3-\epsilon})$ (see Assumption **(A2)**), if we set $\kappa = 1/3$, combining (C.15), (C.16), (C.18) and (C.23), we have (C.10) by setting $\varepsilon_n = 1/(\log q)^{\kappa_0}$ with $\kappa_0 > 0$ sufficiently small. This completes the proof. □

## C.2. Proof of Lemma B.2

***Proof.*** By the definition of $\widehat{Z}_{\beta,ij}$ in (A.11), we obtain

$$\begin{aligned} &\frac{\big(\sum_{\alpha=1}^{n_1} g_{ij}(\boldsymbol{X}_\alpha)/n_1 - \sum_{\alpha=1}^{n_2} g_{ij}(\boldsymbol{Y}_\alpha)/n_2 - u_{1,ij} + u_{2,ij}\big)^2}{\zeta_{1,ij}/n_1 + \zeta_{2,ij}/n_2} \\ &\qquad = \frac{\big(\sum_{\beta=1}^{n_1+n_2} \widehat{Z}_{\beta,ij}\big)^2}{\sum_{\beta=1}^{n_1+n_2} \big(\widehat{Z}_{\beta,ij}\big)^2} \times \frac{\sum_{\beta=1}^{n_1+n_2} (\widehat{Z}_{\beta,ij})^2}{n_2^2 \zeta_{1,ij}/n_1 + n_2 \zeta_{2,ij}}, \end{aligned} \tag{C.24}$$

where $\zeta_{1,ij}$ and $\zeta_{2,ij}$ are variances of $g_{ij}(\boldsymbol{X}_\alpha)$ and $g_{ij}(\boldsymbol{X}_\alpha)$. By the Bernstein's inequality, for any $M > 0$, we have

$$\begin{aligned}
&\mathbb{P}\Big(\max_{(i,j)\in S}\Big|\frac{1}{n_1}\sum_{\beta=1}^{n_1}\big(\widehat{Z}_{\beta,ij}\big)^2 - \frac{n_2^2}{n_1^2}\zeta_{1,ij}\Big| \geq C\frac{\varepsilon_n}{\log q}\Big) = O(q^{-M}),\\
&\mathbb{P}\Big(\max_{(i,j)\in S}\Big|\frac{1}{n_2}\sum_{\beta=n_1+1}^{n_1+n_2}\big(\widehat{Z}_{\beta,ij}\big)^2 - \zeta_{2,ij}\Big| \geq C\frac{\varepsilon_n}{\log q}\Big) = O(q^{-M}),
\end{aligned} \tag{C.25}$$

by letting $\varepsilon_n$ converge to zero sufficiently slowly. Considering $\zeta_{a,ij} \geq r_a > 0$ (see Assumption **(A2)**), by (C.25), we have that the two events

$$\Big\{\Big|\frac{\sum_{\beta=1}^{n_1}\big(\widehat{Z}_{\beta,ij}\big)^2}{n_2^2\zeta_{1,ij}/n_1} - 1\Big| \leq C\frac{\varepsilon_n}{\log q}\Big\} \quad \text{and} \quad \Big\{\Big|\frac{\sum_{\beta=n_1+1}^{n_1+n_2}\big(\widehat{Z}_{\beta,ij}\big)^2}{n_2\zeta_{2,ij}} - 1\Big| \leq C\frac{\varepsilon_n}{\log q}\Big\},$$

happen with probability going to one, as $n, q \to \infty$. Under these two events, we have

$$\begin{aligned}
\Big|\frac{\sum_{\beta=1}^{n_1+n_2}(\widehat{Z}_{\beta,ij})^2}{n_2^2\zeta_{1,ij}/n_1 + n_2\zeta_{2,ij}} - 1\Big| &\leq \Big|\frac{\sum_{\beta=1}^{n_1}(\widehat{Z}_{\beta,ij})^2}{n_2^2\zeta_{1,ij}/n_1} - 1\Big| + \Big|\frac{\sum_{\beta=n_1+1}^{n_1+n_2}\big(\widehat{Z}_{\beta,ij}\big)^2}{n_2\zeta_{2,ij}} - 1\Big|\\
&\leq C\frac{\varepsilon_n}{\log q}.
\end{aligned} \tag{C.26}$$

By the self-normalized large deviation for independent variables in Jing et al. (2003), we have

$$\max_{(i,j)\in S}\mathbb{P}\Big(\big(\sum_{\beta=1}^{n_1+n_2}\widehat{Z}_{\beta,ij}\big)^2\Big/\sum_{\beta=1}^{n_1+n_2}\big(\widehat{Z}_{\beta,ij}\big)^2 \geq t^2\Big) \leq C(1-\Phi(t)), \tag{C.27}$$

uniformly for $t \in [0, O(n^{1/6-\epsilon})]$ for an arbitrary positive number $\epsilon$. Combining (C.24), (C.26) and (C.27), as $n, q \to \infty$, we have

$$\begin{aligned}
&\mathbb{P}\Big(\max_{(i,j)\in\Lambda}\frac{(\frac{1}{n_1}\sum_{\alpha=1}^{n_1}g_{ij}(\boldsymbol{X}_\alpha) - \frac{1}{n_2}\sum_{\alpha=1}^{n_2}g_{ij}(\boldsymbol{Y}_\alpha) - u_{1,ij} + u_{2,ij})^2}{\zeta_{1,ij}/n_1 + \zeta_{2,ij}/n_2} > t^2\Big)\\
&\quad\leq |\Lambda|\max_{(i,j)\in\Lambda}\mathbb{P}\Big(\big(\sum_{\beta=1}^{n_1+n_2}\widehat{Z}_{\beta,ij}\big)^2\Big/\sum_{\beta=1}^{n_1+n_2}\big(\widehat{Z}_{\beta,ij}\big)^2 \geq (1 + C\frac{\varepsilon_n}{\log q})t^2\Big) \leq C|\Lambda|(1-\Phi(t)).
\end{aligned}$$

This completes the proof of Lemma B.2. □

## C.3. Proof of Lemma B.3

***Proof.*** Lemma B.3 is the same as Lemma B.2 except that in Lemma B.3 we use $\sigma_{a,\mathrm{ps}}^2/4$ to replace $\zeta_{a,ij}$ used in Lemma B.2. This motivates us to write

$$\begin{aligned}&\frac{\big(\sum_{\alpha=1}^{n_1} g_{ij}^{\tau}(\boldsymbol{X}_\alpha)/n_1-\sum_{\alpha=1}^{n_2} g_{ij}^{\tau}(\boldsymbol{Y}_\alpha)/n_2-\tau_{1,ij}+\tau_{2,ij}\big)^2}{\sigma_{1,\mathrm{ps}}^2/4n_1+\sigma_{2,\mathrm{ps}}^2/4n_2}\\ &\quad=\frac{\big(\sum_{\alpha=1}^{n_1} g_{ij}^{\tau}(\boldsymbol{X}_\alpha)/n_1-\sum_{\alpha=1}^{n_2} g_{ij}^{\tau}(\boldsymbol{Y}_\alpha)/n_2-\tau_{1,ij}+\tau_{2,ij}\big)^2}{\zeta_{1,ij}/n_1+\zeta_{2,ij}/n_2}\cdot\frac{\zeta_{1,ij}/n_1+\zeta_{2,ij}/n_2}{\sigma_{1,\mathrm{ps}}^2/4n_1+\sigma_{2,\mathrm{ps}}^2/4n_2},\end{aligned}$$

where $g_{ij}^{\tau}$ is defined in (A.17) and $\zeta_{1,ij}$ and $\zeta_{2,ij}$ are variances of $g_{ij}^{\tau}(\boldsymbol{X}_\alpha)$ and $g_{ij}^{\tau}(\boldsymbol{Y}_\alpha)$. By the definitions of $g_{ij}^{\tau}$ and $\zeta_{a,ij}$, under Assumption **(A2)**, we have $0<r_a\le\zeta_{a,ij}\le 1$. We then have

$$0<(C_1)^2\le\frac{\sigma_{1,\mathrm{ps}}^2/4n_1+\sigma_{2,\mathrm{ps}}^2/4n_2}{1/n_1+1/n_2}\le\frac{\sigma_{1,\mathrm{ps}}^2/4n_1+\sigma_{2,\mathrm{ps}}^2/4n_2}{\zeta_{1,ij}/n_1+\zeta_{2,ij}/n_2}. \tag{C.28}$$

Therefore, by Lemma B.2, we have

$$\begin{aligned}&\mathbb{P}\bigg(\max_{(i,j)\in\Lambda}\frac{(\frac{1}{n_1}\sum_{\alpha=1}^{n_1} g_{ij}^{\tau}(\boldsymbol{X}_\alpha)-\frac{1}{n_2}\sum_{\alpha=1}^{n_2} g_{ij}^{\tau}(\boldsymbol{Y}_\alpha)-\tau_{1,ij}+\tau_{2,ij})^2}{\sigma_{1,\mathrm{ps}}^2/(4n_1)+\sigma_{2,\mathrm{ps}}^2/(4n_2)}>t^2\bigg)\\ &\quad\le\mathbb{P}\bigg(\max_{(i,j)\in\Lambda}\frac{(\frac{1}{n_1}\sum_{\alpha=1}^{n_1} g_{ij}^{\tau}(\boldsymbol{X}_\alpha)-\frac{1}{n_2}\sum_{\alpha=1}^{n_2} g_{ij}^{\tau}(\boldsymbol{Y}_\alpha)-\tau_{1,ij}+\tau_{2,ij})^2}{\zeta_{1,ij}/n_1+\zeta_{2,ij}/n_2}>(C_1)^2t^2\bigg)\\ &\quad\le C|\Lambda|(1-\Phi(C_1t)).\end{aligned}$$

This completes the proof of Lemma B.3. □

# Supplement D: Some Useful Technical Lemmas

In this appendix, we introduce the following technical lemmas that we use in our proofs.

**Lemma D.1.** (Hoeffding's inequality (Hoeffding, 1963)). $X_1,\ldots,X_n$ are independent random variables and $X_i$ takes its value in $[a_i,b_i]$. Letting $S=\sum_{i=1}^n(X_i-\mathbb{E}[X_i])$, we have

$$\mathbb{P}(S\ge x)\le\exp\Big(-\frac{2x^2}{\sum_{i=1}^n(b_i-a_i)^2}\Big),$$

for any positive $x\in\mathbb{R}$.

To present some useful results on U-statistics, we set $\boldsymbol{X}_1,\ldots,\boldsymbol{X}_n$ to be independent and identically distributed random vectors in $\mathbb{R}^d$. $\Phi(\mathbf{x}_1,\ldots,\mathbf{x}_m)$ is a kernel function of $m(\le n)$ vectors $\mathbf{x}_\gamma=(x_{\gamma1},\ldots,x_{\gamma d})^T$. Hence, the corresponding U-statistic is defined as $U:=\sum\Phi(\boldsymbol{X}_{i_1},\ldots,\boldsymbol{X}_{i_m})/n(n-1)\cdots(n-m+1)$, where the summation is taken over all $m$-tuples $i_1,\ldots,i_m$ of distinct positive integers not exceeding $n$.

**Lemma D.2.** (Exponential inequality for bounded U-statistics (Hoeffding, 1963)). If the kernel function $\Phi$ is bounded, i.e., $a \leq \Phi(\mathbf{x}_1, \cdots, \mathbf{x}_m) \leq b$, we have

$$\mathbb{P}(U - \mathbb{E}[U] \geq t) \leq \exp(-2kt^2/(b-a)^2),$$

where $k = \lfloor n/m \rfloor$.

We define $\boldsymbol{U} = (U_1, \ldots, U_K) \in \mathbb{R}^K$ as a $K$-dimensional U-statistic with the kernel function $\boldsymbol{\Phi}(\cdot) = \big(\Phi_1(\cdot), \ldots, \Phi_K(\cdot)\big)^T \in \mathbb{R}^K$, where $\Phi_k(\cdot)$ is defined as

$$\Phi_k(\cdot) : \underbrace{\mathbb{R}^d \times \cdots \times \mathbb{R}^d}_{m(k)} \to \mathbb{R} \qquad \text{for} \qquad 1 \leq k \leq K.$$

Letting $\phi_k := \mathbb{E}\big[\Phi_k(\boldsymbol{X}_1, \ldots, \boldsymbol{X}_{m(k)})\big]$, we define

$$\zeta_{k\ell} = \mathbb{E}\Big[\mathbb{E}\big[\Phi_k(\boldsymbol{X}_1, \ldots, \boldsymbol{X}_{m(k)}) - \phi_k | \boldsymbol{X}_1\big]\mathbb{E}\big[\Phi_\ell(\boldsymbol{X}_1, \ldots, \boldsymbol{X}_{m(\ell)}) - \phi_\ell | \boldsymbol{X}_1\big]\Big].$$

Hence, by setting $\boldsymbol{\Sigma}^{\mathbf{U}} = (m(k)m(\ell)\zeta_{k\ell}) \in \mathbb{R}^{K\times K}$, we introduce the central limit theorem for U-statistics.

**Lemma D.3.** (Central limit theorem for U-statistics (Hoeffding, 1948)). We assume $\Phi_k(\cdot)$ is a symmetric kernel for $1 \leq k \leq K$. If $\mathbb{E}\big[\Phi_k^2(\boldsymbol{X}_1, \ldots, \boldsymbol{X}_{m(k)})\big]$ exists for $k = 1, \ldots, K$, as $n \to \infty$, we have

$$\big(\sqrt{n}(U_1 - \phi_1), \ldots, \sqrt{n}(U_K - \phi_K)\big)^T \to N(\mathbf{0}, \boldsymbol{\Sigma}^{\mathbf{U}}),$$

where $N(\mathbf{0}, \boldsymbol{\Sigma}^{\mathbf{U}})$ is the distribution of a normal random vector with mean vector $\mathbf{0}$ and $K \times K$ covariance matrix $\boldsymbol{\Sigma}^{\mathbf{U}}$.

To get the central limit theorem, we need the Hoeffding decomposition. The Hoeffding decomposition divides a U-statistic into two pieces. One is a sum of i.i.d random variables. The other is a small residual term. For notational simplicity, for $\ell = 1, \ldots, n$, we set

$$g(\boldsymbol{X}_\ell) = \mathbb{E}\big[\Phi(\boldsymbol{X}_{\ell_1}, \ldots, \boldsymbol{X}_{\ell_m}) | \boldsymbol{X}_\ell\big] \qquad \text{and} \qquad h(\boldsymbol{X}_\ell) = g(\boldsymbol{X}_\ell) - \phi,$$

where $\phi = E[g(\boldsymbol{X}_\ell)]$ and $\ell \in \{\ell_1, \ldots \ell_m\}$.

**Lemma D.4.** (Hoeffding Decomposition (Hoeffding, 1948)). If $U$ is a U-statistic with a symmetric kernel function $\Phi(\mathbf{x}_1, \mathbf{x}_2, \ldots, \mathbf{x}_m)$, we can decompose $U$ as

$$U = \frac{m}{n}\sum_{\ell=1}^{n} h(\boldsymbol{X}_\ell) + \binom{n}{m}^{-1}\Delta_n,$$

where we set

$$\Delta_n = \sum_{1 \leq \ell_1 < \ell_2 < \ldots < \ell_m \leq n} \Big(\Phi(\boldsymbol{X}_{\ell_1}, \ldots, \boldsymbol{X}_{\ell_m}) - \phi - \sum_{k=1}^{m} h(\boldsymbol{X}_{\ell_k})\Big).$$

In Lemma D.4, $m\sum_{\ell=1}^{n} h(\boldsymbol{X}_\ell)/n$ is a sum of i.i.d random variables and $\binom{n}{m}^{-1}\Delta_n$ is a small residual term. Therefore, we can use $m\sum_{\ell=1}^{n} h(\boldsymbol{X}_\ell)/n$ to approximate $U$. Combining Lemmas D.3 and D.4, we know that $U$ and $m\sum_{\ell=1}^{n} h(\boldsymbol{X}_\ell)/n$ have the same limiting distribution as $n \to \infty$.

Next, we introduce a lemma on the tail probability of a normal distribution.

**Lemma D.5.** If $\xi$ follows a standard normal distribution, we have

$$\frac{t}{t^2+1}\cdot\frac{1}{\sqrt{2\pi}}e^{-t^2/2} \le \mathbb{P}(\xi > t) \le \frac{1}{t}\cdot\frac{1}{\sqrt{2\pi}}e^{-t^2/2},$$

for any $t > 0$.

# Supplement E: Introduction to the Meta-Elliptical Distribution

This section introduces the meta-elliptical distribution, which is essentially the elliptical copula family. It does not require the latent correlation matrix must be positive definite.

**Definition E.1.** Let $\boldsymbol{\mu} \in \mathbb{R}^d$ and $\boldsymbol{\Sigma} \in \mathbb{R}^{d\times d}$ with $\text{Rank}(\boldsymbol{\Sigma}) = q \le d$. A $d$-dimensional random vector $\boldsymbol{X}$ has an elliptical distribution, denoted by $\boldsymbol{X} \sim EC_d(\boldsymbol{\mu}, \boldsymbol{\Sigma}, \xi)$, if it has a stochastic representation: $\boldsymbol{X} = \boldsymbol{\mu} + \xi\mathbf{A}\mathbf{U}$, where $\mathbf{U}$ is a random vector uniformly distributed on the unit sphere in $\mathbb{R}^q$ and $\xi \ge 0$ is a scalar random variable independent of $\mathbf{U}$. $\mathbf{A} \in \mathbb{R}^{d\times q}$ is a deterministic matrix such that $\mathbf{A}\mathbf{A}^T = \boldsymbol{\Sigma}$.

To define the meta-elliptical distribution, we define a set of symmetric matrices:

$$\mathfrak{R}_d = \{\boldsymbol{\Sigma}^\ell \in \mathbb{R}^{d\times d} : (\boldsymbol{\Sigma}^\ell)^T = \boldsymbol{\Sigma}^\ell,\ \text{Diag}(\boldsymbol{\Sigma}^\ell) = \mathbf{I}_d,\ \boldsymbol{\Sigma}^\ell \succeq 0\}.$$

**Definition E.2.** A continuous random vector $\boldsymbol{X} = (X_1, \dots, X_d)^T \in \mathbb{R}^d$ follows a meta-elliptical distribution, denoted by $\boldsymbol{X} \sim ME_d(\boldsymbol{\Sigma}^\ell, \xi; f_1, \dots, f_d)$, if there exist univariate strictly increasing functions $f_1, \dots, f_d$ such that

$$(f_1(X_1), \dots, f_d(X_d))^T \sim EC_d(\mathbf{0}, \boldsymbol{\Sigma}^\ell, \xi),$$

where $\boldsymbol{\Sigma}^\ell = (\sigma_{ij}^\ell) \in \mathfrak{R}_d$. We call $\boldsymbol{\Sigma}^\ell$ the latent generalized correlation matrix.

The following theorem illustrates an important relationship between the population Kendall's tau coefficient matrix $\mathbf{U}^\tau$ and the latent generalized correlation matrix $\boldsymbol{\Sigma}^\ell$.

**Theorem E.3.** (Invariance property of Kendall's tau (Han and Liu, 2014)). If we let $\boldsymbol{X} = (X_1, \dots, X_d)^T \sim ME_d(\boldsymbol{\Sigma}^\ell, \xi; f_1, \dots, f_d)$, and denote $\tau_{ij}$ to be the population Kendall's tau between $X_i$ and $X_j$, we have $\sigma_{ij}^\ell = \sin(\tau_{ij}\pi/2)$.

Theorem E.3 shows that if $\boldsymbol{X}$ and $\boldsymbol{Y}$ follow the meta-elliptical distribution, testing (1.2) is equivalent to testing

$$\mathbf{H_0} : \boldsymbol{\Sigma}_1^\ell = \boldsymbol{\Sigma}_2^\ell \ \text{ v.s. }\ \mathbf{H_1} : \boldsymbol{\Sigma}_1^\ell \ne \boldsymbol{\Sigma}_2^\ell.$$

If $\boldsymbol{X}$ and $\boldsymbol{Y}$ follow multivariate Gaussian distributions, $\boldsymbol{\Sigma}_1^\ell$ and $\boldsymbol{\Sigma}_2^\ell$ are Pearson's correlation coefficient matrices of $\boldsymbol{X}$ and $\boldsymbol{Y}$.

**Table 4.** This table compares $\mathrm{T}_\alpha^{\mathrm{R}}$'s emprical sizes and powers with those of $\mathrm{T}_\alpha^{\mathrm{CLX}}$ and the three proposed tests under the Gaussian distribution in **Model 1** with $\mathbf{\Sigma} = \mathbf{\Sigma}^*$. Results show that the performances of $\mathrm{T}_\alpha^{\mathrm{CLX}}$ and $\mathrm{T}_\alpha^{\mathrm{R}}$ are similar.

| | | Empirical size | | | | Empirical power | | | |
|---|---|---|---|---|---|---|---|---|---|
| $n$ | $d$ | 50 | 100 | 300 | 500 | 50 | 100 | 300 | 500 |
| 500 | $\mathrm{T}_\alpha^{\tau,\mathrm{plug}}$ | 0.05 | 0.05 | 0.07 | 0.07 | 0.85 | 0.81 | 0.76 | 0.76 |
| | $\mathrm{T}_\alpha^{\tau,\mathrm{jack}}$ | 0.05 | 0.05 | 0.06 | 0.06 | 0.81 | 0.80 | 0.75 | 0.75 |
| | $\mathrm{T}_\alpha^{\tau,\mathrm{ps}}$ | 0.05 | 0.04 | 0.05 | 0.04 | 0.83 | 0.79 | 0.73 | 0.73 |
| | $\mathrm{T}_\alpha^{\mathrm{CLX}}$ | 0.04 | 0.04 | 0.05 | 0.04 | 0.83 | 0.78 | 0.74 | 0.72 |
| | $\mathrm{T}_\alpha^{\mathrm{R}}$ | 0.04 | 0.04 | 0.05 | 0.04 | 0.83 | 0.77 | 0.74 | 0.72 |

# Supplement F: More Simulation Results

This section provides another heuristic test (denoted by $\mathrm{T}_\alpha^{\mathrm{R}}$) for testing the equality of Pearson's correlation matrices. The performances of $\mathrm{T}_\alpha^{\mathrm{CLX}}$ and $\mathrm{T}_\alpha^{\mathrm{R}}$ are similar and both outperformed by the proposed rank-based methods.

We first explain the construction of $\mathrm{T}_\alpha^{\mathrm{R}}$. We set $\bar{\boldsymbol{X}} := (\bar{X}_1, \ldots, \bar{X}_d)^T$ and $\bar{\boldsymbol{Y}} := (\bar{Y}_1, \ldots, \bar{Y}_d)^T$, where $\bar{X}_j = (X_{1j}+X_{2j}+\ldots+X_{n_1 j})/n_1$ and $\bar{Y}_j = (Y_{1j}+Y_{2j}+\ldots+X_{n_2 j})/n_2$. We use $\widehat{s}_{1j}$ and $\widehat{s}_{2j}$ to denote the sample standard deviations of $X_j$ and $Y_j$ for $1 \leq j \leq d$. We then construct $\widetilde{\boldsymbol{X}}_i := (\widetilde{X}_{i1}, \ldots, \widetilde{X}_{id})^T$ and $\widetilde{\boldsymbol{Y}}_i := (\widetilde{Y}_{i1}, \ldots, \widetilde{Y}_{id})^T$, where $\widetilde{X}_{ij} = (X_{ij} - \bar{X}_j)/\widehat{s}_{1j}$ and $\widetilde{Y}_{ij} = (Y_{ij} - \bar{Y}_j)/\widehat{s}_{2j}$. $\mathrm{T}_\alpha^{\mathrm{R}}$ uses the same testing procedure as $\mathrm{T}_\alpha^{\mathrm{CLX}}$ in Cai et al. (2013) except that it is built on the rescaled samples $\{\widetilde{\boldsymbol{X}}_i\}_{1\leq i\leq n_1}$ and $\{\widetilde{\boldsymbol{Y}}_i\}_{1\leq i\leq n_2}$ rather than $\{\boldsymbol{X}_i\}_{1\leq i\leq n_1}$ and $\{\boldsymbol{Y}_i\}_{1\leq i\leq n_2}$. Table 4 shows that performances of $\mathrm{T}_\alpha^{\mathrm{R}}$ and $\mathrm{T}_\alpha^{\mathrm{CLX}}$ are similar. We note that ideas similar to $\mathrm{T}_\alpha^{\mathrm{R}}$ have also been considered in Cai and Zhang (2014) for correlation estimation.